\documentclass[a4paper,oneside,12pt]{article}
\usepackage{amsmath,amsfonts,amscd,amssymb,amsthm}
\usepackage[english]{babel}
\usepackage[utf8]{inputenc}
 \usepackage[T1]{fontenc}
\usepackage{graphicx}
\usepackage{abraces}
\usepackage{upgreek}
\usepackage{mathrsfs}   
\usepackage{dsfont}
\usepackage[hang,small,labelfont=bf,up,textfont=it,up]{caption} 
\usepackage{mathtools}
\usepackage{color}
\definecolor{Red}{cmyk}{0,1,1,0}

\definecolor{Blue}{cmyk}{1,1,0,0}

\newtheorem{theorem}{Theorem}
\newtheorem{definition}{Definition}

\newtheorem{corollary}{Corollary}
\newtheorem{lemma}{Lemma}
\newtheorem{proposition}{Proposition}

\newtheorem{remark}{Remark}
\usepackage[hidelinks,pdfpagelabels=true,unicode=true]{hyperref} 
\hypersetup
{
	pdfauthor={Leandro Cioletti and Roberto Vila},
	pdftitle={Graphical Representations for Ising and Potts 
	Models in General External Fields},
	pdfkeywords={Random Cluster Model, Potts Model, Ising Model, 
	Non-Translation Invariant},
	pdfsubject={Non-Translational Invariant Statistical 
	Mechanics models},  
	pdflang={en-US}
}
\title{Graphical Representations for Ising and Potts Models in 
General External Fields}
\author{Leandro Cioletti \qquad and \qquad  Roberto Vila}
\date{\today}
\usepackage{tocloft}

\newcommand{\horrule}[1]{\rule{\linewidth}{#1}}

\addto\captionsenglish
{

}
\begin{document}
\maketitle
%%%%%%%%%%%%%%%%  Abstract  %%%%%%%%%%%%%%%%%%%%%%%%%%%%%%%%%%%
\begin{abstract}
This work is concerned with the theory of graphical 
representation for the Ising and Potts models over general 
lattices with non-translation invariant external field. 
We explicitly describe in terms of the random-cluster 
representation the distribution function and, consequently, the 
expected value of a single spin for the Ising and $q$-state 
Potts models with general external fields. 
We also consider the Gibbs states for the Edwards-Sokal 
representation of the Potts model with non-translation invariant
magnetic field and prove a version of the FKG inequality for the
so called general random-cluster model (GRC model) with free and
wired boundary conditions in the non-translation invariant case. 

Adding the amenability hypothesis on the lattice, we obtain the 
uniqueness of the  infinite connected component and the almost 
sure quasilocality of the Gibbs measures for the GRC model with 
such general magnetic fields. As a final application of the 
theory developed, we show the uniqueness of the Gibbs measures 
for the ferromagnetic Ising model with a positive power-law 
decay magnetic field with small enough power, as conjectured in 
\cite{BCCP15}. 
\end{abstract}
%%%%%%%%%%%%%%%%%%%%%%%%%%%%%%%%%%%%%%%%%%%%%%%%%%%%%%%%%%%%%%%
\setcounter{tocdepth}{1}
\noindent\horrule{1.0pt}
\vspace*{-0.7cm}
\tableofcontents
\vspace*{0.3cm}
\noindent\horrule{1.0pt}
%%%%%%%%%%%%%%%%%%%%%%%%%%%%%%%%%%%%%%%%%%%%%%%%%%%%%%%%%%%%%%%
\section{Introduction}
Graphical representations are extremely useful tools for the 
study of phase transition in Equilibrium Statistical Mechanics. 
Fortuin and Kasteleyn \cite{FK72}, marked the beginning of four 
decades of intense activity that produced a rather complete 
theory for translation invariant systems. These representations 
were successfully employed to obtain non-per\-tur\-ba\-ti\-ve 
and deep results for Ising and Potts models on the hypercubic 
lattice using percolation-type methods, namely the discontinuity
of the magnetization at the phase transition point for the
one-dimensional Ising and Potts models with $1/r^2$ interactions
\cite{ACCN88}, the knowledge of the asymptotic behavior of the 
eigenvalues of the covariance matrix of the Potts model 
\cite{BC96}, the Aizenman-Higuchi Theorem on the Choquet 
decomposition of the two-dimensional Ising and Potts models 
\cite{Aiz80,CDIV14,CV12,GeoHig00,Hig81} and the proof that the 
self-dual point on the square lattice 
$p_{sd}(q)=\sqrt{q}/(1+\sqrt{q})$ is the critical point for 
percolation in the random-cluster model  ($q\geqslant 1$) 
\cite{beffara-duminil}, see also the review \cite{vila}.  
For a detailed introduction to the random-cluster model
we refer the reader to \cite{Dum15,GHM01,Grimmett2,Hag98}.

The relationship between graphical representations and phase 
transitions in Ising/Potts-type models is typically considered 
with respect to the random-cluster model (RC model) and in view 
of the Edwards-Sokal coupling \cite{ES}.
Most papers employing such representations use spin models with 
null or translation invariant magnetic field, whereas we shall 
analyze graphical representations of the Ising and Potts models 
under arbitrary and non-translation invariant external fields, 
which is a significantly more complicated task for several 
reasons: when general boundary conditions are considered,
the FKG property is harder to prove - as previously noticed
 by \cite{BBCK00}, this property does not even hold for certain 
boundary conditions. In the absence of the magnetic field, phase
transitions in the spin system can be directly detected by the 
random-cluster representation, but now this relationship is 
subtle since in some cases the phase transitions 
(in the percolation sense) in the random-cluster model 
does not correspond to a transition 
in the corresponding spin model. Such difficulties also 
appear in the analyses of Dobrushin-like states \cite{GG02}, 
large $q$ order-disorder at the transition temperature 
\cite{CK03} and the effect of ``weak boundary conditions'' 
in the $q$-state Potts model \cite{BKM02}. 

Here the absence of symmetry brings questions regarding the 
color(s) of the infinite connected component(s), which need not 
be addressed in the case of null magnetic field, for instance. 
Furthermore, non-translation invariance causes many technical 
issues when using basic results from the classical theory of 
spin models and Ergodic Theory. 
To avoid confusion, on this paper the terms phase transition and
critical inverse temperature shall hereby be solely employed to 
express changing in the number of the Gibbs measures when the 
temperature varies. 

This paper is motivated by some recent works on ferromagnetic 
Ising model in non-uniform external fields 
\cite{AGB07,BCCP15,BC10,JS99,NOZ99,Navarrete}. 
Here, we are interested in developing the theory of graphical 
representation for non-translation invariant models whilst 
aiming for the problem of classifying which are the positive 
magnetic fields such that the ferromagnetic Ising model on the 
square lattice passes through a first order phase transition,
in terms of its  power law decay exponent. The formal 
Hamiltonian of this model is given by 
\begin{equation}
\label{hamiltoniano-Ising-formal-introducao}
H(\sigma)
=
-J\sum_{ \{i,j\} } \sigma_i\sigma_j
-\sum_{i} h_{i}\sigma_{i},
\end{equation}
where the first sum ranges over the pairs of nearest neighbors.
In this model, if the magnetic field 
$\pmb{h}=(h_i: i\in\mathbb{Z}^d)$ 
satisfies $\liminf h_i>0$, it has been proved \cite{BC10} that 
for any positive temperature the set of the Gibbs measures is a 
singleton, therefore for essentially bounded-from-below
positive magnetic fields the conclusion is similar to the one 
obtained by Lee and Yang \cite{LY}. 
In \cite{BCCP15}, the authors considered a positive, decreasing 
magnetic field and employed the Isoperimetric inequality and a 
Peierls-type argument to show that if the magnetic field is 
given by $h_i=h^{*}/|i|^{\alpha}$, where $h^*$ is a positive 
constant, then the model presents first order phase transition 
in every dimension $d\geqslant 2$, for any fixed exponent 
$\alpha>1$. On the other hand, if $\alpha<1$, they proved by 
means of a contour expansion that the uniqueness of the
Gibbs measures holds at very low and by other methods at very 
high temperatures, and conjectured that the set of Gibbs 
measures at any positive temperature should be a singleton. 
The authors in \cite{BCCP15} justified why the extension of 
their results to any positive temperature is not obvious by 
resorting to most of the known techniques, but we prove as an 
application of the theory to be developed that the conjecture 
holds true. This is done by extending some results of the 
seminal work \cite{BBCK00} to the non-translation invariant 
setting.

The paper is organized in three parts: the first part presents 
the relevant background material, including notation and the 
basic definitions of the models to be treated in subsequent 
parts. The second part is comprised of the theory on general 
finite graphs with free boundary conditions, the main results of 
which are the extension of the Edwards-Sokal coupling for general 
external fields and the explicit computation (in terms of the RC 
model) of the distribution function of a single spin of the 
Ising model with general external field and its expected values. 
These results are also generalized to the $q$-state Potts model   
in general external fields. The third part is concerned with the 
Potts, Edwards-Sokal and General random-cluster models in the 
non-translation invariant external fields setting with general 
boundary conditions.
It is inspired by the reference \cite{BBCK00}, but extends their 
results to non-translation invariant magnetic fields - a task 
that was occasionally non-trivial. In some cases, their results 
were essentially proved for very general fields and our work was 
simply to point out the necessary technical modifications.
Fundamental results such as the FKG inequality required 
non-trivial adaptations and for this reason we presented
its detailed proof for both free and wired boundary conditions 
in the so called general random-cluster model (GRC model) with 
non-translation invariant external field. Even with null external
field, the random cluster measures lacked the key property of 
the quasilocality of the Gibbs measures, although it is possible 
to have the said property almost surely by assuming the 
uniqueness of the infinite connected component. 
For a null magnetic field on the hypercubic lattice, this fact 
was first proved in \cite{Pfister}, however the geometry of 
the graph in this type of question is very important because for 
some non-amenable graphs such as regular trees even almost sure 
quasilocality fails, see \cite{FF,Hag96}. 
For the random-cluster measures with translation invariant 
magnetic field, defined over amenable graphs, almost 
quasilocality was shown in \cite{BBCK00} for those measures 
having almost surely at most one infinite connected component. 
These results were recovered here for GRC models with 
non-translation invariant magnetic fields.
The proofs of both the uniqueness of the
infinite connected component and of the quasilocality of the 
Gibbs measures are given and new ideas are introduced to 
circumvent the lack of translation invariance.   

The conjecture stating the uniqueness of the Gibbs measures
for the Ising model with power-law-decay magnetic
field ($\alpha<1$) is proved in the last section of the third 
part. As a corollary of one of the main results of this part 
(Theorem \ref{teo-unicidade-gibbs-percolacao}), we have obtained 
a characterization of the critical
\footnote
{
in case $\beta_c(\pmb{J},\pmb{h}) = +\infty$ we mean 
that the set of the Gibbs measures is a singleton
for all $0<\beta<+\infty$.
}
inverse temperature $\beta_c(\pmb{J},\pmb{h})$ of the 
ferromagnetic Ising model given by 
\eqref{hamiltoniano-Ising-formal-introducao} where 
$h_i=h^{*}/|i|^{\alpha}$, with $\alpha>1$ on the hypercubic 
lattice. Few facts are known about this inverse critical 
temperature. For example, in the positive external field case of 
the two-dimensional model 
\eqref{hamiltoniano-Ising-formal-introducao}
with the coupling constant \linebreak $\pmb{J}\equiv 1$ and 
$\sum_{i\in\mathbb{Z}^2}h_i<\infty$,
it follows from the Onsager exact solution and a general result 
\cite{Georgii88} about summable perturbations of the Gibbs 
measures that $\beta_c(\pmb{J},\pmb{h})=\log(1+\sqrt{2})$.
From \cite{BC10} it follows that 
$\beta_c(\pmb{J},\pmb{h}) = +\infty$
as long as $\liminf h_i>0$ in any dimension. The last section 
contains the proof that $\beta_c(\pmb{J},\pmb{h})$ is also 
trivial, i.e., $\beta_c(\pmb{J},\pmb{h}) = +\infty$
when $h_i=h^{*}/|i|^{\alpha}$, with $\alpha<1$.
The most interesting cases are those where we do have phase 
transition and the magnetic field is given by 
$h_i=h^{*}/|i|^{\alpha}$, with $1<\alpha<2$
(not summable on entire lattice).
For such cases, to the best of our knowledge, the only known 
fact about this critical point is that 
$\log(1+\sqrt{2})\leqslant \beta_c(1,\pmb{h})$, which is derived 
from the correlation inequalities. It is not known whether the 
Lieb-Simon inequality \cite{Sim,Lieb}, the 
Aizenman-Barsky-Fernández Theorem \cite{ABF} and other 
characterizations of the critical point 
(for example, \cite{DumTas15} )can be extended for the case 
$h_i=h^{*}/|i|^{\alpha}$, with $1<\alpha<2$.
%%%%%%%%%%%%%%%%%%%%%%%%%%%%%%%%%%%%%%%%%%%%%%%%%%%%%%%%%%%%%%%
\part{Basic definitions and models}
%%%%%%%%%%%%%%%%%%%%%%%%%%%%%%%%%%%%%%%%%%%%%%%%%%%%%%%%%%%%%%%
%
\section{Background in graph theory}
We say that a graph $G=(V,E)$ is a countable graph if its vertex
set $V$ is countable. 
As usual, a {\bf path} $\gamma$ on $G$ is an alternating sequence
of vertices and edges 
$\gamma=(v_0,e_1,v_1,e_2,\ldots,e_n,v_n)$, such that 
$v_i\neq v_j$ for all $0\leqslant i,j\leqslant n-1$,
$v_n\in V\setminus\{v_1,v_2,\ldots,v_{n-1}\}$
and $e_j=\{v_{j-1},v_j\}$ for $1\leqslant j\leqslant n$.
In case $v_0=v_n$ we say that $\gamma$ 
is a closed path or a circuit. The vertices $v_0$ and $v_n$ 
of $\gamma$ are called initial and final vertices, respectively. 
We say that $x,y\in V$ are {\bf connected} if $x=y$ or there is 
a path $\gamma$ on $G$ so that $x=v_0$ and $y=v_n$, denoted 
$
x \longleftrightarrow_{\hspace*{-0.49cm}\scriptscriptstyle G}\ y
$, 
and whenever it is clear from context, we shall remove the 
subscript $G$ from the notation. The length of a path  
$\gamma=(v_0,e_1,v_1,e_2,\ldots,e_n,v_n)$ is defined as  
$\vert \gamma \vert=n$. 

A graph $G$ is said to be a {\bf connected graph} if any two 
vertices $i,j\in V$ are connected, otherwise we say that $G$ is
disconnected. The connected component of $x\in V$ is the vertex 
set $C_x\equiv \{y\in V: y\leftrightarrow x\}$. 
The {\bf distance} $d_G(x,y)$ between $x,y\in \mathbb{V}$ is 
defined by $d_G(x,y)=0$ if $x=y$; $d_G(x,y)=+\infty$ if 
$x\notin C_y$ and    
\(
d_G(x,y) =
\inf\{\vert \gamma \vert: 
\gamma\ \text{is a path connecting}\ x\ \text{to}\ y 
\}
\),
if $x\in C_y$.

A graph $H=(\tilde{V},\tilde{E})$ is a subgraph of $G=(V,E)$ 
if $\tilde{V}\subset V$ and $\tilde{E}\subset E$. 
A subgraph $H$ of $G$ is an {\bf induced subgraph}
if it has the same vertex set as $G$ and a random subgraph of 
$G$ is an induced subgraph such that the edges are chosen 
randomly.

Any infinite countable connected graph
$\mathbb{L}\equiv (\mathbb{V}, \mathbb{E})$
will be called a lattice and from now on a 
finite subgraph of $\mathbb{L}$ will be denoted by $G=(V,E)$.
The vertex set $V$ will sometimes be called {\it the volume}
$V$.

There are several definitions for the boundary of a vertex set 
$V$ contained in $\mathbb{L}$. 
In this work, the boundary is defined as follows. 
\begin{definition}[Boundary of $V$]
	\label{fronteira para vertices}
	The boundary of an arbitrary vertex set $V$ in
	$\mathbb{L}$ is defined by
	$
	\partial V
	\equiv
	\{i\in \mathbb{V}\setminus V: \text{d}_{\mathbb{L}}(i,V)=1\},
	$ 
	where 
	$d_{\mathbb{L}}$ is the distance on the lattice $\mathbb{L}$. 
	See Figure \ref{cfl}.
\end{definition}
\begin{center} 
	\begin{minipage}{\linewidth}
	\makebox[\linewidth]{%
	\includegraphics[scale=0.253,keepaspectratio=true]{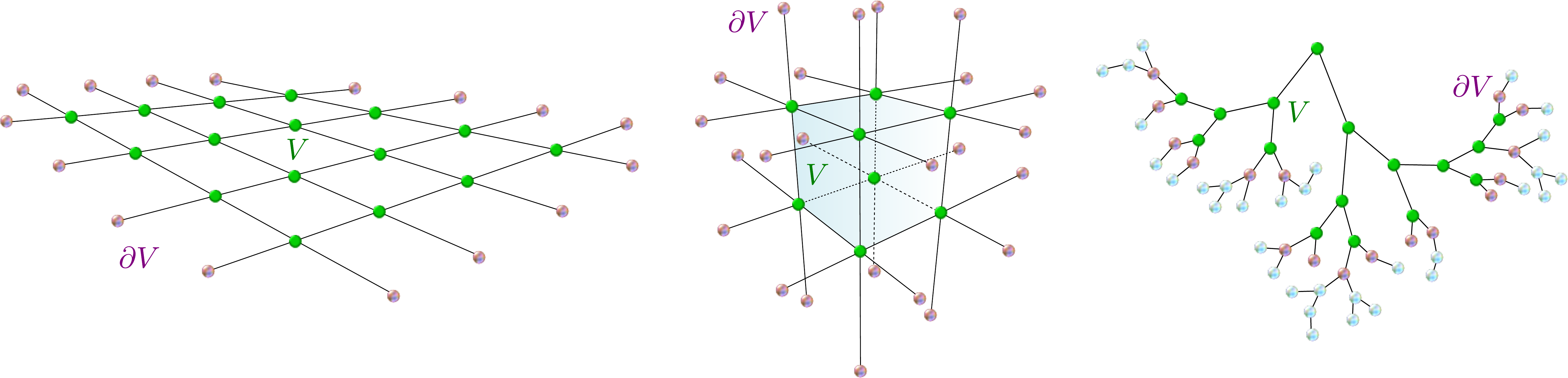}
	}
	\captionof{figure}
	{Examples of boundary of $V$ in three different lattices.
	The boundary of $V$ in each case is the vertex set colored 
	pink.}
	\label{cfl}
	\end{minipage}
\end{center}
\section{The Ising model on countable graphs}
\noindent
Let $\mathbb{L}=(\mathbb{V}, \mathbb{E})$ be an arbitrary 
lattice and $\Sigma$ the standard configuration space of the
Ising model, i.e.,
\[
\Sigma\equiv 
\big\{ 
\sigma
=
(\sigma_i: i\in \mathbb{V}): 
\sigma_i\in \{-1,+1\}, \ \forall i\in \mathbb{V}
\big\}
=
\{-1,+1\}^\mathbb{V}.
\]
This configuration space has a standard metric,
for which the distance between any pair of configurations 
$\sigma,\omega\in\Omega$ is given by
\[
d(\sigma,\omega)
=
\frac{1}{2^{R}}, 
\ \text{and}\
R
=
\inf 
\left\{
r>0:\ 
\begin{array}{c}
\sigma_i
=
\omega_i,\ \forall i\in B(o,r)
\ \text{and} \\ \ \exists\ j\in\partial B(o,r)\  \text{such that}\   
\sigma_j\neq \omega_j
\end{array}
\right\},
\]
where $B({o,r})$ is the open ball in $\mathbb{L}$ 
of center $o\in \mathbb{V}$ (fixed) and radius $r>0$.
Since the metric $d$ induces the product topology on $\Sigma$, 
it follows from Tychonoff's Theorem that $(\Sigma,d)$ is a 
compact metric space. As a measure space, we always consider 
$\Sigma$ endowed with the Borel $\sigma$-algebra 
$\mathscr{B}(\Sigma)$, which is generated by the open sets on 
$(\Sigma,d)$.  
\paragraph{The Hamiltonian of the Ising model}on a finite volume  
$V\subset \mathbb{L}$ with boundary condition $\mu\in\Sigma$ is 
given by
\begin{equation}\label{hamiltoniano-uni}
\mathscr{H}^{\mu,\mathrm{{Ising}}}_{\pmb{h},V}(\sigma)
\equiv 
- \sum_{
	    \substack{ i,j\in V  \\ \{i,j\}\in \mathbb{E} }
		}
J_{ij}\, \sigma_i \sigma_j
- 
\sum_{i\in V} h_i\, \sigma_i
-
\!\!\!\! 
\sum_{
\substack{ i\in V, \ j\in \partial V \\ \{i,j\}\in \mathbb{E} }
	  }
\!\!\!\!\!\!
J_{ij}\, \sigma_i \mu_j,
\end{equation}
where the  coupling constant
$\pmb{J}\equiv (J_{ij}:{\{i,j\}\in \mathbb{E}})\in 
[0,+\infty)^{\mathbb{E}}$ satisfies the regularity condition
$
\sum_{j\in \mathbb{V}} J_{ij} <+\infty, 
\ \forall i\in \mathbb{V}
$
and the magnetic field is
$
\pmb{h}\equiv (h_i:{i\in \mathbb{V}})\in \mathbb{R}^{\mathbb{V}}
$.
\paragraph{Gibbs measures.}
For any fixed finite volume $V$ and boundary condition $\mu$,
we define the (finite) set of configurations compatible with 
$\mu$ outside $V$ as being the set of configurations
$
\Sigma_{V}^{\mu}
\equiv
\{ 
\sigma\in \Sigma: \sigma_i = \mu_i \ \text{for} \ 
i\in \mathbb{V}\setminus V
\}.
$
The Gibbs measures of the Ising model on the finite volume $V$ 
with boundary condition $\mu$ at the inverse temperature 
$\beta>0$ is the probability measure
$
\lambda^{\mu}_{\beta, \pmb{h},V}:
\mathscr{B}(\Sigma)\to \mathbb{R}
$
given by
\begin{eqnarray*}
\lambda^{\mu}_{\beta, \pmb{h},V}(\{\sigma\})
=
\left\{
\begin{array}{lll} \label{i1}
{1\over \mathscr{Z}^{\mu,\mathrm{{Ising}}}_{\beta, \pmb{h},V}}
\exp
\big(-\beta \mathscr{H}_{\pmb{h},V}^{\mu,\mathrm{{Ising}}}
(\sigma)
\big), 
& \text{if} \ 
\sigma\in \Sigma_ V^{\mu}
\\\\
0, \ &\text{otherwise} 
\end{array}
\right. 
\end{eqnarray*}
where $\mathscr{Z}_{\beta, \pmb{h},V}^{\mu,\mathrm{{Ising}}}$ 
is a normalizing constant called the {\bf partition function} 
given by 
\[
\mathscr{Z}_{\beta, \pmb{h},V}^{\mu,\mathrm{{Ising}}}
=
\sum_{\sigma\in\Sigma_ V^{\mu}}
\exp
\big(
-\beta \mathscr{H}_{\pmb{h},V}^{
	\mu,\mathrm{{Ising}}}(\sigma)
\big).
\]
We denote by 
$\mathscr{G}^{\mathrm{{Ising}}}_{\beta}(\pmb{J},\pmb{h})$ 
the set of infinite-volume Gibbs measures which is given by the 
closure of the convex hull of the set of all the weak limits 
$
\lim_{V_n \uparrow \mathbb{V}}
\lambda^{\mu}_{\beta, \pmb{h},V},
$
where $V_n\subset V_{n+1}$ and $\mu$ runs over all possible 
sequences of boundary conditions.
\paragraph{The Ising model with free boundary condition.}
The Gibbs measure of the Ising model on a finite subgraph 
$G\subset \mathbb{L}$ with free boundary condition is given by
\begin{eqnarray*}
\lambda_{\beta, \pmb{h},V}(\{\sigma\})
=
{1\over \mathscr{Z}^{\mathrm{{Ising}}}_{\beta, \pmb{h},V}}
\exp
\big(
-\beta \mathscr{H}_{\pmb{h},V}^{\mathrm{{Ising}}}
(\sigma)
\big),
\end{eqnarray*}
where $\mathscr{Z}_{\beta,\pmb{h},V}^{\mathrm{ Ising}}$
is the partition function and the Hamiltonian is given by
\[
\mathscr{H}_{\pmb{h},V}^{\mathrm{ Ising}}
=
- \sum_{
	\{i,j\}\in E
}
J_{ij}\, \sigma_i \sigma_j
- 
\sum_{i\in V} h_i\, \sigma_i.
\] 
The expected value of 
a random variable $f:\Sigma\to\mathbb{R}$, 
with respect to $\lambda^{\mu}_{\beta, \pmb{h},V}$ is given by
\[
\lambda^{\mu}_{\beta, \pmb{h},V}(f)
\equiv 
\sum_{\sigma\in\Sigma^{\mu}_{V}} 
f(\sigma) 
\lambda^{\mu}_{\beta, \pmb{h},V}(\{\sigma\}).
\]
\section{The Potts model with inhomogeneous magnetic field}
\label{sec-Potts-com-campo-cc-livre}
\noindent
Let $q\in \mathbb{Z}^{+}$ be a fixed positive integer. 
The state space of the $q$-state Potts model on the lattice 
$\mathbb{L}$ is defined as
\[
\Sigma_{q}\equiv 
\big\{ 
\hat{\sigma}=(\hat{\sigma}_i: i\in \mathbb{V}): 
\hat{\sigma}_i\in \{1,2, \ldots, q\}, \ 
\forall i\in \mathbb{V}
\big\}
=
\{1,2, \ldots, q\}^\mathbb{V}.
\]
To define a $q$-state Potts model with magnetic field, we fix a
family of coupling constants 
$
\pmb{J}\equiv \left(J_{ij}:{\{i,j\}\in \mathbb{E}}\right)\in 
[0,\infty)^{\mathbb{E}}
$ 
and magnetic fields
$
\pmb{\hat{h}}\equiv (h_{i,p}:i\in \mathbb{V};
p=1,\ldots,q)\in \mathbb{R}^{\mathbb{V}}
\times 
\cdots\times \mathbb{R}^{\mathbb{V}}
$.
The Hamiltonian on a finite volume $G$ with
boundary condition $\hat{\mu}\in\Sigma_{q}$ is given by 
\begin{eqnarray}\label{Potts-Geral}
\mathscr{H}^{\hat{\mu},\mathrm{{Potts}}}_{\pmb{\hat{h}},q,V}
(\hat{\sigma})
\equiv 
- \sum_{
	\substack{ i,j\in V \\ \{i,j\}\in \mathbb{E} }
		}
J_{ij}\delta_{\hat{\sigma}_i, \hat{\sigma}_j}
- 
\sum_{p=1}^q
\sum_{i\in V} {h_{i,p}\over q}\delta_{\hat{\sigma}_{i,p}}
- 
\sum_{
	\substack{ 
		i\in V, \ j\in \partial V \\ \{i,j\}\in \mathbb{E} 
			}
	}
J_{ij}
\delta_{\hat{\sigma}_i, \hat{\mu}_j},
\end{eqnarray}
where $\delta_{\hat{\sigma}_i,\hat{\sigma}_j}$ is the Kronecker 
delta function.

The Gibbs measure of Potts model on a finite volume $G$ with
boundary condition $\hat{\mu}$ is defined analogously to the 
Ising model. We consider the set of all compatible 
configurations with the boundary condition $\hat{\mu}$, i.e.,
$
\Sigma_{q,V}^{\hat{\mu}}
\equiv
\{ 
\hat{\sigma}\in \Sigma_q: \hat{\sigma}_i = \hat{\mu}_i \ 
\text{for} \ 
i\in \mathbb{V}\setminus V
\}
$
and define the Gibbs measure on the volume $G$  with boundary 
condition $\hat{\mu}$ as the probability measure   
$\pi^{\hat{\mu}}_{\beta, \pmb{\hat{h}},q, V}$ on 
$\big(\Sigma_q, \mathscr{B}(\Sigma_q)\big)$ 
such that 
\[
\pi^{\hat{\mu}}_{\beta, \pmb{\hat{h}},q,V}(\hat{\sigma})
=
\left\{
\begin{array}{lll}
{
1
\over 
\mathscr{Z}^{\hat{\mu},\mathrm{{Potts}}}_{\beta,\pmb{\hat{h}},q,V}
}
\exp
\big(
	-\beta 
	\mathscr{H}_{\pmb{\hat{h}},q,V}^{\hat{\mu},\mathrm{{Potts}}}
	(\hat{\sigma})
\big),
& \text{if} \ 
\hat{\sigma}\in \Sigma_{q, V}^{\hat{\mu}}
\\
\\
0, \ &\text{otherwise}
\end{array}
\right.
\]
where 
$
\mathscr{Z}_{\beta,\pmb{\hat{h}},q,V}^{\hat{\mu},\mathrm{{Potts}}}
$ 
is the partition function. The free boundary condition case can 
be treated similarly to the previous section.
The expected value of a random variable $f:\Sigma_{q}\to\mathbb{R}$ 
in this model is denoted by 
$
\pi^{\hat{\mu}}_{\beta, \pmb{\hat{h}},q,V}(f)
$.
The set of infinite-volume Gibbs measures is defined 
similarly to the previous section and denoted by 
$\mathscr{G}^{\mathrm{{Potts}}}_{\beta}(\pmb{J},\pmb{\hat{h}})$.
\begin{remark}
In general, we use $\pmb{\hat{h}}$ 
to denote the magnetic field. 
In the special case where $q=2$ and the 
magnetic field satisfies $h_{i,1}=-h_{i,2}\equiv h_i$ 
we drop the hat from notation $\pmb{\hat{h}}$ and
write the Hamiltonian, 
for example in the free boundary condition case,  
as follows
	\begin{eqnarray}
	\label{rescrever-Potts}
	\mathscr{H}^{\mathrm{{Potts}}}_{\pmb{h},2,V}(\hat{\sigma})
	\equiv 
	- \sum_{
		\{i,j\}\in E
	}
	J_{ij}\delta_{\hat{\sigma}_i, \hat{\sigma}_j}
	- 
	\sum_{i\in V} {h_i\over 2}
	(
	\delta_{\hat{\sigma}_{i,1}}
	-
	\delta_{\hat{\sigma}_{i,2}}
	).
	\end{eqnarray}
\end{remark}
\begin{proposition}\label{Medeq}
	Fix a finite graph $G=(V,E)$
	and assume that the magnetic field 
	of the $2$-state Potts model satisfies 
	$h_{i,1}=-h_{i,2}\equiv h_i$ for all $i\in V$.  
	If $\hat{\sigma}\in \{1,2\}^V$ denotes the configuration 
	obtained from $\sigma\in\{-1,1\}^V$ using the spins 
	identification
	$-1\leftrightarrow 2$ and $1\leftrightarrow 1$, then
	we have for any $\beta>0$ that 
	\[
	\lambda_{{\beta},\pmb{h},V}^{\mu}(\{\sigma\})
	=
	\pi_{{2\beta},\pmb{h},2,V}^{\hat{\mu}}(\{\hat{\sigma}\}).
	\]	
\end{proposition}
\section{The random-cluster model with external field}
\label{secao-def-modelo-RC}
This section is devoted to the $q=2$ 
inhomogeneous random-cluster models  
on a finite graph $G=(V,E)$. 
The general random-cluster model in external field will be
introduced in the Section \ref{secao-potts-geral},
more precisely by the expression \eqref{def-rcm-q-geral}.  
The state space over which these models are defined 
is the cartesian product  $\{0,1\}^{E}$.
A generic element of this space will be denoted 
by $\omega$ and called an edge configuration.
We say that an edge $e$ is open in the configuration
$\omega$ if $\omega_e=1$, and we otherwise  
say $e$ is closed.  
Given $\omega\in \{0,1\}^{E}$, its set of open edges
is denoted by 
$\eta(\omega)=\{e\in E:  \omega_e=1\}$.
We say that a path $\gamma:=(v_0,e_1,v_1,e_2,\ldots,e_n,v_n)$
on the graph $G$ is an open path on $\omega$  
if all of its edges belong to $\eta(\omega)$,
i.e., $\omega_{e_i}=1, \ \forall i=1,\ldots,n$.  

Two distinct vertices $x,y\in V$ are said to be connected in 
$\omega$ if there exists an open
path $\gamma:=(v_0,e_1,v_1,e_2,\ldots,e_n,v_n)$ 
on this edge configuration such that $v_0=x$ and $v_n=y$.
If $x,y\in V$ are connected on $\omega$, 
we write $x \leftrightarrow y$.
A subgraph $H$ of $G$ is connected on $\omega$ if any pair 
of vertices of $H$ can be connected through a open path entirely 
contained in $H$. The open connected component of a vertex 
$x\in V$ is defined by 
$
C_x(\omega) 
\equiv 
\{y\in V: x\leftrightarrow y \ \text{in}\ \omega \}
\cup\{x\}.
$ 
The set $C_x(\omega)$ is called the open connected component of 
$x$ in the configuration $\omega$. 

To define the probability measure of the random-cluster model 
with external field we fix two families 
$
\pmb{p}
\equiv
(p_{ij}\in[0,1]:\{i,j\}\in E)\in [0,1]^{E}
$
and $\pmb{h}\equiv (h_i:{i\in V})\in \mathbb{R}^{V}$.
For convenience we will assume that the 
family $\pmb{p}$ is given by a family of coupling constants
$\pmb{J}=(J_{ij}\in [0,+\infty]: \{i,j\}\in E )$ 
and the inverse temperature $\beta>0$ so that 
$p_{ij}= 1-\exp(-2\beta J_{ij})$. Following \cite{CMR98},
the probability measure of the random-cluster model with external
field $\pmb{h}$ on the finite volume $G$ is defined for each
$\omega\in \{0,1\}^{E}$ by  
\begin{equation*}
\phi_{\pmb{p},\pmb{h},G}(\omega)
=
{
	1
	\over
	\mathscr{Z}^{\mathrm{{RC}}}_{\pmb{p},\pmb{h},G}
}
B_{\pmb{J}}(\omega)
\displaystyle
\prod_{\alpha =1}^{k({\omega,G})} 
2\cosh \big(\pmb{h}(K_\alpha(\omega))\big),
\end{equation*}
where $\pmb{h}(K_\alpha(\omega))\equiv \beta 
{\sum_{i\in K_\alpha(\omega)}} h_i$, with the sets 
$K_1(\omega),\ldots, K_{k({\omega,G})}(\omega)$ 
being composed by the connected components
of $(V,\eta(\omega))$,
$B_{\pmb{J}}(\omega)$ representing the Bernoulli factors
\begin{eqnarray}\label{fator de bernoulli}
\displaystyle{
	B_{\pmb{J}} (\omega)
	\equiv
	\prod_{\{i,j\}: \omega_{ij}=1}
	p_{ij}
	\prod_{\{i,j\}: \omega_{ij}=0}
	(1-p_{ij}) 
			}
\end{eqnarray}
and $\mathscr{Z}^{\mathrm{{RC}}}_{\pmb{p},\pmb{h},G}$ 
being the partition function
\[
\mathscr{Z}^{\mathrm{{RC}}}_{\pmb{p},\pmb{h},G}
=
\sum_{\omega\in \{0,1\}^{E}}
B_{\pmb{J}}(\omega)
\prod_{\alpha =1}^{k({\omega,G})} 
2\cosh \big(\pmb{h}(K_\alpha(\omega))\big).
\]
From now on, in order to ease the notation, we shall omit the
$\omega$-dependence from the components
$K_1(\omega),\ldots, K_{k({\omega,G})}(\omega)$ 
and simply write $K_1,\ldots, K_{k({\omega,G})}$.
%
%%%%%%%%%%%%%%%%%%%%%%%%%%%%%%%%%%%%%%%%%%%%%%%%%%%%%%%%%%%%%%%
\part{Free boundary conditions}
%%%%%%%%%%%%%%%%%%%%%%%%%%%%%%%%%%%%%%%%%%%%%%%%%%%%%%%%%%%%%%%
\section[The Edwards-Sokal coupling]{The Edwards-Sokal Coupling}
In this section we present the Edwards-Sokal model on a finite 
graph $G=(V,E)$.
The configuration space of this model is given by the 
cartesian product $\{-1,+1\}^V\times\{0,1\}^E$. 
A pair of configurations $\sigma\in \{-1,+1\}^V$ and 
$\omega\in\{0,1\}^{E}$ are deemed consistent if
$
\omega_{ij}=1 
\Rightarrow  \sigma_i=\sigma_j, \
\forall \ \{i,j\}\in E.
$
The indicator function of the consistency of a pair 
$(\sigma,\omega)\in \{-1,+1\}^V \times\{0,1\}^{E}$,
is denoted by 
\[
\Delta(\sigma, \omega)\equiv 
\mathds{1}_{\{(\xi, \eta)\in \{-1,+1\}^V \times\{0,1\}^{E} 
	:\ \text{if}\ \eta_{ij}=1
	\ \text{then}\  \xi_i=\xi_j \}}(\sigma, \omega)
.
\]
\vspace*{0.3cm}
\begin{center} 
	\begin{minipage}{\linewidth}
	\makebox[\linewidth]{%
	\includegraphics[scale=0.38,keepaspectratio=true]{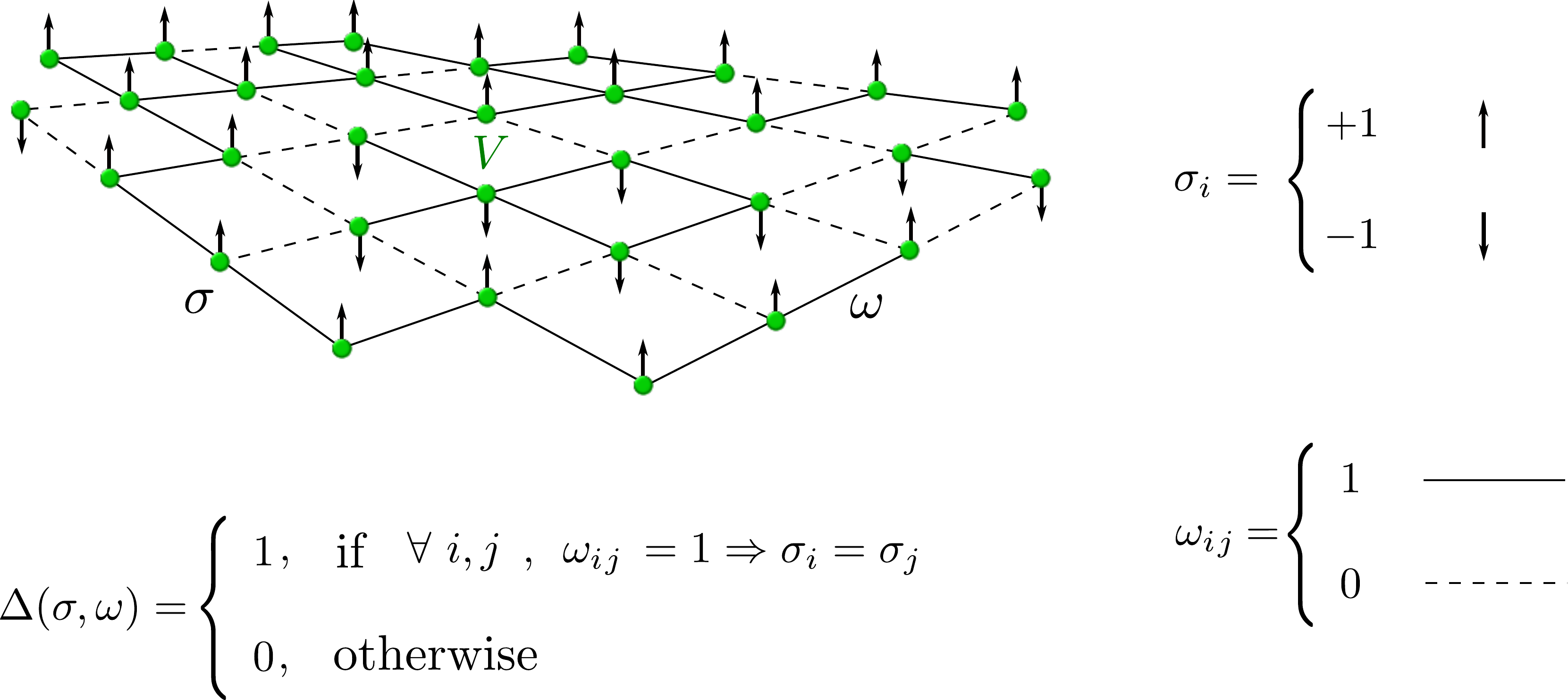}}
	\captionof{figure}{An example of spin-edge compatible 
	configuration.}
	\label{f2}
	\end{minipage}
\end{center}
\vspace*{0.3cm}
Similarly to the previous section, we fix $\beta>0$, coupling 
constants $\pmb{J}$ and magnetic field $\pmb{h}$, and put 
$
p_{ij}\equiv 1-\exp(-2\beta J_{ij}).
$
In the Edwards-Sokal model, the probability of a configuration 
$(\sigma, \omega)$ on a finite volume $G=(V,E)\subset \mathbb{L}$ 
is defined by
\begin{eqnarray*}
{\nu}_{\pmb{p},\pmb{h},G}(\sigma, \omega)
\equiv  
{1 \over \mathscr{Z}^{\mathrm{ ES}}_{\pmb{p},\pmb{h},G}}
\displaystyle{
	B_{\pmb{J}}(\omega)
	\Delta(\sigma, \omega)
	\times
	\exp \big( \beta 
	\sum_{i\in V} h_i \left(\delta_{\sigma_i,1}
	-\delta_{\sigma_i,-1} \right)
	\big),
}
\end{eqnarray*}
where $B_{\pmb{J}}$ represents the Bernoulli factors
introduced in \eqref{fator de bernoulli},
$\delta_{\sigma_i,\sigma_j}$ is Kronecker's delta function
and $\mathscr{Z}^{\mathrm{ ES}}_{\pmb{p},\pmb{h},G}$ is the 
partition function.
\begin{lemma}\label{lema-fev1}
Let $G=(V,E)$ be a finite graph and
consider the $2$-state Potts model on $G$ with free 
boundary condition and Hamiltonian
given by \eqref{rescrever-Potts}. 
Suppose that $p_{ij}\equiv 1-\exp(-2\beta J_{ij})$
and $\hat{\sigma}$ is obtained from $\sigma$ as in the 
Proposition \ref{Medeq}. Then
	\begin{eqnarray*}
	&\exp\big(
	-2\beta
	(\mathscr{H}^{\mathrm{ Potts}}
	_{\pmb{h},2,V}(\hat{\sigma}) +\sum_{\{i,j\}\in E}J_{ij} 
	)
	\big)&
	\\
	&=&
	\\
	&
	\displaystyle
	\sum_{\omega}
	\big(
	\prod_{\{i,j\}: \omega_{ij}=1}
	p_{ij} \delta_{\sigma_i, \sigma_j}
	\prod_{\{i,j\}: \omega_{ij}=0}
	(1-p_{ij})
	\big)
	\times
	{\textstyle 
	\exp
	\big( 
	\beta\sum_{i\in V} h_i 
	(\delta_{\sigma_i,1}-\delta_{\sigma_i,-1})
	\big)
	}.
	\end{eqnarray*}
\end{lemma}
\begin{proof} 
	Using the relation between $\hat{\sigma}$ and $\sigma$,
	we first obtain the following equality
	\[
	\begin{array}{c}
	\exp
	\big(
	-2\beta
	(\mathscr{H}^{\mathrm{ Potts}}
	_{\pmb{h},2,V}(\hat{\sigma}) +\sum_{\{i,j\}\in E}J_{ij}
	)
	\big)
	\\
	=
	\\
	\textstyle
	\exp
	\big( 2\beta 
	(\sum_{\{i,j\}\in E}J_{ij} (\delta_{\sigma_i,\sigma_j}-1)
	+
	\sum_{i\in V} h_i {1\over 2} 
	(
	\delta_{\sigma_i,1}
	-\delta_{\sigma_i,-1} 
	)
	)
	\big).
	\end{array}
	\]
	By using that $p_{ij}=1-\exp(-2\beta J_{ij})$ 
	and the elementary properties of the exponential, 
	a straightforward computation shows that 
	the above expression equals
	\[
	\prod_{\{i,j\}\in E}
	\big(
	p_{ij} \delta_{\sigma_i, \sigma_j} + (1-p_{ij})
	\big)
	\times
	{\textstyle
	\exp
	\big( \beta 
	\sum_{i\in V} h_i 
	(\delta_{\sigma_i,1}-\delta_{\sigma_i,-1} )
	\big)
	}.
	\]
	By expanding the product in the above expression we get
	\begin{eqnarray*}
	\sum_{E'\subset E}
	\!\!
	\big(
	\prod_{\{i,j\}\in E'}
	\!\!p_{ij} \delta_{\sigma_i, \sigma_j}
	\!\!
	\prod_{\{i,j\}\in E\setminus E'}
	\!\!\!\!\!(1-p_{ij})
	\big) 
	\times
	{\textstyle
	\exp
	\big( \beta 
	\sum_{i\in V} h_i 
	(\delta_{\sigma_i,1}-\delta_{\sigma_i,-1})
	\big)
	}
	. 
	\end{eqnarray*}
	Since the collection of all the induced subgraphs of $G$
	is in bijection with $\{0,1\}^E$ we can rewrite 
	the last expression as follows
	\begin{eqnarray*}
		\sum_{\omega}
		\big(
		\prod_{\{i,j\}: \omega_{ij}=1}
		\!\!
		p_{ij} \delta_{\sigma_i, \sigma_j}
		\prod_{\{i,j\}: \omega_{ij}=0}
		\!\!\!(1-p_{ij})
		\big)
		\times
		{\textstyle
		\exp
		\big( \beta 
		\sum_{i\in V} h_i 
		(\delta_{\sigma_i,1}-\delta_{\sigma_i,-1})
		\big)
		}
	\end{eqnarray*}
	which completes  the proof.
\end{proof}
\begin{lemma} \label{Lema5dez}
Under the hypothesis of Lemma \ref{lema-fev1}
we can show that there exists a constant 
$C\equiv C(\beta,G)>0$ so that
	\begin{eqnarray*}
	\mathscr{Z}^{\mathrm{{Potts}}}_{2\beta,\pmb{h},2,V}
	=
	C 
	\mathscr{Z}^{\mathrm{{ES}}}_{\pmb{p},\pmb{h},G}.
	\end{eqnarray*}
\end{lemma}
\begin{proof}
Observe that
	\[ 
	\begin{array}{rcl}
	\mathscr{Z}^{\mathrm{{Potts}}
	}_{2\beta, \pmb{h },2,V}
	&=&
	\displaystyle
	\sum_{\hat{\sigma}\in\{1,2\}^V}
	\exp
	\big(
	-2\beta \mathscr{H}_{\pmb{h},2,V}^{\mathrm{ Potts}}
	(\hat{\sigma})
	\big)
	\\\\[0,2cm]
	&=&
	{\displaystyle
		{
			1
			\over 
			C
		}
		\sum_{\hat{\sigma}\in\{1,2\}^V}
	}
	\exp
	\big(
	-2\beta 
	(\mathscr{H}_{\pmb{h},2,V}^{\mathrm{ Potts}}
	(\hat{\sigma})
	+\sum_{\{i,j\}\in E}J_{ij} 
	)
	\big),
	\end{array}
	\]
where 	
$
C\equiv 
\exp
(
2\beta\sum_{\{i,j\}\in E}J_{ij}
)
>0.
$
From Lemma \ref{lema-fev1} it follows that 
the right-hand-side above is equal to
	\[
	\begin{array}{c}
		\displaystyle
		C
		\sum_{\sigma, \omega}
		\!
		\big(
		\!
		\prod_{\{i,j\}: \omega_{ij}=1}
		p_{ij} \delta_{\sigma_i, \sigma_j}
		\prod_{\{i,j\}: \omega_{ij}=0}
		\!\!\!\!(1-p_{ij})
		\big)
		\!\!
		\times
		\!
		{\textstyle
		\exp
		\big( \beta 
		\sum_{i\in V} h_i 
		(\delta_{\sigma_i,1}-\delta_{\sigma_i,-1})
		\big)
		}
		\\
		=
		\\[0,2cm]
		\displaystyle
		C
		\sum_{\sigma, \omega}
		\!
		\big(
		\!
		\prod_{\{i,j\}: \omega_{ij}=1}
		\!\!\! p_{ij}
		\prod_{\{i,j\}: \omega_{ij}=0}
		\!\!\!\!(1-p_{ij})
		\big)
		\Delta(\sigma, \omega)
		\!\!
		\times
		\!
		{\textstyle
		\exp
		\big( \beta 
		\sum_{i\in V} h_i 
		(\delta_{\sigma_i,1}-\delta_{\sigma_i,-1})
		\big)
		}
		\\
		=
		\\[0,2cm]
		\displaystyle
		C
		\sum_{\sigma,\omega}
		B_{\pmb{J}}(\omega)
		\Delta(\sigma, \omega)
		\times
		{\textstyle
		\exp
		\big( \beta 
		\sum_{i\in V} h_i 
		(\delta_{\sigma_i,1}-\delta_{\sigma_i,-1})
		\big)
		}
		=
		C \mathscr{Z}_{\pmb{p},\pmb{h},G}^{\mathrm{ ES}}.
		\qedhere
	\end{array}
	\]
\end{proof}
\begin{lemma}\label{Lema-cont} 
	Let $G$ be a finite graph and fix an edge configuration 
	$\omega\in\{0,1\}^E$. 
	If $\pmb{h}(K_\alpha)\equiv \beta 
	{\sum_{i\in K_\alpha}} h_i$, where 
	$K_1,\ldots, K_{k({\omega,G})}$ 
	denote the connected components of $(V,\eta(\omega))$ we 
	have
	\begin{eqnarray*}
		\sum_{\sigma\in\{-1,+1\}^V}
		\Delta(\sigma, \omega)
		\times
		\exp 
		\big( \beta 
		\sum_{i\in V} h_i 
		(\delta_{\sigma_i,1}-\delta_{\sigma_i,-1})
		\big)
		=
		\prod_{\alpha=1}^{k(\omega,G)}
		2 \cosh\left(\textbf{h}(K_\alpha) \right).
	\end{eqnarray*}
\end{lemma}
\begin{proof}
	For a fixed $\omega$, if $K_1,\ldots,K_{k(\omega,G)}$ 
	denotes the decomposition of $V$ on its connected components,
	we get	
	\[
	\sum_{i\in V} h_i 
	(\delta_{\sigma_i,1}-\delta_{\sigma_i,-1})
	=
	\sum_{\alpha =1}^{k(\omega,G)} 
	\sum_{i\in K_{\alpha }} 
	h_i 
	(
	\delta_{\sigma_i,1}-\delta_{\sigma_i,-1} 
	).
	\]
	For each spin configuration $\sigma\in\{-1,1\}^{V}$ 
	so that the pair 
	$
	(\sigma,\omega)\in \{-1,1\}^{V}\times\{0,1\}^{E}
	$
	satisfies $\Delta(\sigma,\omega)=1$, we have
	that the value of all the spins in the same component
	has the same sign, see Figure \ref{f2}.
	From the elementary properties of the exponential function
	we obtain the following equality 
	\[
	\begin{array}{c}
		\Delta(\sigma, \omega)
		\times 
		\exp 
		\big( \beta  
		\sum_{i\in V} h_i 
		(\delta_{\sigma_i,1}-\delta_{\sigma_i,-1})
		\big)
		\\[0.3cm]
		=
		\\
		{\displaystyle 
		\Delta(\sigma, \omega)
		\times
		\prod_{\alpha=1}^{k(\omega,G)}
		}
		\exp
		\big(\beta 
		\sum_{i\in K_\alpha} h_i 
		(
		\delta_{\sigma_i,1}-\delta_{\sigma_i,-1} 
		)
		\big)
		.
	\end{array}
	\]
	Since 
	$V=\sqcup_{\alpha=1}^{k(\omega,G)}K_{\alpha}$,
	we have a natural bijection between the following spaces:
	\begin{eqnarray*}
		\{-1,+1\}^V \quad \text{and} \quad 
		\prod_{\alpha=1}^{k(\omega,G)}\{-1,1\}^{K_\alpha}.
	\end{eqnarray*}
	For sake of simplicity, we denote a generic element  
	of the cartesian product  
	$\prod_{\alpha=1}^{k(\omega,G)}\{-1,1\}^{K_\alpha}$
	by $(\sigma_{K_1},\ldots,\sigma_{K_{k(\omega,G)}})$,
	where $\sigma_{K_j}\equiv({\sigma_i:i\in K_j}),$ $\, 
	\forall j=1,\ldots, k(\omega,G).$
	In this way we can simply write
	$
	\sigma= (\sigma_{K_1},\ldots,\sigma_{K_{k(\omega,G)}}).
	$
	By abusing the notation, we write 
	\[
	\Delta(\sigma,\omega)= \prod_{\alpha=1}^{k(\omega,G)}
	\Delta(\sigma_{K_\alpha},\omega).
	\]
	Using the previous observations and 
	$
	h_i 
	(
	\delta_{\sigma_i,1}-\delta_{\sigma_i,-1} 
	)
	=
	h_i\sigma_i
	$
	we obtain
	\begin{multline*}
	\sum_{\sigma}
	\Delta(\sigma, \omega)
	\times
	{\textstyle 
	\exp 
	\big( \beta  
	\sum_{i\in V} h_i 
	(\delta_{\sigma_i,1}-\delta_{\sigma_i,-1})
	\big)
	}	
	\\	
	=
	\sum_{(\sigma_{K_1},\ldots,\sigma_{K_{k(\omega,G)}})} 
	\prod_{\alpha=1}^{k(\omega,G)}
	\Delta(\sigma_{K\alpha}, \omega)
	{\textstyle
	\exp
	\big(
	\beta 
	\sum_{i\in K_\alpha} h_i \sigma_i
	\big)
	}
	\\[0,2cm]
	=
	\prod_{\alpha=1}^{k(\omega,G)}
	\sum_{\sigma_{K_\alpha}} 
	\Delta(\sigma_{K\alpha}, \omega)
	{\textstyle
	\exp
	\big(
	\beta 
	\sum_{i\in K_\alpha} h_i \sigma_i
	\big).
	}
	\end{multline*}
	Because of the consistency condition,
	for each fixed $\alpha$,  
	the sum appearing on the right-hand-side 
	above has exactly two non zero terms 
	where the spins in $K_{\alpha}$ 
	take the values $\pm 1$.
	Therefore the above expression reduces to
	\begin{eqnarray*}
		\prod_{\alpha=1}^{k(\omega,G)}
		2\cosh
		\big(
		\beta 
		\sum_{i\in K_\alpha} h_i 
		\big),
	\end{eqnarray*}
	thus the lemma is proved.
\end{proof}
\begin{lemma}\label{lema-fev2}
	For any finite graph $G=(V,E)$, 
	$\pmb{p}=(p_{ij}:i,j\in \mathbb{V})$, $\pmb{h}=(h_i:i\in V)$
	and $\beta>0$ we have
	\[
	\mathscr{Z}_{\pmb{p},\pmb{h},G}^{\mathrm{ ES}}
	=
	\mathscr{Z}_{\pmb{p},\pmb{h},G}^{\mathrm{ RC}}.
	\]
\end{lemma}
\begin{proof} 
	The proof of this lemma is trivial given the above results.
	It is enough to change the sum order of the partition 
	function of the Edwards-Sokal model, apply Lemma 
	\ref{Lema-cont} and finally use the definition of the 
	partition function of the random-cluster model as we show 
	below
	\begin{eqnarray*} 
		\mathscr{Z}_{\pmb{p},\pmb{h},G}^{\mathrm{ ES}}
		&=&
		\sum_\omega
		B_{\pmb{J}}(\omega)
		\sum_\sigma
		\Delta(\sigma, \omega)
		\times
		{\textstyle
		\exp
		\big( 
		\beta 
		\sum_{i\in V} h_i 
		(\delta_{\sigma_i,1}-\delta_{\sigma_i,-1})
		\big)
		}
		\\\\[0,2cm]
		&=&
		\sum_\omega
		B_{\pmb{J}}(\omega)
		\prod_{\alpha=1}^{k(\omega,G)}
		2 \cosh\left(\textbf{h}(K_\alpha) \right)
		=
		\mathscr{Z}_{\pmb{p},\pmb{h},G}^{\mathrm{ RC}}.
	\end{eqnarray*}
\end{proof}
In the sequel we prove the main result of this section.
The technicalities of the proof were worked out in the previous 
lemmas and now the remaining task is to simply piece them 
together. 
\begin{theorem}[The marginals of ${\nu}_{\pmb{p},\pmb{h},G}$] 
	\label{lemin}
	Let $G=(V,E)$ be a finite graph, $\beta>0$, 
	$\pmb{p}=(p_{ij}:i,j\in \mathbb{V})$
	as above and $\pmb{h}=(h_i:i\in V)$ an external field. Then 
	\noindent
	\begin{itemize}
		\item[(1)] 
		$
		\begin{aligned}
		\displaystyle{\sum_{\omega\in \{0,1\}^E} 
			{\nu}_{\pmb{p},\pmb{h},G}(\sigma, \omega) 
			=\lambda_{\beta,\pmb{h},V}(\sigma)}
		\quad\qquad\quad \ 
		(\text{spin-marginal of} \ {\nu}_{\pmb{p},\pmb{h},G} )
		\end{aligned}
		$
		\item[(2)] 
		$
		\begin{aligned}
		\displaystyle{\sum_{\sigma\in\{-1,+1\}^V} 
			{\nu}_{\pmb{p},\pmb{h},G}(\sigma, \omega)
			=\phi_{\pmb{p},\pmb{h},G}(\omega)}.
		\quad\qquad
		(\text{edge-marginal of} \ {\nu}_{\pmb{p},\pmb{h},G})
		\end{aligned}
		$
	\end{itemize}
\end{theorem}
\begin{proof} We first prove \textit{(1)}.
	Using the definition of Bernoulli factors $B_{\pmb{J}}$,
	and Lemma \ref{Lema5dez} with  
	$C\equiv\exp (2\beta\sum_{\{i,j\}\in E} J_{ij})>0$,
	we obtain
	\begin{multline*}
		\sum_{\omega}
		{\nu}_{\pmb{p},\pmb{h},G}(\sigma, \omega)
		\\[0,1cm] 
		=
		{
		1 
		\over 
		\mathscr{Z}^{\mathrm{{ES}}}_{\pmb{p},\pmb{h},G}
		}
		\sum_{\omega}
		B_{\pmb{J}}(\omega)
		\Delta(\sigma, \omega)
		\times
		{\textstyle 
		\exp
		\big( 
		\beta 
		\sum_{i\in V} h_i 
		(\delta_{\sigma_i,1}-\delta_{\sigma_i,-1})
		\big)
		}
		\\[0,2cm]
		=
		{
		C 
		\over 
		\mathscr{Z}^{\mathrm{{Potts}}}_{2\beta,\pmb{h},2,V}
		}
		\sum_{\omega}
		{\textstyle 
		\big(
		\prod_{\{i,j\}: \omega_{ij}=1}
		p_{ij} \delta_{\sigma_i, \sigma_j}
		\prod_{\{i,j\}: \omega_{ij}=0}
		(1-p_{ij})
		\big)
		}\times
		\\[0,2cm]
		\times
		{\textstyle
		\exp
		\big( \beta 
		\sum_{i\in V} h_i 
		(\delta_{\sigma_i,1}-\delta_{\sigma_i,-1})
		\big)
		}.
	\end{multline*}
	By applying Lemma \ref{lema-fev1}, it follows that the rhs 
	above is equal to
	\begin{eqnarray*} 
	{C \over \mathscr{Z}^{\mathrm{{Potts}}}_{2\beta,\pmb{h},V}}
	\exp
	{\textstyle
	\big(
		-2\beta
		(
			\mathscr{H}_{\pmb{h},2,V}^{\mathrm{ Potts}}
			(\hat{\sigma})
			+ \sum_{\{i,j\}\in E}J_{ij} 
		) 
	\big)
	}
	=
	\pi_{{2\beta},\pmb{h},2,V}(\hat{\sigma})
	=
	\lambda_{\beta, \pmb{h}, V}(\sigma).
	\end{eqnarray*}
	To prove \textit{(2)} it is enough to use Lemmas 
	\ref{Lema-cont} and \ref{lema-fev2} as follows 
	\begin{eqnarray*}\label{importante_1}
		\sum_{\sigma}
		{\nu}_{\pmb{p},\pmb{h},G}(\sigma, \omega)
		&=&
		{
			B_{\pmb{J}}(\omega) 
		\over 
		\mathscr{Z}^{\mathrm{{ES}}}_{\pmb{p},\pmb{h},G}
		}
		\sum_{\sigma}
		\Delta(\sigma, \omega)
		\times
		{\textstyle
		\exp
		\big( \beta 
		\sum_{i\in V} h_i 
		(\delta_{\sigma_i,1}-\delta_{\sigma_i,-1})
		\big)
		}\nonumber
		\\[0,2cm]
		&=&
		{
			1 
			\over 
			\mathscr{Z}^{\mathrm{{RC}}}_{\pmb{p},\textbf{h},G}
		}
		\displaystyle{
			B_{\pmb{J}}(\omega)
			\prod_{\alpha=1}^{k(\omega,G)}
			2
			\cosh\left(\pmb{h}(K_\alpha) \right)
		}
		=
		\phi_{\pmb{p},\pmb{h},G}(\omega).
	\end{eqnarray*}
\end{proof}
\begin{corollary}[Conditional measure of $\nu_{\beta,\pmb{h},G}$] 
	\label{medida-conditional1}
    Let $\omega\in\{0,1\}^E$ be a fixed edge configuration. 
	For each $\sigma\in\{-1,+1\}^V$ we have that
	\[
	\nu_{\pmb{p},\pmb{h},G}(\sigma\vert \omega)=
	{
		\Delta(\sigma, \omega)
		\over
		\prod_{\alpha=1}^{k(\omega,G)}
		2
		\cosh\left(\pmb{h}(K_\alpha) \right)
	}
	\times
	{\textstyle
	\exp
	\big( \beta 
	\sum_{i\in V} h_i 
	(\delta_{\sigma_i,1}-\delta_{\sigma_i,-1})
	\big).
	}
	\]
\end{corollary}
\begin{proof} 
	From Proposition \ref{Medeq} and
	Lemma \ref{lema-fev1} we have,
	for any random variable $g:\{-1,1\}^V\to \mathbb{R}$,
	that 
	\begin{multline*}
		\lambda_{\beta,\pmb{h},V}(g)
		=
		\sum_{\hat{\sigma}\in\{1,2\}^V}g(\hat{\sigma}) 
		\pi_{{2\beta},\pmb{h},2,V}(\hat{\sigma})
		\\[0,2cm]
		=
		{
		C
		\over 
		\mathscr{Z}^{\mathrm{{Potts}}}_{2\beta,\pmb{h},2,V}
		}
		\sum_{\sigma\in\{-1,+1\}^V}
		g(\sigma)
		\sum_{\omega\in\{0,1\}^E}
		\big(
		\prod_{\{i,j\}: \omega_{ij}=1}
		p_{ij} \delta_{\sigma_i, \sigma_j}
		\prod_{\{i,j\}: \omega_{ij}=0}
		(1-p_{ij})
		\big)\times
		\\[0,2cm]
		\hspace*{5cm}
		\times
		\exp
		\big( \beta 
		\sum_{i\in V} h_i 
		(\delta_{\sigma_i,1}-\delta_{\sigma_i,-1})
		\big),
	\end{multline*}
	where 
	$
	C
	\equiv 
	\exp
	(
	2\beta\sum_{\{i,j\}\in E}J_{ij}
	)
	>
	0
	$.
	By changing the order of the sums in the last expression we 
	get
	\begin{eqnarray*}
		{
			C
			\over
			 \mathscr{Z}^{\mathrm{{Potts}}}_{2\beta,\pmb{h},2,V}
		}
		\sum_{\omega\in\{0,1\}^E}
		B_{\pmb{J}}(\omega) 
		\!\!
		\sum_{\sigma\in\{-1,+1\}^E}
		\!\!\!\!\!\!
		g(\sigma)\,
		\Delta(\sigma,\omega)\,
		{\textstyle
		\exp
		\big( \beta 
		\sum_{i\in V} h_i 
		(\delta_{\sigma_i,1}-\delta_{\sigma_i,-1})
		\big)
		}.
	\end{eqnarray*}
	According to Lemmas \ref{Lema5dez} and \ref{lema-fev2},
	we have 
	$C^{-1} \mathscr{Z}^{\mathrm{{Potts}}}_{2\beta,\pmb{h},2,V}
	=
	\mathscr{Z}^{\mathrm{{ES}}}_{\pmb{p},\pmb{h},G}
	=
	\mathscr{Z}^{\mathrm{{RC}}}_{\pmb{p},\pmb{h},G}
	$.
	By introducing the product appearing in the definition 
	of the random-cluster model, we can see that the above 
	expression is equal to 
	\begin{multline*}
		\sum_{\omega}
		{
			B_{\pmb{J}}(\omega) 
			\prod_{\alpha=1}^{k(\omega,G)}
			2
			\cosh\left(\pmb{h}(K_\alpha) \right)
			\over 
			\mathscr{Z}^{\mathrm{{RC}}}_{\pmb{p},\pmb{h},G}
		}\times
		\\
		\times\sum_{\sigma}
		g(\sigma)
		{
			\Delta(\sigma,\omega)
			\times
			\exp
			\big(
			\beta 
			\sum_{i\in V} h_i 
			(\delta_{\sigma_i,1}-\delta_{\sigma_i,-1})
			\big)
			\over
			\prod_{\alpha=1}^{k(\omega,G)}
			2
			\cosh\left(\pmb{h}(K_\alpha) \right)
		}
		\\
		\\[0,2cm]
		\!\!\!\!
		=
		\sum_{\omega,\sigma}
		g(\sigma)
		{
			\Delta(\sigma,\omega)
			\times
			\exp
			\big( \beta 
			\sum_{i\in V} h_i 
			(\delta_{\sigma_i,1}-\delta_{\sigma_i,-1})
			\big)
			\over
			\prod_{\alpha=1}^{k(\omega,G)}
			2
			\cosh\left(\pmb{h}(K_\alpha) \right)
		}\phi_{\pmb{p},\pmb{h},G}(\omega).
	\end{multline*}
	Therefore
	\begin{eqnarray*}
	\lambda_{\beta,\pmb{h},V}(g)=
	\sum_{\omega}
	\left[
	\sum_\sigma
	g(\sigma)
	{
		\Delta(\sigma,\omega)
		\exp
		\big( 
		\beta 
		\sum_{i\in V} h_i 
		(\delta_{\sigma_i,1}-\delta_{\sigma_i,-1})
		\big)
		\over
		\prod_{\alpha=1}^{k(\omega,G)}
		2
		\cosh\left(\pmb{h}(K_\alpha) \right)
	}
	\right]
	\phi_{\pmb{p},\pmb{h},G}(\omega).
	\end{eqnarray*}
	On the other hand, we get from Theorem \ref{lemin} that
	\begin{eqnarray}\label{mart2}
	\lambda_{\beta,\pmb{h},V}(g)
	=
	\sum_{\omega,\sigma}
	g(\sigma) \nu_{\pmb{p},\pmb{h},G}(\sigma,\omega)  \nonumber
	=
	\sum_\omega
	\Big[
	\sum_\sigma g(\sigma)
	\nu_{\pmb{p},\pmb{h},G}(\sigma\vert\omega)
	\Big]
	\phi_{\pmb{p},\pmb{h},G}(\omega).
	\end{eqnarray}
	The proof is completed upon comparison of the two previous 
	expressions.
\end{proof}
\section{Two-point function}
The {\bf two-point function} of the $q$-state Potts model 
is defined by
	\[
	\tau_{\beta,\pmb{\hat{h}},q,V}(x,y)
	\equiv
	\pi_{\beta,\pmb{\hat{h}},q,V}
	(\hat{\sigma}_x=\hat{\sigma}_y)
	-
	{1\over q}. 
	\]
The term $1/q$ represents the probability that 
two independent spins uniformly chosen have the same value.
In the random-cluster model, the {\bf connectivity function} 
plays the role of the two-point function of the Potts model.
This function is precisely the probability, with 
respect to  $\phi_{\pmb{p},\pmb{h},G}$, 
that $x$ and $y$ are in the same connected component, notation
$
\phi_{\pmb{p},\pmb{h},G}(x\leftrightarrow y).
$
\begin{lemma}\label{lema28fev} 
	Let $G=(V,E)$ be a finite graph and $x,y\in V$ 
	two distinct vertices. 
	Fix an edge configuration $\omega\in\{0,1\}^{E}$.
	If $x\not\leftrightarrow y$ in $\omega$, 
	then 
	\begin{multline*} 
	\sum_{\sigma\in\{-1,+1\}^V}
	\mathds{1}_{\{\sigma_x=\sigma_y\}}\Delta(\sigma,\omega)
	{\textstyle 
	\exp
	\big(
	\beta\sum_{i\in V}h_i\sigma_i
	\big)
	}
	\\
	=
	2\cosh\left(\pmb{h}(K_t)+\pmb{h}(K_u)\right)
	\prod_{\substack{\alpha=1 \\ \alpha\not=t,u}}^{k(\omega,G)}
	2\cosh\big(\pmb{h}(K_\alpha)\big),
	\end{multline*}
	where $K_t\equiv K_t(\omega)$ and $K_u\equiv K_u(\omega)$ are
	two disjoint connected components containing the vertices
	$x$ and $y$, respectively.
\end{lemma}
\begin{proof}
	The basic ideas used to prove this lemma are the 
	same we employed to prove Lemma \ref{Lema-cont}, 
	which we once more present for sake of completeness.
	Let  
	$V = \sqcup_{\alpha=1}^{k(\omega,G)}K_{\alpha}$
	be a decomposition in terms of the connected components of 
	the graph $(V,\eta(\omega))$. We recall that   
	$
		\{-1,+1\}^V 
		\cong 
		\prod_{\alpha=1}^{k(\omega,G)}\{-1,1\}^{K_\alpha} 
	$
	and its elements are denoted by  
	$(\sigma_{K_1},\ldots,\sigma_{K_{k(\omega,G)}})$,
	with 
	$\sigma_{K_j}\equiv({\sigma_i:i\in K_j}),$ 
	$\ \forall j=1,\ldots, k(\omega,G).$
	We also use the natural identification
	$
	\sigma = (\sigma_{K_1},\ldots,\sigma_{K_{k(\omega,G)}}).
	$
	
	Suppose that $x\not\leftrightarrow y$ in $\omega$. 
	Denote $K_t$ and $K_u$ the components containing 
	the vertices $x$ and $y$, respectively. 
	Taking into account the decomposition of $V$ mentioned above,
	we have 
	\begin{multline*}
		\sum_{\sigma}
		\mathds{1}_{\{\sigma_x=\sigma_y\}}
		\Delta(\sigma,\omega)
		{\textstyle
		\exp
		\big(
		\beta\sum_{i\in V}h_i\sigma_i
		\big)
		}
		\\
		=
		\sum_{\sigma}
		\mathds{1}_{\{\sigma_x=\sigma_y=+1\}}
		\Delta(\sigma,\omega)
		\prod_{\alpha=1}^{k(\omega,G)}
		{\textstyle
		\exp
		\big(
		\beta\sum_{i\in K_\alpha}h_i\sigma_i
		\big)
		}
		\\
		\hspace*{4cm}+
		\sum_{\sigma}
		\mathds{1}_{\{\sigma_x=\sigma_y=-1\}}
		\Delta(\sigma,\omega)
		\prod_{\alpha=1}^{k(\omega,G)}
		{\textstyle
		\exp
		\big(
		\beta\sum_{i\in K_\alpha}h_i\sigma_i
		\big).
		}
	\end{multline*}
	As previously observed, 
	$
	\Delta(\sigma,\omega)
	=
	\prod_{\alpha=1}^{k(\omega,G)}
	\Delta(\sigma_{K_\alpha},\omega)
	$,
	so from a simple computation we get that the expression 
	above is equal to
	\begin{eqnarray}\label{bld}
	2\cosh\big(\pmb{h}(K_t)+\pmb{h}(K_u)\big)
	\prod_{\substack{\alpha=1 \\ \alpha\not=t,u}}^{k(\omega,G)}
	\sum_{\sigma_{K_\alpha}}
	\Delta(\sigma_{K_\alpha},\omega)
	\exp
	\big(
	\beta\sum_{i\in K_\alpha}h_i\sigma_i
	\big).
	\end{eqnarray}
	Because of the consistency condition, the sums above over 
	$\sigma_{K_\alpha}$ have actually two non-null terms.
	In each of such term the value of the spins is constant 
	and therefore the product simplifies to 
	\begin{multline}
	\prod_{\substack{\alpha=1 \\ \alpha\not=t,u}}^{k(\omega,G)}
	\sum_{\sigma_{K_\alpha}}
	\Delta(\sigma_{K_\alpha},\omega)
	{\textstyle 
	\exp
	\big(
	\beta\sum_{i\in K_\alpha}h_i\sigma_i
	\big)
	}
	\\
	=
	\prod_{\substack{\alpha=1 \\ \alpha\not=t,u}}^{k(\omega,G)}
	\big(
	\exp(\pmb{h}(K_\alpha))
	+
	\exp(-\pmb{h}(K_\alpha))
	\big) \nonumber
	=
	\prod_{\substack{\alpha=1 \\ \alpha\not=t,u}}^{k(\omega,G)}
	2\cosh\big( \pmb{h}(K_\alpha)\big).
	\end{multline}
	Finally, by replacing the last expression in \eqref{bld}, 
	we end the proof.
\end{proof}
Lemma \ref{lema28fev} is vital for the most important
result of this section, which is the next theorem. 
We state below the theorem for the $2$-state Potts model,
but in fact the theorem is valid for much more general Potts
models. The general case is treated in the last section.
\begin{theorem}[Correlation-connectivity] 
	\label{Correlação/quase-conectividade}
	Let $G=(V,E)$ be a finite graph and $x,y$
	two distinct vertices in $V$. 
	Then
	\begin{multline*}
		\hspace*{-0.32cm}\tau_{2\beta,\pmb{h},2,V}(x,y)
		%\\
		\!
		=
		\! 
		{1\over 2}
		\phi_{\pmb{p},\pmb{h},G}(x\leftrightarrow y)
		+
		{1\over 2}
		\phi_{\pmb{p},\pmb{h},G}
		\big(
		\mathds{1}_{\{x\not\leftrightarrow y\}}
		\cdot
		\tanh\left(\pmb{h}(K_t)\right)
		\cdot 
		\tanh\left(\pmb{h}(K_u)\right)
		\!\big),
	\end{multline*}
	where $K_t\equiv K_t(\omega)$ e $K_u\equiv K_u(\omega)$ 
	are two disjoint connected components containing the
	vertices $x$ and $y$, respectively.
\end{theorem}
\begin{proof} 
	By using the definition of the two-point function
	and Theorem \ref{lemin} we get 
	\begin{eqnarray*}
	\tau_{2\beta,\pmb{h},2,V}(x,y)
	&=&
	\pi_{2\beta,\pmb{h},2,V}(\hat{\sigma}_x=\hat{\sigma}_y)
	-
	{1\over 2}
	\nonumber
	\\[0,2cm]
	&=&
	\sum_{\hat{\sigma}\in\{1,2\}^V} 
	\Big(\mathds{1}_{\{\hat{\sigma}_x=\hat{\sigma}_y\}}
	-
	{1\over 2}\Big)
	\pi_{2\beta,\pmb{h},2,V}(\hat{\sigma}) 
	\\[0,2cm]
	&=&
	\sum_{(\sigma,\omega)\in\{-1,+1\}^V\times\{0,1\}^E} 
	\Big(\mathds{1}_{\{\sigma_x=\sigma_y\}}-{1\over 2}\Big)
	\nu_{\beta,\pmb{h},G}(\sigma,\omega) \nonumber
	\\[0,2cm]
	&=&
	\sum_{\omega\in\{0,1\}^E}
	\Big[
	\sum_{\sigma\in\{-1,+1\}^V} 
	\Big(\mathds{1}_{\{\sigma_x=\sigma_y\}}-{1\over 2}\Big)
	\nu_{\beta,\pmb{h},G}(\sigma\vert \omega)
	\Big]
	\phi_{\pmb{p},\pmb{h},G}(\omega).\nonumber
	\end{eqnarray*}
	Since 
	$
	\delta_{\sigma_i,\sigma_j}
	=
	{1\over 2}(1+\sigma_i\sigma_j)
	$,
	it follows from Corollary \ref{medida-conditional1} that the 
	rhs above is 
	\begin{multline}\label{I-1+I-2}
	=
	\sum_{\omega}
	\left[
	\sum_{\sigma} 
	\Big(\mathds{1}_{\{\sigma_x=\sigma_y\}}-{1\over 2}\Big)
	{
		\Delta(\sigma, \omega)
		\exp
		\big( 
		\beta 
		\sum_{i\in V} h_i\sigma_i
		\big)
		\over 
		\prod_{\alpha=1}^{k(\omega,G)}
		2
		\cosh\left(\pmb{h}(K_\alpha) \right)
	}
	\right]
	\phi_{\pmb{p},\pmb{h},G}(\omega)\nonumber
	\\[0,3cm]
	=
	\sum_{\omega}
	\mathds{1}_{\{x\leftrightarrow y\}}(\omega)
	\left[
	\sum_{\sigma} 
	\Big(\mathds{1}_{\{\sigma_x=\sigma_y\}}-{1\over 2}\Big)
	{
		\Delta(\sigma, \omega)
		\exp
		\big( \beta 
		\sum_{i\in V} h_i\sigma_i
		\big)
		\over 
		\prod_{\alpha=1}^{k(\omega,G)}
		2
		\cosh\left(\pmb{h}(K_\alpha) \right)
	}
	\right]
	\phi_{\pmb{p},\pmb{h},G}(\omega)\nonumber
	\\[0,3cm]
	+
	\sum_{\omega}
	\mathds{1}_{\{x\not\leftrightarrow y\}}(\omega)
	\left[
	\sum_{\sigma} 
	\Big(\mathds{1}_{\{\sigma_x=\sigma_y\}}-{1\over 2}\Big)
	{
		\Delta(\sigma, \omega)
		\exp
		\big( \beta 
		\sum_{i\in V} h_i\sigma_i
		\big)
		\over 
		\prod_{\alpha=1}^{k(\omega,G)}
		2
		\cosh\left(\pmb{h}(K_\alpha) \right)
	}
	\right]
	\phi_{\pmb{p},\pmb{h},G}(\omega)\nonumber
	\\[0,3cm]
	\equiv I_1+I_2.
	\end{multline}
	Notice that, as long as $x\leftrightarrow y$ in $\omega$ 
	and the pair $(\sigma,\omega)$ is consistent, then 
	$\sigma_x=\sigma_y$.
	From this observation and Lemma \ref{Lema-cont}, it follows 
	that
	\begin{eqnarray*}
	I_1={1\over 2}\phi_{\pmb{p},\pmb{h},G}(x\leftrightarrow y).
	\end{eqnarray*}
	On the other hand, applying Lemma \ref{Lema-cont} again 
	yields
	\begin{eqnarray}\label{I-2}
	I_2&=&
	-{1\over 2}
	\phi_{\pmb{p},\pmb{h},G}(x\not\leftrightarrow y)
	\nonumber
	\\[0,2cm]
	&&\hspace*{0,2cm}
	+
	\sum_{\omega}
	\mathds{1}_{\{x\not\leftrightarrow y\}}(\omega)
	\left[
	\sum_{\sigma} 
	\mathds{1}_{\{\sigma_x=\sigma_y\}}
	{
		\Delta(\sigma, \omega)
		\exp
		\big(
		\beta 
		\sum_{i\in V} h_i\sigma_i
		\big)
		\over 
		\prod_{\alpha=1}^{k(\omega,G)}
		2
		\cosh\left(\pmb{h}(K_\alpha) \right)
	}
	\right]
	\phi_{\pmb{p},\pmb{h},G}(\omega)\nonumber
	\\[0,2cm]
	&\equiv&
	-{1\over 2}\phi_{\pmb{p},\pmb{h},G}(x\not\leftrightarrow y)
	+\tilde{I}_2.
	\end{eqnarray}
	Now we work on $\tilde{I}_2$. 
	By using Lemma \ref{lema28fev}, we have
	\begin{multline}
	\tilde{I}_2
	=
	\sum_{\omega}
	\mathds{1}_{\{x\not\leftrightarrow y\}}
	\left[
	\sum_{\sigma} 
	\mathds{1}_{\{\sigma_x=\sigma_y\}}
	{
		\Delta(\sigma, \omega)
		\exp
		\big(
		\beta 
		\sum_{i\in V} h_i\sigma_i
		\big)
		\over 
		\prod_{\alpha=1}^{k(\omega,G)}
		2
		\cosh\left(\pmb{h}(K_\alpha) \right)
	}
	\right]
	\phi_{\pmb{p},\pmb{h},G}(\omega)
	\nonumber
	\\[0,3cm]
	=
	\sum_{\omega}
	\mathds{1}_{\{x\not\leftrightarrow y\}}(\omega)
	\!\!
	\left[
	{
	2\cosh\left(\pmb{h}(K_t)+\pmb{h}(K_u)\right)
	\prod_{\substack{\alpha=1 \\ \alpha\not=t,u}}^{k(\omega,G)}
	2\cosh\big(\pmb{h}(K_\alpha)\big)
	\over 
	\prod_{\alpha=1}^{k(\omega,G)}
	2
	\cosh\left(\pmb{h}(K_\alpha) \right)
	}
	\right]
	\!\!
	\phi_{\pmb{p},\pmb{h},G}(\omega)
	\nonumber
	\\[0,3cm]
	=
	{1\over 2}
	\phi_{\pmb{p},\pmb{h},G}
	\Big(
	\mathds{1}_{\{x\not\leftrightarrow y\}}
	\cdot
	{
		\cosh\left(\pmb{h}(K_t)+\pmb{h}(K_u)\right)
		\over 
		\cosh\left(\pmb{h}(K_t)\right)\cdot\cosh
		\left(\pmb{h}(K_u)\right)
	}
	\Big)
	\nonumber
	\\[0,3cm]
	=
	{1\over 2}
	\phi_{\pmb{p},\pmb{h},G}
	\big(
	\mathds{1}_{\{x\not\leftrightarrow y\}}
	\cdot
	\big\{
	1
	+
	\tanh\left(\pmb{h}(K_t)\right)
	\cdot
	\tanh\left(\pmb{h}(K_u)\right)
	\big\}
	\big).
	\end{multline}
	Replacing the last expression in \eqref{I-2}
	we get that
	\begin{eqnarray}\label{I-2-1}
	I_2
	=
	{1\over 2}
	\phi_{\pmb{p},\pmb{h},G}
	\big(
	\mathds{1}_{\{x\not\leftrightarrow y\}}
	\cdot
	\tanh\left(\pmb{h}(K_t)\right)
	\cdot 
	\tanh\left(\pmb{h}(K_u)\right)
	\big).
	\end{eqnarray}
	Since $\tau_{2\beta,\pmb{h},2,V}(x,y)=I_1+I_2$, the theorem 
	follows.
\end{proof}
\begin{remark}
	Notice that in the absence of the magnetic field, i.e. 
	$\pmb{h}\equiv 0$, the conclusion 
	of the Theorem \ref{Correlação/quase-conectividade}
	reduces to  
	\[
	\tau_{2\beta,0,2,V}(x,y)=
	{1\over 2}\phi_{\pmb{p},0,G}(x\leftrightarrow y),
	\quad \forall x,y\in V
	\]
	which is a well known identity for the
	Ising/Potts model with $q=2$, see \cite{Grimmett2} 
	Theorem 1.16, p. 11.
\end{remark}
\begin{corollary} 
	\label{col}
	The spin-spin correlation of the Ising model on the finite
	volume $V$ satisfies the 
	following identity for any magnetic field 
	$\pmb{h}\in \mathbb{R}^{V}$
	\[
	\lambda_{{\beta},\pmb{h},V}(\sigma_x\sigma_y )
	=
	\phi_{\pmb{p},\pmb{h},G}(x\leftrightarrow y)
	+
	\phi_{\pmb{p},\pmb{h},G}
	\big(
	\mathds{1}_{\{x\not\leftrightarrow y\}}
	\cdot
	\tanh\left(\pmb{h}(K_t)\right)
	\cdot 
	\tanh\left(\pmb{h}(K_u)\right)
	\big).
	\]
\end{corollary}
\begin{proof}
	This follows easily from the definition of the expected 
	value and Theorem \ref{Correlação/quase-conectividade}
	since 
	$
	\lambda_{{\beta},\pmb{h},V}(\sigma_x\sigma_y )
	=
	\lambda_{{\beta},\pmb{h},V}(\sigma_i=\sigma_j)
	-
	\lambda_{{\beta},\pmb{h},V}(\sigma_i\neq \sigma_j)
	=
	2\lambda_{{\beta},\pmb{h},V}(\sigma_i=\sigma_j)-1
	=
	2
	\left[
	\pi_{2\beta,\pmb{h},2,V}(\hat{\sigma}_x
	=\hat{\sigma}_y)-{1\over 2} 
	\right]
	=
	2\tau_{2\beta,\pmb{h},2,V}(x,y)
	$.
\end{proof}
\begin{remark}
	If we consider the Ising model on $G$ without magnetic field, 
	from Corollary \ref{col} we get 
	\(
	\lambda_{{{\beta},0,V}} (\sigma_x\sigma_y )
	=
	\phi_{\pmb{p},0,G}(x\leftrightarrow y), 
	\ \forall x,y\in V.
	\)
\end{remark}
\section{Applications}
\paragraph{Spin-spin correlations.}
Corollary \ref{col} can be used to obtain some correlation
inequalities. Keeping the notation of Theorem 
\ref{Correlação/quase-conectividade} 
and supposing that $h_i\geqslant 0$ for all $i\in V$, it follows 
from the monotonicity of the hyperbolic tangent that
$\tanh(\beta h_x)\leqslant \tanh(\pmb{h}(K_t))$ 
and 
$\tanh(\beta h_y)\leqslant \tanh(\pmb{h}(K_u))$.
These estimates together with Corollary \ref{col} 
give us the following lower bound 
$
\phi_{\pmb{p},\pmb{h},G}\ (x\not\leftrightarrow y)\allowbreak
\tanh(\beta h_x)\tanh(\beta h_y)\allowbreak
\leqslant \allowbreak
\lambda_{{\beta},\pmb{h},V}(\sigma_x\sigma_y).
$
A simple computation shows that 
$
P_{\pmb{p}}(x\not\leftrightarrow y)
\leqslant
\phi_{\pmb{p},\pmb{h},G}(x\not\leftrightarrow y)
$,
where $P_{\pmb{p}}$ is the probability measure of the 
independent bond percolation model with parameter $\pmb{p}$.
Supposing that $\pmb{p}\equiv p$ (the homogeneous model) 
and $p<p_c(\mathbb{V})$, for any given $\varepsilon>0$, if the 
distance between $x$ and $y$ is large enough then
$
(1-\varepsilon)
\tanh(\beta h_x)\tanh(\beta h_y)
\leqslant
\lambda_{{\beta},\pmb{h},V}(\sigma_x\sigma_y) ,
$
which, of course, can also be (better) obtained by 
the GKS inequality. 

Under the above assumptions, Corollary \ref{col}
also gives us an upper bound in terms of the iid Bernoulli
bond percolation model, which is 
	\(
	\lambda_{{\beta},\pmb{h},V} (\sigma_x\sigma_y)\allowbreak
	\leqslant\allowbreak
	e^{-C(\beta) d_{G}(x,y)}\allowbreak
	+
	P_p
	\big(
	\tanh\left(\pmb{h}(K_t)\right)
	\cdot 
	\tanh\left(\pmb{h}(K_u)\right)
	\big),
	\)
where at this point we are assuming $J_{ij}\equiv J$ and 
$p=1-e^{-\beta J}$.
To obtain the asymptotic behavior of the second term in the rhs 
above, one needs to impose extra conditions on the geometry of 
the graph and the decay ratio of the magnetic field.
\paragraph{Expected value and distribution function of a single 
spin.}
\begin{lemma}\label{lema anteprincipal}
	Consider a finite graph $G=(V,E)$,  $x\in V$ and 
	$\omega\in\{0,1\}^E$ a fixed edge configuration. Then
	\begin{multline*} 
	\sum_{\sigma}
	\!
	\mathds{1}_{\{\sigma_x=\pm 1\}}
	\Delta(\sigma,\omega)
	{\textstyle
	\exp
	\big(\beta\sum_{i\in V}h_i\sigma_i\big)
	}
	\\
	{\textstyle
	\!=\!
	\exp\left( \pm \pmb{h}(K_t)\right)
	\prod_{\substack{\alpha=1 \\ \alpha\not=t}}^{k(\omega,G)}
	}
	2\cosh\big(\pmb{h}(K_\alpha)\big),
	\end{multline*}
	where $K_t\equiv K_t(\omega)$ is the connected component
	containing the vertex $x$.
\end{lemma}
\begin{proof}
	To prove this lemma we proceed, mutatis mutandis, 
	as in the proof of Lemma \ref{lema28fev}.
\end{proof}
\begin{theorem}[Distribution function] 
	\label{teorema de medida} Let $G=(V,E)$ be a finite graph.
	We have, for any fixed $x\in V$, that
	\[
	\lambda_{\beta,\pmb{h},V}(\sigma_x=\pm 1)
	=
	{1\over 2}
	\pm
	{1\over 2}
	\phi_{\pmb{p},\pmb{h},G}
	\big(
		\tanh\left(\pmb{h}(K_t)\right)
	\big),
	\]
	where $K_t(\omega)\equiv K_t$ is the connected component 
	containing $x$.
\end{theorem}
\begin{proof}
	From Theorem \ref{lemin}, it follows that
	\begin{multline}\label{espera1}
	\lambda_{\beta,\pmb{h},V}(\sigma_x=\pm1)
	=
	\sum_{(\sigma,\omega)
	\in
	\{-1,+1\}^V\times\{0,1\}^E} 
	\mathds{1}_{\{\sigma_x=\pm 1\}}
	\nu_{\beta,\pmb{h},G}(\sigma,\omega) \nonumber
	\\[0,2cm]
	=
	\sum_{\omega\in\{0,1\}^E}
	\Big[
	\sum_{\sigma\in\{-1,+1\}^V} 
	\mathds{1}_{\{\sigma_x=\pm1\}}
	\nu_{\beta,\pmb{h},G}(\sigma\vert \omega)
	\Big]
	\phi_{\pmb{p},\pmb{h},G}(\omega).\nonumber
	\end{multline}
	Using Corollary \ref{medida-conditional1}, the above
	expression can be rewritten as
	\[
	\sum_{\omega\in\{0,1\}^E}\left[
	\sum_{\sigma\in\{-1,+1\}^V}
	\mathds{1}_{\{\sigma_x=\pm 1\}}{\Delta(\sigma,\omega)
		\exp
		\big(
		\beta\sum_{i\in V}h_i\sigma_i
		\big)
		\over
		\prod_{\alpha=1}^{k(\omega,G)}2
		\cosh\left(\pmb{h}(K_\alpha)\right)}\right]
		\phi_{\pmb{p},\pmb{h},G}(\omega).
	\]
	Using now Lema \ref{lema anteprincipal}, we can see that the 
	above expression is equal to
	\begin{multline*}
		\sum_{\omega\in\{0,1\}^E}\left[
		{
		\exp\left( \pm \pmb{h}(K_t)\right)
		\prod_{\substack{\alpha=1 \\ \alpha\not=t}}^{k(\omega,G)}
		2\cosh\big(\pmb{h}(K_\alpha)\big)
		\over
		\prod_{\alpha=1}^{k(\omega,G)}
		2\cosh\big(\pmb{h}(K_\alpha)\big)
		}
		\right]
		\phi_{\pmb{p},\pmb{h},G}(\omega)
		\\[0,2cm]
		=
		{1\over 2}
		\sum_{\omega\in\{0,1\}^E}\left[
		1\pm\tanh\left(\pmb{h}(K_t)\right)
		\right]
		\phi_{\pmb{p},\pmb{h},G}(\omega).
	\end{multline*}
\end{proof}
\begin{corollary}
    Under the hypothesis of Theorem \ref{teorema de medida}, we 
    have that 
	\[
	\lambda_{\beta,\pmb{h},V}(\sigma_x)
	=
	\phi_{\pmb{p},\pmb{h},G}
	\big( 
		\tanh\left(\pmb{h}(K_t)\right)
	\big).
	\]
\end{corollary}
\begin{proof}
The proof follows directly from Theorem \ref{teorema de medida}.
\end{proof}
\section{General Potts models in external fields}
\label{secao-potts-geral}
In this last section we state two propositions establishing a 
graphical representation for the two-point function of the 
$q$-state Potts model with general external fields, defined in 
the Section \ref{sec-Potts-com-campo-cc-livre},  in terms of the 
connectivity of the random-cluster model introduced below.
The techniques employed to prove these results are similar to the 
ones we used in the previous section and therefore the proofs are 
omitted. 

Given a finite graph $G=(V,E)$, coupling constants 
$\pmb{J}=(J_{ij}\geqslant 0:\{i,j\}\in E)$   
and $\pmb{\hat{h}}$ a magnetic field as defined in the Section 
\ref{sec-Potts-com-campo-cc-livre}, for each 
$\omega\in \{0,1\}^{E}$ we define the finite-volume Gibbs measure 
of the (general) random-cluster model in external field by
\begin{equation}\label{def-rcm-q-geral}
\phi_{\pmb{p},\pmb{\hat{h}},q,G}(\omega)
=
{
	1
	\over
	\mathscr{Z}^{\mathrm{{RC}}}_{\pmb{p},\pmb{\hat{h}},q,G}
}
B_{\pmb{J},q}(\omega)
\prod_{\alpha =1}^{k({\omega,G})} 
\sum_{p=1}^q 
\exp
\big(
\beta\sum_{i\in K_\alpha}h_{i,p}
\big),
\end{equation}
where $K_\alpha$ is defined exactly 
as in the Section \ref{secao-def-modelo-RC}
and $B_{\pmb{J},q}(\omega)$ is similar to the Bernoulli 
factor of the Section \ref{secao-def-modelo-RC}
with exception that $p_{ij}=1-\exp(-q\beta J_{ij})$.

The Edwards-Sokal measure is generalized to
\begin{eqnarray}\label{ES1-q}
{\nu}_{\pmb{p},\pmb{\hat{h}},q,G}(\hat{\sigma}, \omega)
\equiv  
{1 \over \mathscr{Z}^{\mathrm{ ES}}_{\pmb{p},\pmb{\hat{h}},q,G}}
\displaystyle{
	B_{\pmb{J},q}(\omega)
	\Delta_q(\hat{\sigma}, \omega)
	\times
	\exp
	\big( \beta 
	\sum_{i\in V}\sum_{p=1}^q h_{i,p}\delta_{\hat{\sigma}_i,p}
	\big).
}
\end{eqnarray}
\begin{proposition}
	\label{Correlação/quase-conectividade-geral}
	Consider the Potts model with  Hamiltonian 
	given by
	\eqref{Potts-Geral}, densities 
	$p_{ij}\equiv1-\exp(-q\beta J_{ij})$ and 
	$q\in\{2,3, \ldots\}$ fixed. For any pair 
	of vertices $x,y\in V$ we have that 
	\begin{multline*}
		\tau_{q\beta,\pmb{\hat{h}},q,V}(x,y)
		=
		\Big(1-{1\over q}\Big)
		\phi_{\pmb{p},\pmb{\hat{h}},q,G}(x\leftrightarrow y)
		\\
		+
		\phi_{\pmb{p},\pmb{\hat{h}},q,G}
		\Big(
		\mathds{1}_{\{x\not\leftrightarrow y\}}
		\cdot
		\big\{
		H_{\pmb{\hat{h}}}(K_t,K_u)
		-
		{1\over q}
		\big\}
		\Big),
	\end{multline*}
	where the random variable 
	$H_{\pmb{\hat{h}}}(K_t,K_u)$ is given by  
	\[
	H_{\pmb{\hat{h}}}(K_t,K_u)
	\equiv
	{
		\sum_{r=1}^q 
		\exp
		\big(
		\beta\sum_{i\in K_t}h_{i,r}+ 
		\beta\sum_{i\in K_u}h_{i,r}
		\big)
		\over 
		\sum_{r=1}^q 
		\exp
		\big(
		\beta\sum_{i\in K_t}h_{i,r}
		\big)
		\cdot
		\sum_{r=1}^q 
		\exp
		\big(
		\beta\sum_{i\in K_u}h_{i,r}
		\big)
	}
	,
	\]
	with $K_t\equiv K_t(\omega)$ and 
	$K_u\equiv K_u(\omega)$ 
	being the disjoint connected components 
	containing the vertices 
	$x$ and $y$, respectively.
\end{proposition}
\begin{proof}
	We omit the proof of this proposition because it 
	is similar to the one given for Theorem 
	\ref{Correlação/quase-conectividade}.
\end{proof}
\begin{remark}
	Notice that in case $\pmb{\hat{h}}\equiv 0$, 
	we have for any $\omega\in \{0,1\}^{E}$ that  
	\[
	H_0(K_t,K_u)(\omega)= {q\over q^2}={1\over q},
	\]
	so Proposition \ref{Correlação/quase-conectividade-geral}
	gives us the following identity 
	\[
	\tau_{q\beta,0,q,V}(x,y)
	=
	\Big(1-{1\over q}\Big)
	\phi_{\pmb{p},0,q,G}(x\leftrightarrow y).
	\]
	This is also a very well know identity, as can be seen 
	in \cite{Grimmett2} Theorem 1.16, p. 11.
	Furthermore in case $q=2$ and  
	$h_{i,1}=-h_{i,2}={h_i}$ for all $i\in V$, 
	we have for any pair $x,y\in V$ that 
	\[
	H_{\pmb{h}}(K_t,K_u)
	=
	{1\over 2}
	\Big\{
	1+\tanh
	\big({\beta}\sum_{i\in K_t}h_i\big).
	\tanh\big({\beta}\sum_{i\in K_u}h_i\big)  
	\Big\}.
	\]
	In other words, Proposition 
	\ref{Correlação/quase-conectividade-geral} 
	generalizes Theorem \ref{Correlação/quase-conectividade}.
\end{remark}
\begin{proposition}\label{theor2}
	Let $G=(V,E)$ be a finite graph and  $x\in V$. 
	For each $m\in\{1,\ldots,q\}$ with $q\geqslant 1$, we have
	\begin{eqnarray*}
		\pi_{q\beta,\pmb{\hat{h}},q,V}(\hat{\sigma}_x=m)
		=
		\phi_{\pmb{p},\pmb{\hat{h}},q,G}
		\Big(
		{
			\exp
			(
			\beta\sum_{i\in K_t}h_{i,m}
			)
			\over 
			\sum_{p=1}^q
			\exp
			(
			\beta\sum_{i\in K_t}h_{i,p}
			)
		}
		\Big),
	\end{eqnarray*}
	where $K_t \equiv K_t(\omega)$ is the connected component of 
	$x$.
\end{proposition}
\noindent
	\textit{Sketch of the Proof.}
	To prove this theorem one needs to 
	compute the marginals of the Edwards-Sokal coupling given in
	\eqref{ES1-q}.
	The computation is similar to the one presented in the 
	previous sections. The next step is to prove the identity
	\[
	\pi_{q\beta,\pmb{\hat{h}},q,V}(\hat{\sigma}_x=m)
	=
	\sum_\omega
	\Big[\sum_{\hat{\sigma}} 
	\mathds{1}_{\{\hat{\sigma}_x=m\}}
	\nu_{\pmb{p},\pmb{\hat{h}},q,G}(\hat{\sigma}\vert \omega)
	\Big]
	\phi_{\pmb{p},\pmb{\hat{h}},q,G}(\omega)
	\]
	and then one proves that the rhs above is exactly
	\begin{eqnarray*}
		\sum_\omega
		\left[
		\sum_{\hat{\sigma}} 
		\mathds{1}_{\{\hat{\sigma}_x=\hat{\sigma}_y\}}
		{
			\Delta_q(\hat{\sigma},\omega)
			\exp
			\big(
			\beta\sum_{i\in V}
			\sum_{p=1}^q h_{i,p}\delta_{\hat{\sigma}_i,p}
			\big)
			\over 
			\prod_{\alpha=1}^{k(\omega,G)}
			\sum_{p=1}^q
			\exp 
			\big(
			\beta\sum_{i\in K_\alpha}h_{i,p}
			\big)
		}
		\right]
		\phi_{\pmb{p},\pmb{\hat{h}},q,G}(\omega).
	\end{eqnarray*}
	From this point, the result follows from the combinatorial 
	arguments presented before.
\qed
%
%%%%%%%%%%%%%%%%%%%%%%%%%%%%%%%%%%%%%%%%%%%%%%%%%%%%%%%%%%%%%%%
\part{General boundary conditions}
%%%%%%%%%%%%%%%%%%%%%%%%%%%%%%%%%%%%%%%%%%%%%%%%%%%%%%%%%%%%%%%
%
%
%
\section{The general random-cluster model}
In this section we define the so called general random-cluster 
model on the lattice $\mathbb{L}=(\mathbb{V},\mathbb{E})$ 
(this terminology, GRC model, comes from \cite{BBCK00})
with inhomogeneous magnetic field  of the form 
$
\pmb{\hat{h}}
\equiv
(h_{i,p}:i\in \mathbb{V};\ p=1,\ldots,q)
\in 
\mathbb{R}^{\mathbb{V}}\times \cdots\times \mathbb{R}^{\mathbb{V}}
$
and boundary conditions.

The Bernoulli factors introduced before will be replaced in this 
section by  (abusing notation)
\begin{eqnarray}\label{nova def. Bernoulli}
B_{\pmb{J}}(\omega)\equiv \prod_{\{i,j\}:\omega_{ij}=1}
r_{ij},
\end{eqnarray}
where $\pmb{J}=(J_{ij}\geqslant 0:\{i,j\}\in E)$,
$r_{ij}\equiv \exp(q\beta J_{ij})- 1$ and $q\in\mathbb{Z}^+$ 
fixed.
Although $r_{ij}\geqslant 0$, in general, they are not bounded
by one, but mind that the random-cluster measure obtained with 
such ``Bernoulli factors'' is the same one gets when considering 
the old Bernoulli factors, since the weights in both cases are 
related by an overall normalization factor that cancels out 
because of the partition function.

Fix a random subgraph $G=(V,E)$ on the lattice $\mathbb{L}$, let 
$\partial E
=
\{ e\in\mathbb{E}:\ e\cap V\neq \emptyset\ 
\text{and}\ e\cap \partial V\neq \emptyset \}.
$
We denote by $\mathbb{B}_0(V)$
the set of all edges $\{x,y\}\in\mathbb{E}$ so that 
$\{x,y\}\subset V$. 
With this definition we have $\mathbb{B}_0(V)=E$.
We use the notation $\mathbb{B}(V)$ to denote the set of all 
edges with at least one vertex in $V$. Note that 
$\mathbb{B}(V)=E\cup\partial E$. For
any $\widetilde{E}\subset\mathbb{B}_0(\mathbb{V})$, we define 
$\mathbb{V}(\widetilde{E})$ as the set of sites which 
belong to at least one edge in $\widetilde{E}$.
\paragraph{GRC model with general boundary condition.}
Fix a finite subgraph $G=(V,E)$ of the lattice $\mathbb{L}$. 
For each $i\in\mathbb{V}$ we define
$   
h_{i,\mathrm{ max}}
\equiv 
\max\{h_{i,p}: p=1,\ldots, q\}.
$
If
$\omega\in \{0,1\}^{\mathbb{E}}$ and $C(\omega)$ denotes a 
generic connected component on $(\mathbb{V},\eta(\omega))$,
the GRC measure with general boundary condition is obtained by 
normalizing the followings weights
\begin{eqnarray}\label{definição de general boundary condition}
\mathcal{W}^{\mathrm{ GRC}}_{E}(\omega_E|\omega_{E^c})
\equiv
B_{\pmb{J}}(\omega)
\!\!\!\!\!\!\!\!
\prod_{
\substack{C(\omega):\\ \mathbb{V}(C(\omega))\cap V\not=\emptyset}} 
\sum_{p=1}^q 
q_p \ 
\exp
\big(
{-\beta\!\!\!\underset{i\in C(\omega)}{\sum}
(h_{i,\mathrm{ max}}-h_{i,p})}
\big),
\end{eqnarray}
where $\{q_p:p=1,\ldots, q\}$ are  positive constants, 
$B_{\pmb{J}}(\omega)$ is given by \eqref{nova def. Bernoulli} 
and the product runs over all the connected components 
$C(\omega)$ of the graph $(\mathbb{V},\eta(\omega))$.
In the above expression we are using the convention 
$e^{-\infty}=0$. This measure is denoted by
$
\phi^{\mathrm{ GRC}}_{E}.
$
\paragraph{GRC model with free boundary condition.}
Let $G=(V,E)$ be a finite graph and $\omega\in\{0,1\}^E$
a configuration. If $C(\omega)$ denotes a generic 
connected component on $(V,\eta(\omega))$, we define 
\begin{eqnarray*}
\Theta_{V,\mathrm{ free}}(C(\omega))
\equiv 
\sum_{p=1}^q 
q_p\ 
{\textstyle
\exp
\big(
\beta\underset{i\in C(\omega)}{\sum}h_{i,p}
\big).
}
\end{eqnarray*}
The GRC measure with free boundary condition is obtained by 
normalizing the weights 
\begin{eqnarray} \label{free peso}
\mathcal{W}^{\mathrm{ GRC}}_{V,\mathrm{ free}}(\omega)
\equiv
B_{\pmb{J}}(\omega)
\prod_{C(\omega)} \Theta_{V,\mathrm{ free}}(C(\omega)),
\end{eqnarray}
where $B_{\pmb{J}}(\omega)$ is given by 
\eqref{nova def. Bernoulli} and the product runs over all the 
connected components $C(\omega)$ of the graph $(V,\eta(\omega))$.
This measure is denoted by
$\phi^{\mathrm{ GRC}}_{V,\mathrm{ free}}$ 
and for each $\omega\in\{0,1\}^E$ it satisfies
$
\phi^{\mathrm{ GRC}}_{V,\mathrm{ free}}(\omega)
\propto 
\mathcal{W}^{\mathrm{ GRC}}_{V,\mathrm{ free}}(\omega),
$
where the proportionality constant is 
exactly the (inverse of the) partition function of the GRC model.
\paragraph{GRC model with wired boundary condition.}
Fix $\mathrm{m}\in \{1,\ldots, q\}$
and a finite subgraph $G=(V,E)$ of the lattice $\mathbb{L}$. 
If for each $\omega\in \{0,1\}^{E\cup\partial E}$.
$C(\omega)$ denotes a connected 
component on $(V\cup\partial V,\eta(\omega))$, then we define 
\begin{eqnarray*}
\Theta_{V,\mathrm{ m}}(C(\omega))
\equiv 
\begin{cases}
\Theta_{V,\mathrm{ free}}(C(\omega)),
& \text{if}\ C(\omega)\cap \partial V=\emptyset
\\[0,2cm]
\exp
{\textstyle
\big(\beta\underset{i\in C(\omega)}{\sum}h_{i,\mathrm{m}}\big),
}
&  \text{otherwise}.
\end{cases}
\end{eqnarray*}
Similarly, the GRC measure with m-wired 
boundary condition is obtained by normalizing
the weights
\begin{eqnarray}\label{def.rcm-wired condition}
\mathcal{W}^{\mathrm{ GRC}}_{V,\mathrm{ m}}(\omega)
\equiv
B_{\pmb{J}}(\omega)\prod_{C(\omega)} 
\Theta_{V,\mathrm{ m}}(C(\omega)),
\end{eqnarray}
where the product runs over all the connected components 
$C(\omega)$ of the graph 
$(V\cup\partial V,\eta(\omega))$. This measure is denoted by
$
\phi^{\mathrm{ GRC}}_{V,\mathrm{ m}}.
$
\begin{remark}
One can easily see that when 
$E\equiv\mathbb{B}(\Lambda)$, $\Lambda\subset \mathbb{V}$ finite, 
in \eqref{definição de general boundary condition} 
\begin{align*}
\mathcal{W}^{\mathrm{ GRC}}_{\mathbb{B}
(\Lambda)}(\omega_{\mathbb{B}(\Lambda)}
|
\omega_{\mathbb{B}(\Lambda)^c}^{(1)})
&=
e^{-\beta\sum_{i\in \Lambda}h_{i,\mathrm{ max}}} q_{\mathrm{ m}}
\cdot
B_{\pmb{J}}(\omega)\prod_{C(\omega)}
\Theta_{\Lambda,\mathrm{ m}}(C(\omega))
\\[0,3cm]
&\stackrel{\eqref{def.rcm-wired condition}}{=}
e^{-\beta\sum_{i\in \Lambda}h_{i,\mathrm{ max}}} q_{\mathrm{ m}}
\cdot
\mathcal{W}^{\mathrm{ GRC}}_{\Lambda,\mathrm{ m}}(\omega),
\end{align*}
where $\omega^{(i)}$ is the configuration with 
$\omega_{e}^{(i)}=i$ for all
$e\in\mathbb{B}_0(\mathbb{V})$ $(i=0,1)$. Therefore
\[
\phi^{\mathrm{ GRC}}_{\mathbb{B}(\Lambda)}
(\omega_{\mathbb{B}(\Lambda)}
|
\omega_{\mathbb{B}(\Lambda)^c}^{(1)})
=
\phi^{\mathrm{ GRC}}_{\Lambda,\mathrm{ m}}(\omega).
\]
Similarly we obtain in 
\eqref{definição de general boundary condition} with 
$E\equiv\mathbb{B}_0(\Lambda)$
\begin{align*}
\mathcal{W}^{\mathrm{ GRC}}_{\mathbb{B}_0
(\Lambda)}(\omega_{\mathbb{B}_0(\Lambda)}
|
\omega_{\mathbb{B}_0(\Lambda)^c}^{(0)})
&=
e^{-\beta\sum_{i\in \Lambda}h_{i,\mathrm{ max}}}
\cdot
B_{\pmb{J}}(\omega)\prod_{C(\omega)}
\Theta_{\Lambda,\mathrm{ free}}(C(\omega))
\\[0,3cm]
&\stackrel{\eqref{free peso}}{=}
e^{-\beta\sum_{i\in \Lambda}h_{i,\mathrm{ max}}}
\cdot
\mathcal{W}^{\mathrm{ GRC}}_{\Lambda,\mathrm{ free}}(\omega),
\end{align*}
then
\[
\phi^{\mathrm{ GRC}}_{\mathbb{B}_0(\Lambda)}
(\omega_{\mathbb{B}_0(\Lambda)}|\omega_{\mathbb{B}_0(\Lambda)^c}^{(0)})
=
\phi^{\mathrm{ GRC}}_{\Lambda,\mathrm{ free}}(\omega).
\]
\end{remark}
\subsection{The FKG inequality}
Throughout this section we assume that 
$\{q_p:p=1,\ldots,q\}$ introduced in 
\eqref{definição de general boundary condition} 
and the magnetic field $\pmb{\hat{h}}$ 
satisfy
\begin{eqnarray}\label{qmmaiorque1}
\sum_{p\in 
{\cap}_{i\in\mathbb{V}}\mathcal{Q}_{i,\mathrm{ max}}
(\pmb{\hat{h}})}q_p \geqslant 1,
\end{eqnarray}
where
$
\mathcal{Q}_{i,\mathrm{ max}}(\pmb{\hat{h}})
\equiv \big\{p\in\{1,\ldots,q\}:h_{i,p}=h_{i,\mathrm{ max}}\big\}.
$
We consider as usual the partial order on $\{0,1\}^\mathbb{E}$ 
where  
$
\omega\preceq \tilde{\omega} 
\ \Longleftrightarrow \ 
\omega_e\leqslant \tilde{\omega}_e, \ \forall \ e\in \mathbb{E}.
$
We also use the standard notations 
$\omega_1\vee\omega_2$ and $\omega_1\wedge\omega_2$ for 
$
(\omega_1\vee\omega_2)_{e}=\max\{\omega_1(e),\omega_2(e)\}
\ \text{and}\ 
(\omega_1\wedge\omega_2)_{e}=\min\{\omega_1(e),\omega_2(e)\}
$
with $ e\in \mathbb{E}$, respectively.
\begin{definition}[FKG property] 
	Let $(\Omega,\preceq)$ be a partially ordered space. 
	A measure $\mu$ over $\Omega$ said to have the  
	$\mathrm{FKG}$ property if 
	\[
	\mu(fg)\geqslant \mu(f)\mu(g) 
	\]
	for any increasing (with respect to $\preceq$) 
	measurable functions $f,g:\Omega\to \mathbb{R}$. 
	Furthermore, if $\Omega$ is a cartesian product 
	$\Omega=\prod_{e\in B}\Omega_e$, with $|\Omega_{e}|<\infty$,
	 then $\mu$ is said to have 
	the {\bf strong FKG property}, if 
	$\mu(\cdot\vert A)$ has the $\mathrm{FKG}$ property for each 
	cylinder event  
	$
	A
	=
	\{
	\omega\in\Omega:\omega_e=\alpha_e, \ \forall e\in \tilde{B}
	\}
	$,
	where $\tilde{B}\subset{B}$ is finite and 
	$\alpha_e\in\Omega_e$ for all $e\in\tilde{B}$.
\end{definition}
\begin{remark} 
	If $\mathrm{m},\mathrm{\widetilde{m}}\in$ 
	${\cap}_{i\in \mathbb{V}}
	\mathcal{Q}_{i,\mathrm{ max}}(\pmb{\hat{h}})$, 
	then
	$
	\Theta_{V, \mathrm{\widetilde{m}}}(C)
	=
	\Theta_{V,\mathrm{m}}(C)
	$
	and therefore
	$\phi^{\mathrm{ GRC}}_{V,\mathrm{\widetilde{m}}}=
	\phi^{\mathrm{ GRC}}_{V,\mathrm{m}}$.
	This measure is denoted by 
	$\phi^{\mathrm{ GRC}}_{V,\mathrm{ max}}$.
\end{remark}
\begin{theorem}[Strong FKG property] \label{Strong FKG Property}
	Let $q\in \mathbb{Z}^+$, $\beta\geqslant 0$, 
	$
	\pmb{J}
	=
	(J_{ij}:\{i,j\}\in \mathbb{E})$ $\in [0,\infty)^\mathbb{E}
	$, 
	$
	\pmb{\hat{h}}
	=
	(
	h_{i,p}\in\mathbb{R}: i\in\mathbb{V}, 1\leqslant p\leqslant q
	)
	$ 
	and $\{q_p: p=1,\ldots,q\}$  satisfying \eqref{qmmaiorque1}. 
	Then for any finite subgraph $G=(V,E)$ of $\mathbb{L}$, the 
	measures 
	$\phi^{\mathrm{ GRC}}_{V,\mathrm{ free}}$ and
	$\phi^{\mathrm{ GRC}}_{V,\mathrm{ max}}$ 
	have the strong $\mathrm{FKG}$ property.
\end{theorem}
\begin{proof}
	For simplicity  
	we assume that the magnetic field
	we are dealing with satisfies 
	the following inequalities
	\begin{eqnarray}\label{crescente-campos}
	h_{i,1}\leqslant h_{i,2}\leqslant \ldots \leqslant h_{i,q},
	\quad
	\forall \ i\in\mathbb{V}.
	\end{eqnarray}

	The FKG lattice condition for the 
	$\phi^{\mathrm{ GRC}}_{V,\mathrm{ free}}$ is equivalent to
	\[
	\mathcal{W}^{\mathrm{ GRC}}_{V,\mathrm{ free}}
	(\omega^{(1)}\vee \omega^{(2)})
	\mathcal{W}^{\mathrm{ GRC}}_{V,\mathrm{ free}}
	(\omega^{(1)}\wedge \omega^{(2)})
	\geqslant 
	\mathcal{W}^{\mathrm{ GRC}}_{V,\mathrm{ free}}(\omega^{(1)})
	\mathcal{W}^{\mathrm{ GRC}}_{V,\mathrm{ free}}(\omega^{(2)}),
	\]
	where $\omega^{(1)}$ and $\omega^{(2)}$ are arbitrary 
	configurations.
	Similarly for $\phi^{\mathrm{ GRC}}_{V,\mathrm{ max}}$. 
	It is well known that such condition 
	implies the strong FKG property, see for example 		
	Theorem $2.19$, p. $25$ in \cite{Grimmett2}. 
	By defining 
	\[
	\mathcal{R}(\xi,\omega)
	\equiv
	{
		\mathcal{W}^{\mathrm{ GRC}}_{V,\mathrm{ free}}
		(\xi\vee \omega)
		\over 
		\mathcal{W}^{\mathrm{ GRC}}_{V,\mathrm{ free}}(\xi)
	},
	\]
	one can see that the lattice condition holds if
	\begin{eqnarray}\label{min}
	\mathcal{R}(\omega^{(1)},\omega^{(2)})
	\geqslant 
	\mathcal{R}(\omega^{(1)}\wedge\omega^{(2)},\omega^{(2)}).
	\end{eqnarray}
	For a fixed configuration $\omega$, 
	we chose an arbitrary order for $\eta(\omega)$ and 
	represent these open edges 
	as $(e_1,\ldots,e_{\vert \eta(\omega)\vert})$.
	So for any configuration $\xi\in\{0,1\}^{E}$ 
	we have that 
	\[
	\mathcal{R}(\xi,\omega)=\prod_{k=1}^{\vert\eta(\omega)\vert}
	\mathcal{R}(\xi\vee\omega^{(e_1)}\vee\cdots 
	\vee 
	\omega^{(e_{k-1})},\omega^{(e_{k})}),
	\]
	where $ (\omega^{(e)})_{e'}\equiv \delta_{e,e'}$.
	Therefore it is enough to prove \eqref{min} 
	for configurations $\xi$, $\omega^{(1)}$ and $\omega^{(2)}$ 
	such that $\xi$ has at least two zero coordinates or 
	at most one zero and 	 
	$\omega^{(1)}\equiv \xi\vee\omega^{(b)}$
	and $\omega^{(2)}\equiv \xi\vee\omega^{(b')}$.
	Let us begin assuming that 
	$\xi$ has at least two zero coordinates
	and   
	\[
	\xi\equiv (*,\ldots,*,\aunderbrace[D]{0}_{ b-\mathrm{th}},
	*,\ldots,*,\aunderbrace[D]{0}_{b'\mathrm{-th}},*,\ldots,*),
	\]
	where $b,b'\in E\cup \partial E $, $b\not=b'$ and the stars 
	indicate generic elements in $\{0,1\}$ (not necessarily 
	equal). If we define  	  
	\[
	\xi^b\equiv (*,\ldots,*,\aunderbrace[D]{1}_{b\mathrm{-th}},
	*,\ldots,*,\aunderbrace[D]{0}_{ b'\mathrm{-th}},*,\ldots,*)
	\]
	and
	\[
	\xi^{b'}\equiv (*,\ldots,*,\aunderbrace[D]{0}_{b\mathrm{-th}},
	*,\ldots,*,\aunderbrace[D]{1}_{b'\mathrm{-th}},*,\ldots,*),
	\]
	then we have that 
	$
	\omega^{(1)}=\xi\vee\omega^{(b)}=\xi^b, 
	\omega^{(2)}=\xi\vee\omega^{(b')}=\xi^{b'}
	\text{and}\allowbreak\ 
	\omega^{(1)}\wedge\omega^{(2)}=\xi
	.
	$
	So in order to prove \eqref{min} it is enough to prove that 
	\begin{eqnarray}\label{por provar}
	\mathcal{R}(\xi^b,\xi^{b'})\geqslant\mathcal{R}(\xi,\xi^{b'}),
	\quad \mbox{with} \ b\not = b'.
	\end{eqnarray}
	Now we concentrate on proving \eqref{por provar}.
	To do this we first observe that if 
	$\prod_{\{i,j\}:\xi_{ij}=1}r_{ij}=k$, then
	\[
	\prod_{\{i,j\}:(\xi^b\vee\xi^{b'})_{ij}=1}
	\!\!\!\!\!\!\!\! r_{ij}
	=
	r_b r_{b'} k,
	\quad 
	\prod_{\{i,j\}:\xi^b_{ij}=1}
	\!\!\!\! r_{ij}
	=
	r_b k
	\quad 
	\text{and}
	\quad
	\prod_{\{i,j\}:(\xi\vee\xi^{b'})_{ij}=1}
	\!\!\!\!\!\!\!\! r_{ij}
	=
	r_{b'} k.
	\]
	So it follows from the definition \eqref{nova def. Bernoulli} 
	that
	\[
	{B_{\pmb{J}}(\xi^b\vee\xi^{b'}) \over B_{\pmb{J}}(\xi^b)}
	= 
	{
		\prod_{\{i,j\}:(\xi^b\vee\xi^{b'})_{ij}=1}r_{ij}
		\over
		\prod_{\{i,j\}:\xi^b_{ij}=1}r_{ij}
	}
	=
	r_{b'}
	=
	{
		\prod_{\{i,j\}:(\xi\vee\xi^{b'})_{ij}=1}r_{ij}
		\over
		\prod_{\{i,j\}:\xi_{ij}=1}r_{ij}
	}
	=
	{B_{\pmb{J}}(\xi\vee\xi^{b'}) \over B_{\pmb{J}}(\xi)}.
	\]
	Because of the above observation and the definitions of 
	$\mathcal{W}^{\mathrm{ GRC}}_{V,\mathrm{ free}}$
	and
	$\mathcal{W}^{\mathrm{ GRC}}_{V,\mathrm{ m}}$,
	the proof of \eqref{por provar} reduces to 
	\begin{eqnarray}\label{Desigualdade mt importante}
	{
		\Theta_{V,\mathrm{ \#}}(C(\xi^b\vee\xi^{b'}))
		\over 
		\Theta_{V,\mathrm{ \#}}(C(\xi^b))
	}
	\geqslant 
	{
		\Theta_{V,\mathrm{ \#}}(C(\xi\vee\xi^{b'}))
		\over 
		\Theta_{V,\mathrm{ \#}}(C(\xi))
	},
	\end{eqnarray}
	where  $\mathrm{\#}$ stands for ``free''or ``m''.
\noindent
\\\\
{\sc Free boundary condition case.}
\\
We broke the proof of \eqref{Desigualdade mt importante} 
in several cases. Let $A_1,A_2,B_1$ and $B_2$ be connected 
components of $(V,\eta(\xi))$ and consider the cases 
showed in the picture below 
	\begin{center} 
		\begin{minipage}{\linewidth}
			\makebox[\linewidth]{%
				\includegraphics[scale=0.4,keepaspectratio=true]
				{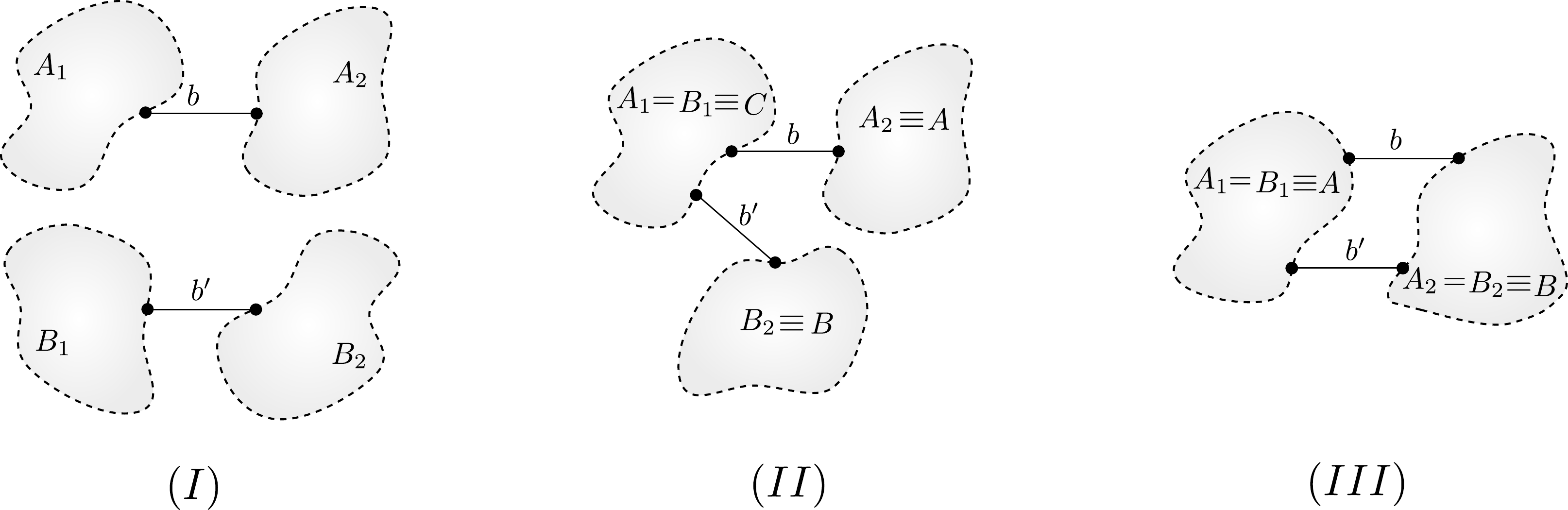}}
			%\captionof{figure}{}
			\label{q15}
		\end{minipage}
	\end{center}
	The case
	$(I)$ represents that the end vertices of $b$ belong
	to $A_1$ and $A_2$ and the end vertices of $b'$ 
	belong to $B_1$ and $B_2$. In this case 
	the left and right sides of 
	\eqref{Desigualdade mt importante}
	are equal, since 
	\begin{multline*}
	{
		\Theta_{V,\mathrm{ free}}(A_1\cup A_2)
		\Theta_{V,\mathrm{ free}}(B_1\cup B_2)
		\over 
		\Theta_{V,\mathrm{ free}}(A_1\cup A_2)
		\Theta_{V,\mathrm{ free}}(B_1)
		\Theta_{V,\mathrm{ free}}(B_2)
	}
	\\[0.2cm]
	=
	{
		\Theta_{V,\mathrm{ free}}(B_1\cup B_2)
		\Theta_{V,\mathrm{ free}}(A_1)
		\Theta_{V,\mathrm{ free}}(A_2)
		\over 
		\Theta_{V,\mathrm{ free}}(A_1)
		\Theta_{V,\mathrm{ free}}(A_2)
		\Theta_{V,\mathrm{ free}}(B_1)
		\Theta_{V,\mathrm{ free}}(B_2)
	}.
	\end{multline*}For the case $(II)$, we should prove that 
	\[
	{
		\Theta_{V,\mathrm{ free}}(A\cup B\cup C)
		\over 
		\Theta_{V,\mathrm{ free}}(C\cup A)
		\Theta_{V,\mathrm{ free}}(B)
	}
	\geqslant 
	{
		\Theta_{V,\mathrm{ free}}(C\cup B)
		\Theta_{V,\mathrm{ free}}(A)
		\over 
		\Theta_{V,\mathrm{ free}}(A)
		\Theta_{V,\mathrm{ free}}(B)
		\Theta_{V,\mathrm{ free}}(C)
	},
	\]
	which is equivalent to 
	\begin{eqnarray}\label{seq,nerv}
	\Theta_{V,\mathrm{ free}}(C)
	\Theta_{V,\mathrm{ free}}(A\cup B\cup C)
	\geqslant 
	\Theta_{V,\mathrm{ free}}(C\cup A)
	\Theta_{V,\mathrm{ free}}(C\cup B).
	\end{eqnarray}
	To help us prove inequality \eqref{seq,nerv}, 
	we define for each 
	$m\in \{1,\ldots,q\}$ the following numbers
	\begin{eqnarray*}
		a_m
		\equiv
		\exp \big(\beta{\sum}_{i\in A} \!h_{i,m}\big), 
		\,
		b_m
		\equiv
		\exp \big(\beta{\sum}_{i\in B} \!h_{i,m}\big)
		\ 
		\text{and} \ 
		c_m
		\equiv
		\exp \big(\beta{\sum}_{i\in C} \!h_{i,m}\big).
	\end{eqnarray*}

	The hypothesis 	\eqref{crescente-campos}
	implies immediately that $(a_m)$ and $(b_m)$
	are non-decreasing in $m$.
	Using this notation, \eqref{seq,nerv} reads
	\begin{eqnarray}\label{serie}
	\sum_{m=1}^q q_m c_m
	\sum_{m'=1}^q q_{m'}a_{m'}b_{m'}c_{m'}
	\geqslant 
	\sum_{m=1}^q q_m a_mc_m
	\sum_{m'=1}^q q_{m'}b_{m'}c_{m'}.
	\end{eqnarray}
	Both sides of the above inequality can be written using 
	a bilinear form 
	\[
	\varphi(a,b)\equiv\sum_{m,m'=1}^q r_{m,m'}a_{m'}b_{m},
	\]
	where $r_{m,m'}\equiv q_mc_mq_{m'}c_{m'}$ , 
	$a\equiv(a_1,\ldots,a_q)$
	and $b\equiv(b_1,\ldots,b_q)$.
	Note that $\varphi$ is a symmetric bilinear form 
	and \eqref{serie} can be written as
	\begin{eqnarray}\label{var}
	\varphi(1,c)\geqslant \varphi(a,b),\quad 
	\text{where}\ c\equiv (a_1b_1,\ldots,a_qb_q).
	\end{eqnarray}
	Therefore it remains to prove \eqref{var}, which clearly
	holds since 
	\[
	r_{m,m'}(a_{m'}-a_m)(b_{m'}-b_m)\geqslant 0
	\ \ \Longleftrightarrow \ \ 
	\varphi(1,c)-\varphi(a,b)-\varphi(b,a)+\varphi(c,1)\geqslant 0.
	\]

	We proceed with \eqref{Desigualdade mt importante} for the 
	case $(III)$. Now we have to prove that 
	\begin{equation*}
	1=
	{
		\Theta_{V,\mathrm{ free}}(A\cup B)
		\over 
		\Theta_{V,\mathrm{ free}}(A\cup B)
	}
	\geqslant
	{ 
		\Theta_{V,\mathrm{ free}}(A\cup B)
		\over 
		\Theta_{V,\mathrm{ free}}(A)
		\Theta_{V,\mathrm{ free}}(B)
	},
	\end{equation*}
	in other words 
	$	
	\Theta_{V,\mathrm{ free}}(A)
	\Theta_{V,\mathrm{ free}}(B)
	\geqslant 
	\Theta_{V,\mathrm{ free}}(A\cup B)
	$,
	or equivalently
	\[ \textstyle
	\sum_{m=1}^qq_ma_m\sum_{m'=1}^qq_{m'}b_{m'}
	\geqslant 
	\sum_{m=1}^qq_ma_mb_m.
	\]
	This last inequality is actually true since 
	\[ 
	\sum_{m=1}^qq_ma_m \sum_{m'=1}^qq_{m'}b_{m'}
	\geqslant 
	\sum_{m=1}^qq_ma_m b_{\text{max}}
	\sum_{m'\in {\cap}_{i\in\mathbb{V}}
	\mathcal{Q}_{i,\mathrm{ max}}(\pmb{\hat{h}})}^qq_{m'}
	\geqslant 
	\sum_{m=1}^qq_ma_mb_m.
	\]

	For the cases where the end vertices of $b$ or $b'$ are 
	contained 
	in the same connected component, the inequality is trivial.
\noindent
\\\\
{{\sc Max wired boundary condition case.}}
\\
Suppose that 
$
\mathrm{\widetilde{m}}\in 
{\cap}_{{i\in\mathbb{V}}}\mathcal{Q}_{i,\mathrm{max}}
(\pmb{\hat{h}})
$.
To prove the inequality \eqref{Desigualdade mt importante} 
we have to analyze again the three cases above.
For case $(I)$,	analogously to the free boundary 
condition case, we have
	\[
	{
		\Theta_{V,\mathrm{ \widetilde{m}}}(A_1\cup A_2)
		\Theta_{V,\mathrm{ \widetilde{m}}}(B_1\cup B_2)
		\over 
		\Theta_{V,\mathrm{ \widetilde{m}}}(A_1\cup A_2)
		\Theta_{V,\mathrm{ \widetilde{m}}}(B_1)
		\Theta_{V,\mathrm{ \widetilde{m}}}(B_2)
	}
	=
	{
		\Theta_{V,\mathrm{ \widetilde{m}}}(B_1\cup B_2)
		\Theta_{V,\mathrm{ \widetilde{m}}}(A_1)
		\Theta_{\mathrm{ \widetilde{m}}}(A_2)
		\over 
		\Theta_{V,\mathrm{ \widetilde{m}}}(A_1)
		\Theta_{V,\mathrm{ \widetilde{m}}}(A_2)
		\Theta_{V,\mathrm{ \widetilde{m}}}(B_1)
		\Theta_{V,\mathrm{ \widetilde{m}}}(B_2)
	},
	\]
	independently on whether the components $A_1,A_2,B_1$ and 
	$B_2$  and the possible combinations among them intersect 
	$V^c$.
	
	For the case $(II)$ and all the configurations sketched 
	on the figure below 
	\begin{center} 
		\begin{minipage}{\linewidth}
			\makebox[\linewidth]{%
				\includegraphics[scale=0.32,keepaspectratio=true]
				{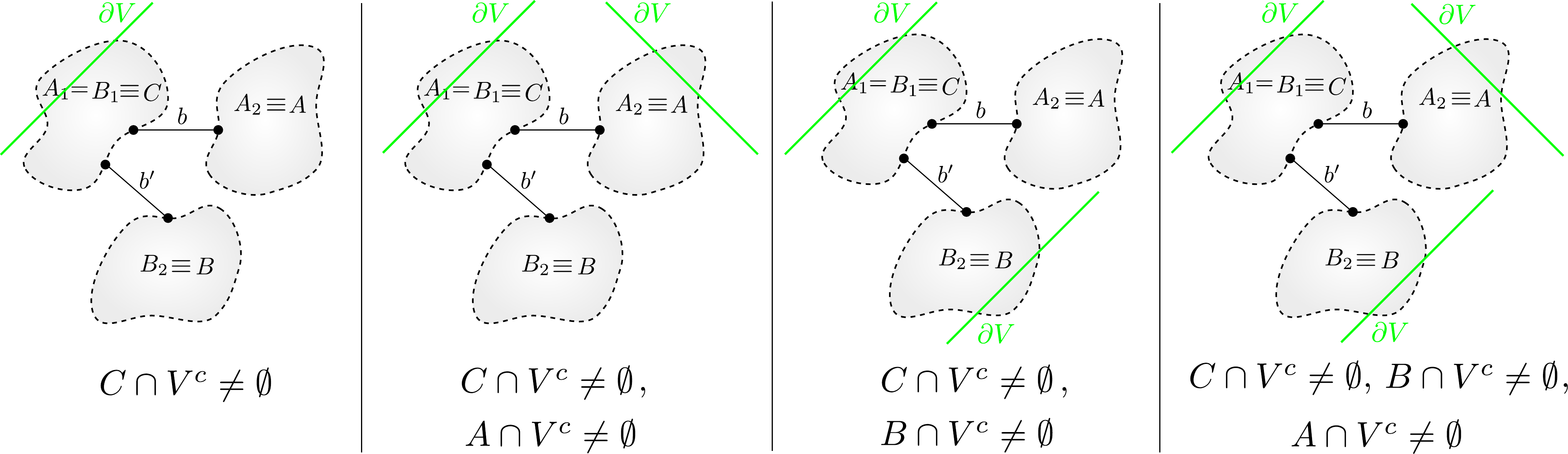}}
			%\captionof{figure}{}
		\end{minipage}
	\end{center}
	we have from the definition \eqref{def.rcm-wired condition} 
	of $\Theta_{V,\mathrm{ \widetilde{m}}}$ that the following equality 
	holds
	\[
	c_{\mathrm{\widetilde{m}}}a_{\mathrm{\widetilde{m}}}
	b_{\mathrm{\widetilde{m}}}c_{\mathrm{\widetilde{m}}}
	=
	a_{\mathrm{\widetilde{m}}}c_{\mathrm{\widetilde{m}}}
	b_{\mathrm{\widetilde{m}}}c_{\mathrm{\widetilde{m}}}.
	\]
	For the following configurations that also appears in the 
	case $(II)$ : 
	\begin{center} 
		\begin{minipage}{\linewidth}
			\makebox[\linewidth]{%
				\includegraphics[scale=0.42,keepaspectratio=true]
				{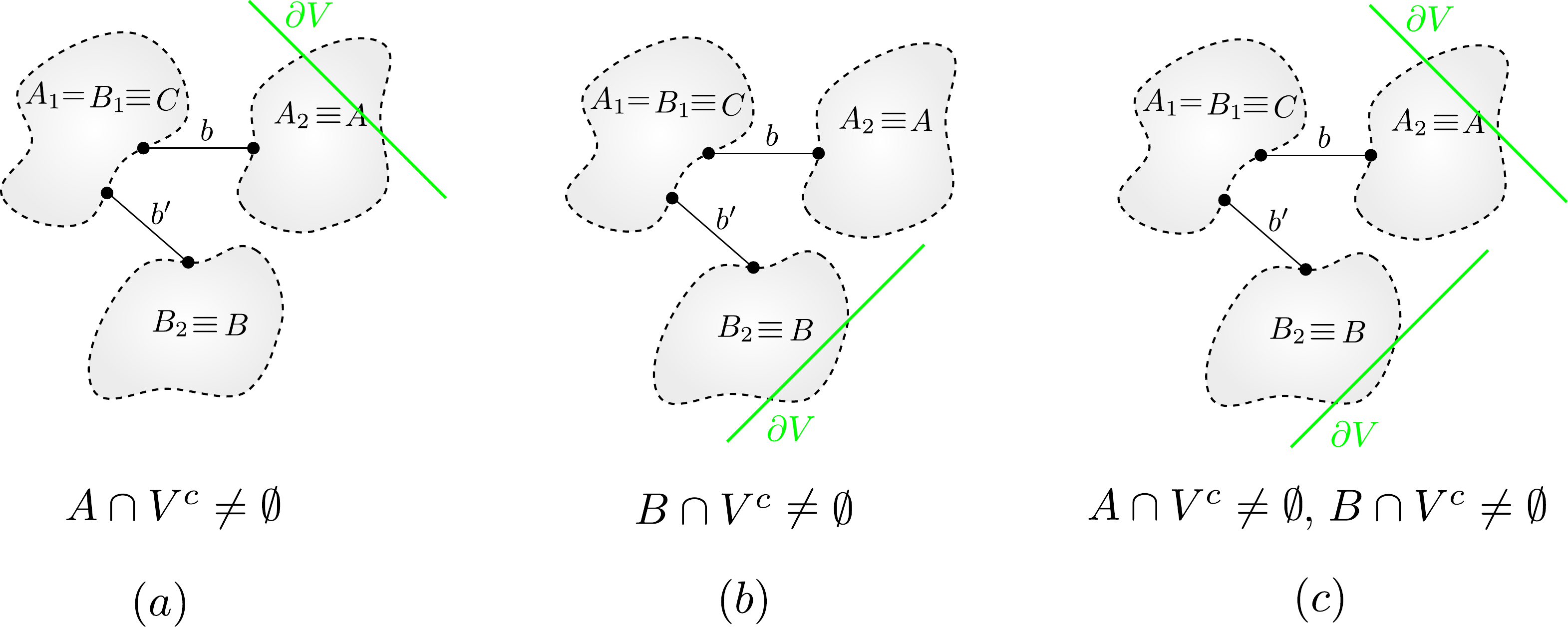}}
			%\captionof{figure}{}
		\end{minipage}
	\end{center}
	For $(a)$, the inequality \eqref{Desigualdade mt importante}
	comes from 
	\[
	{\textstyle
	\left(\sum_{m=1}^q q_m c_m\right)
	a_{\mathrm{\widetilde{m}}}b_{\mathrm{\widetilde{m}}}
	c_{\mathrm{\widetilde{m}}}
	\geqslant 
	a_{\mathrm{\widetilde{m}}}c_{\mathrm{\widetilde{m}}}
	\left(\sum_{m'=1}^q q_{m'}b_{m'}c_{m'}\right),
	}
	\]
	which is always valid since we have that  
	$b_{\mathrm{\widetilde{m}}}\geqslant b_{m'}, \ 
	\forall m'=1,\ldots,q.$
	In $(b)$ inequality \eqref{Desigualdade mt importante},
	comes from 
	\[
	{\textstyle
	\left(\sum_{m=1}^q q_m c_m \right)
	a_\mathrm{\widetilde{m}}b_\mathrm{\widetilde{m}}
	c_\mathrm{\widetilde{m}}
	\geqslant
	\left(\sum_{m=1}^q q_m a_m c_m\right)
	b_\mathrm{\widetilde{m}} c_\mathrm{\widetilde{m}},
	}
	\]
	which is also true because
	$a_\mathrm{\widetilde{m}}\geqslant a_{m}, \ 
	\forall m=1,\ldots,q.$ 
	Finally, in $(c)$ inequality 
	\eqref{Desigualdade mt importante},
	is a consequence of
	\[
	{\textstyle
	\left(\sum_{m=1}^q q_mc_m\right)
	a_\mathrm{\widetilde{m}} b_\mathrm{\widetilde{m}} 
	c_\mathrm{\widetilde{m}}
	\geqslant
	a_\mathrm{\widetilde{m}} c_\mathrm{\widetilde{m}}
	b_\mathrm{\widetilde{m}} c_\mathrm{\widetilde{m}}
	}
	\]
	and the validity of this inequality is ensured by 
	\[
	{\textstyle
	\sum_{m=1}^q q_mc_m
	\geqslant
	\underset{   
	m\in\cap_{i\in \mathbb{V}}\mathcal{Q}_{i,\mathrm{ max}}
	(\pmb{\hat{h}})   
	}{\sum}
	q_mc_{\mathrm{ max}}
	\geqslant
	c_{\mathrm{ max}}=c_\mathrm{\widetilde{m}},
	}
	\]
	which follows from \eqref{qmmaiorque1}.	
	
	Now we consider the case $(III)$, by splitting  its analysis
	in the following sub-cases
	\begin{center} 
		\begin{minipage}{\linewidth}
			\makebox[\linewidth]{%
				\includegraphics[scale=0.39,keepaspectratio=true]
				{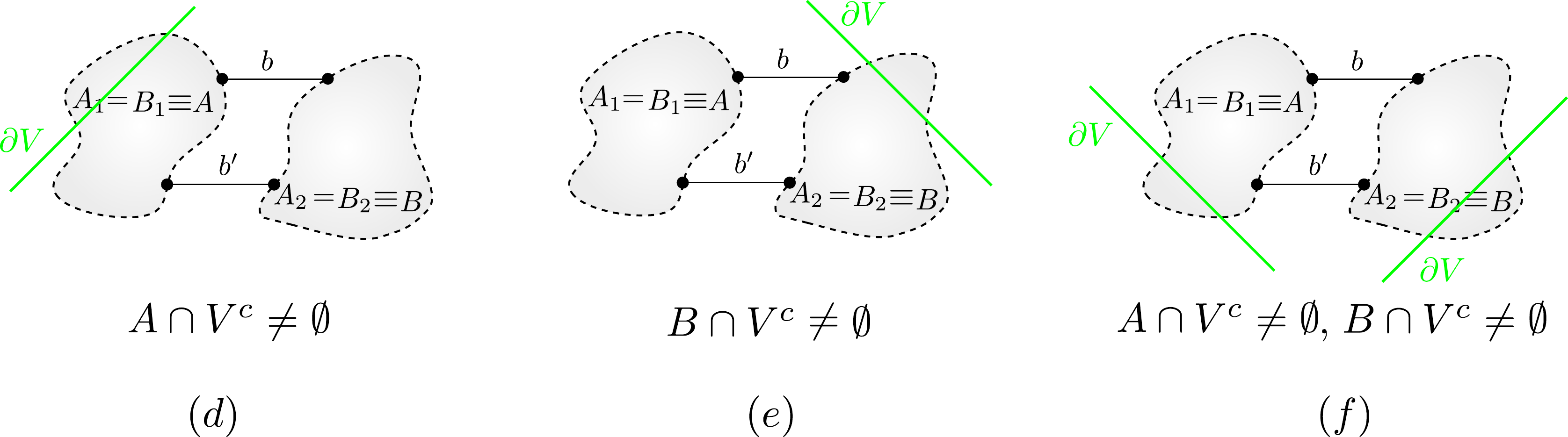}}
			%\captionof{figure}{}
		\end{minipage}
	\end{center}
	For $(d)$ the inequality \eqref{Desigualdade mt importante},
	is valid as long as
	\[
	{\textstyle
	a_\mathrm{\widetilde{m}}
	\left(
	\sum_{m'=1}^q\ q_{m'}b_{m'}
	\right)
	\geqslant
	a_\mathrm{\widetilde{m}}b_\mathrm{\widetilde{m}}.
	}
	\]
	This is in fact true because 
	\[
	{\textstyle
	\sum_{m'=1}^q\ q_{m'}b_{m'}
	\geqslant 
	\underset{   
	m'
	\in
	\cap_{i\in \mathbb{V}}
	\mathcal{Q}_{i,\mathrm{ max}}(\pmb{\hat{h}})   
	}{\sum} 
	q_{m'}b_{\mathrm{ max}}
	\geqslant 
	b_{\mathrm{ max}}=b_\mathrm{\widetilde{m}}.
	}
	\]
	For $(e)$, the desired inequality follows from
	\[
	{\textstyle
	\left(\sum_{m=1}^q q_m a_m\right)b_\mathrm{\widetilde{m}}
	\geqslant 
	a_\mathrm{\widetilde{m}} b_\mathrm{\widetilde{m}},
	}
	\]
	but this inequality holds because 
	\[
	\sum_{m=1}^q\ q_{m}a_{m}
	\geqslant 
	\underset{   
	m\in\cap_{i\in \mathbb{V}}\mathcal{Q}_{i,\mathrm{ max}}
	(\pmb{\hat{h}})   
	}{\sum}
	\ q_{m}a_{\mathrm{ max}}
	\geqslant 
	a_{\mathrm{ max}}=a_\mathrm{\widetilde{m}}.
	\]
	For the last sub-case $(f)$, we have to prove that
	$
	a_\mathrm{\widetilde{m}} b_\mathrm{\widetilde{m}}
	=
	a_\mathrm{\widetilde{m}} b_\mathrm{\widetilde{m}},
	$
	which is obviously true.
	
	In the max wired boundary conditions, if 
	the end vertices of $b$ or $b'$ belong to
	the same component, the result follows. 
\\\\
\noindent
To finish the proof we need to address the case 
where $\xi$ has at most one zero and 
$\omega^{(1)}\equiv \xi\vee\omega^{(b)}$
and $\omega^{(2)}\equiv \xi\vee\omega^{(b')}$.
Suppose that
	\[
	\xi\equiv (1,\ldots,1,\aunderbrace[D]{0}_{b\mathrm{-th}},
	1,\ldots,1,\aunderbrace[D]{1}_{b'\mathrm{-th}},1,\ldots,1)
	\]
	where $b,b'\in E\cup \partial E$ with $b\not=b'$.
	By defining  
	\[
	\xi^b\equiv (1,\ldots,1,\aunderbrace[D]{1}_{b\mathrm{-th}},
	1,\ldots,1,\aunderbrace[D]{1}_{b'\mathrm{-th}},1,\ldots,1)
	\]
	and
	\[
	\xi^{b'}\equiv (1,\ldots,1,\aunderbrace[D]{0}_{b\mathrm{-th}},
	1,\ldots,1,\aunderbrace[D]{1}_{b'\mathrm{-th}},1,\ldots,1),
	\]
	we can see that 
	$
	\omega^{(1)}=\xi\vee\omega^{(b)}=\xi^b \text{,} \quad
	\omega^{(2)}=\xi\vee\omega^{(b')}=\xi^{b'}
	\quad \text{and} \quad
	\omega^{(1)}\wedge\omega^{(2)}=\xi.
	$
	In this case, to prove \eqref{min} for both 
	``free'' and ``max'' wired boundary conditions, it is enough 
	to prove that
	$
	\mathcal{R}(\xi^b,\xi^{b'})=\mathcal{R}(\xi,\xi^{b'}) \ 
	\mathrm{with} \ b\not = b',
	$
	but this is trivial since 
	\begin{multline*}
	\mathcal{R}(\xi^b,\xi^{b'})=
	{
		\mathcal{W}^{\mathrm{ GRC}}_{V,\mathrm{ \#}}
		(\xi^b\vee\xi^{b'})
		\over 
		\mathcal{W}^{\mathrm{ GRC}}_{V,\mathrm{ \#}}(\xi^b)
	}
	=
	{
		\mathcal{W}^{\mathrm{ GRC}}_{V,\mathrm{ \#}}(\xi^b)
		\over 
		\mathcal{W}^{\mathrm{ GRC}}_{V,\mathrm{ \#}}(\xi^b)
	}
	=
	{
		\mathcal{W}^{\mathrm{ GRC}}_{V,\mathrm{ \#}}(\xi)
		\over 
		\mathcal{W}^{\mathrm{ GRC}}_{V,\mathrm{ \#}}(\xi)
	}
	\\
	=
	{
		\mathcal{W}^{\mathrm{ GRC}}_{V,\mathrm{ \#}}
		(\xi\vee\xi^{b'})
		\over 
		\mathcal{W}^{\mathrm{ GRC}}_{V,\mathrm{ \#}}(\xi)
	}
	=
	\mathcal{R}(\xi,\xi^{b'}).
	\end{multline*}
\end{proof}
\section{Edwards-Sokal model}
\paragraph{Edwards-Sokal model with general boundary condition.}
Fix $q\in\mathbb{Z}^+$, for any finite set
$V\subset\mathbb{V}$ and any fixed configurations 
$\sigma_{V^c},\omega_{\mathbb{B}(V)^c}$ prescribed outside of 
$V$, we define de Edwards-Sokal measure 
$\phi^{\mathrm{ ES}}_{V,\mathbb{B}(V)}$ as the normalization of 
the following weights
\begin{eqnarray*}
\mathcal{W}(\sigma_V,\omega_{\mathbb{B}(V)}|
\sigma_{V^c},\omega_{\mathbb{B}(V)^c})
=
\prod_{\substack{\{i,j\}\in\mathbb{B}(V)\\ \omega_{ij}=1}}r_{ij}
\delta_{\sigma_i,\sigma_j} 
\times 
\exp
\big(
\beta \underset{i\in V}{\sum} 
\sum_{p=1}^q h_{i,p}\delta_{\sigma_i,p}
\big),
\end{eqnarray*}
where $r_{ij}$ has been defined  in \eqref{nova def. Bernoulli}.
\paragraph{Edwards-Sokal model with wired and free boundary 
conditions.}
From the previous definition we can observe that, 
for any finite volume $V\subset\mathbb{V}$, the state 
$
\phi^{\mathrm{ ES}}_{V,\mathbb{B}(V)}
(
\cdot|
\sigma_{V^c},\omega_{\mathbb{B}(V)^c}
)
$
is independent of $\omega_{\mathbb{B}(V)^c}$, and we define 
\[
\phi^{\mathrm{ ES}}_{V,\mathrm{ m}}(\cdot)
\equiv
\phi^{\mathrm{ ES}}_{V,\mathbb{B}(V)}
(
\cdot|
\sigma_{V^c}^\mathrm{m},\omega_{\mathbb{B}(V)^c}
),
\]
where $\sigma^\mathrm{m}$ is the constant configuration, 
$\sigma^\mathrm{m}_i=\mathrm{m}$ for
all $i\in \mathbb{V}$, with $\mathrm{m}\in\{1,\ldots,q\}$ fixed. 
This state is known as the $\mathrm{m}$-wired boundary condition 
state.

By similar reasons we have that 
$
\phi^{\mathrm{ ES}}_{V,\mathbb{B}_0(V)}
(
\cdot
|
\sigma_{V^c},\omega_{\mathbb{B}_0(V)^c}
)
$
does not depend on $\sigma_{V^c}$, provided that the
$\omega$-boundary condition is chosen as 
$\omega_{\mathbb{B}_0(V)^c}=\omega^0_{\mathbb{B}_0(V)^c}$,
where $\omega^0$ denotes the configuration with $\omega^0_{ij}=0$ 
for all $\{i,j\}\in\mathbb{B}(\mathbb{V})$. In this case we 
introduce the notation
\[
\phi^{\mathrm{ ES}}_{V,\mathrm{ free}}(\cdot)
\equiv
\phi^{\mathrm{ ES}}_{V,\mathbb{B}_0(V)}
(
\cdot|
\sigma_{V^c},\omega^0_{\mathbb{B}_0(V)^c}
).
\] 
\section{Gibbs states and limit states}
\paragraph{Gibbs states.}
Let $\mathscr{P}(\Omega)$ denote the set of probability measures 
defined on some probability space $\Omega$.
Since the families $\{\phi_\mathbb{B}^{\mathrm{ GRC}}\}$ and 
$\{\phi_{V,\mathbb{B}(V)}^{\mathrm{ ES}}\}$ are specifications 
(see \cite{Georgii88}), we can define as usual the set of the 
Gibbs measures compatible with these specifications as follows 
\begin{eqnarray}\label{Gibbs-rcm}
\mathscr{G}^{\mathrm{ GRC}}
\equiv
\left\{
\phi\in\mathscr{P}(\Omega): 
\begin{array}{l}
\phi(f)
\stackrel{\mathrm{ DLR}}{=}
\displaystyle
\int \phi_\mathbb{B}^{\mathrm{ GRC}}(f|\omega_{\mathbb{B}^c})
\phi(\textrm{d}\omega), 
\\[0.3cm]
\qquad\qquad \mathrm{supp}(f)\subset \mathbb{B}
\end{array}
\right\}
\end{eqnarray}
and
\[ 
\mathscr{G}^{\mathrm{ ES}}
\equiv
\left\{
\nu\in\mathscr{P}(\Omega\times\Sigma): 
\begin{array}{l}
\nu(f)
\stackrel{\mathrm{ DLR}}{=}
\displaystyle\int \phi_{V,\mathbb{B}(V)}^{
\mathrm{ ES}}(f|\sigma_{V^c},\omega_{\mathbb{B}(V)^c})\,
d\nu(\sigma,\omega), 
\\[0.3cm] 
\ \qquad\qquad
\mathrm{supp}(f)\subset V\times\mathbb{B}(V)
\end{array}
\right\}.
\]
That is, 
$\mathscr{G}^{\mathrm{ GRC}}$ and $\mathscr{G}^{\mathrm{ ES}}$
are the class of probability measures (Gibbs measures) 
that are preserved for their respective probability kernels.
\paragraph{Limit states.}
On the other hand, we define the set of the thermodynamic limits 
of the specification 
$\{\phi_{\mathbb{B}_n}^{\mathrm{ GRC}}\}$, where
$\{\mathbb{B}_n\}$ is a cofinal collection in $\mathbb{E}$ :
\begin{eqnarray}\label{Gibbs-limit}
\mathscr{G}^{\mathrm{ GRC}}_{\mathrm{ lim}}
\equiv
\big\{
\phi\in\mathscr{P}(\{0,1\}^\mathbb{\mathbb{E}}): 
\phi\stackrel{\mathrm{ weak}}{=}\lim_{n\to\infty}
\phi_{\mathbb{B}_n}^{\mathrm{ GRC}}(\cdot|\omega_n)
\big\}.
\end{eqnarray}
In general, it is not easy to relate the sets
$\mathscr{G}^{\mathrm{ GRC}}$
and 
$\mathscr{G}^{\mathrm{ GRC}}_{\mathrm{ lim}}$
due to the lack of quasilocality of the
specifications $\phi_\mathbb{B}^{\mathrm{ GRC}}$.
One case where these sets can be related  
is the case in which we assume the existence of at most one 
connected component with probability one. As a consequence of
Lemma \ref{quasilocality}, one can prove the following relation:
$
\mathscr{G}^{\mathrm{ GRC}}_{\mathrm{ lim}}
\subset 
\mathscr{G}^{\mathrm{ GRC}}
$,
see Lemma \ref{continencia-Gibbs} below. For more
details see \cite{BBCK00}.

By using the FKG property for the GRC model 
and the previous definitions, one can prove the following
theorem which ensures the existence of thermodynamic limit.
\begin{theorem}[Monotonicity and existence of limit states]
	\label{limite termodinamico}
	Let $\beta\geqslant 0$,  
	$\pmb{J}=(J_{ij}:\{i,j\}\in \mathbb{E})$ 
	$\in [0,\infty)^\mathbb{E}$ and 
	$
	\pmb{\hat{h}}
	\equiv 
	(h_{i,p}\in\mathbb{R}:i\in \mathbb{V};\ p=1,\ldots,q)
	$. 
	For each increasing quasilocal 
	function $f$ (see \cite{Georgii88}), 
	\begin{itemize}
	\item[(i)] The following limits exist
	\[\phi_{\mathrm{{max}}}^{\mathrm{{GRC}}}(f)\equiv 
	\lim_{V\uparrow \mathbb{V}}
	\phi_{V, \mathrm{{max}}}^{\mathrm{{GRC}}}(f)
	\quad
	\ \mathrm{and} \
	\quad
	\phi_{\mathrm{{free}}}^{\mathrm{{GRC}}}(f)\equiv 
	\lim_{V\uparrow \mathbb{V}}
	\phi_{V, \mathrm{{free}}}^{\mathrm{{GRC}}}(f).
	\]

	\item[(ii)] 
	If in addition, 
	$
	\mathrm{m}
	\in 
	{\cap}_{i\in\mathbb{V}}\mathcal{Q}_{i,\mathrm{ max}}
	(\pmb{\hat{h}})
	$, 
	then the following limits exist
\[ 
	\phi_{\mathrm{{max}}}^{\mathrm{{ES}}}(f)
	\equiv 
	\lim_{V\uparrow \mathbb{V}}
	\phi_{V, \mathrm{{max}}}^{\mathrm{{ES}}}(f)
	\quad
	\ \mathrm{and} \
	\quad
	\phi_{\mathrm{{free}}}^{\mathrm{{ES}}}(f)
	\equiv 
	\lim_{V\uparrow \mathbb{V}}
	\phi_{V, \mathrm{{free}}}^{\mathrm{{ES}}}(f).
\]

	\item[(iii)] 
	If $\phi\in\mathscr{G}_{\mathrm{ lim}}^{\mathrm{ GRC}}$ or 
	$\phi\in\mathscr{G}^{\mathrm{ GRC}}$, then
	for each increasing quasilocal function $f$ 
	we have
\[
	\phi_{\mathrm{ free}}^{\mathrm{ GRC}}(f)
	\leqslant 
	\phi(f)
	\leqslant \phi_{\mathrm{ max}}^{\mathrm{ GRC}}(f).
\]
	\end{itemize}
\end{theorem}
\begin{proof}
The proof is similar to the proof of Theorem III.1 in 
\cite{BBCK00}.
\end{proof}
Since we are also interested in monotonicity properties with
respect to the magnetic field, 
it is needed to introduce a partial order between two 
fields \cite{BBCK00}. 
Given two arbitrary magnetic fields
$\pmb{\hat{h}}$ and $\pmb{\hat{h}'}$, we say that
\begin{eqnarray}\label{relação de ordem nos campos}
\pmb{\hat{h}}\prec \pmb{\hat{h}'} 
\ \ 
\Longleftrightarrow 
\ \
\forall \ i\in\mathbb{V}: \ 
h_{i,k}-h_{i,l}\leqslant h'_{i,k}-h'_{i,l},
\ \ \
k,l=1,\ldots,q
\end{eqnarray}
whenever $h_{i,k}-h_{i,l}>0.$
\begin{theorem}[Monotonicity with respect to the magnetic field]
	\label{teo-monotonicidade-campo-externo-grc}
	Let $\pmb{\hat{h}}$ and
	$\pmb{\hat{h}'}$ be two arbitrary magnetic fields
	such that $\pmb{\hat{h}}\prec\pmb{\hat{h}'}$. 
	Denote by  
	$\phi^{\mathrm{ GRC},\pmb{\hat{h}}}_{\mathrm{ \#}}$
	and $\phi^{\mathrm{ GRC},\pmb{\hat{h}'}}_{\mathrm{ \#}}$ 
	their respective measures
	defined in Theorem \ref{limite termodinamico},
	where  $\mathrm{\#}$ \ stands for ``$\mathrm{free}$'' 
	or ``$\mathrm{max}$''.	
	Then, for any quasilocal increasing function $f$ we have 
	\begin{eqnarray*}
	\phi^{\mathrm{ GRC},\pmb{\hat{h}}}_{\mathrm{ free}}(f)
	\leqslant 
	\phi^{\mathrm{ GRC},\pmb{\hat{h}'}}_{\mathrm{ free}}(f)
	\qquad\text{and}\qquad
	\phi^{\mathrm{ GRC},\pmb{\hat{h}}}_{\mathrm{ max}}(f)
	\leqslant 
	\phi^{\mathrm{ GRC},\pmb{\hat{h}'}}_{\mathrm{ max}}(f).
	\end{eqnarray*}
\end{theorem}
\begin{proof}
	By the  Holley Theorem, 
	the stochastic domination claimed 
	in the statement of the theorem is proved as long as
	the following lattice condition is satisfied
	\begin{eqnarray}\label{cond. latice}
	\phi^{\mathrm{ GRC},\pmb{\hat{h}}}_{V,\mathrm{ \#}}
	(\omega^{(1)}\vee \omega^{(2)})
	\phi^{\mathrm{ GRC},\pmb{\hat{h}'}}_{V,\mathrm{ \#}}
	(\omega^{(1)}\wedge \omega^{(2)})
	\geqslant 
	\phi^{\mathrm{ GRC},\pmb{\hat{h}'}}_{V,\mathrm{ \#}}
	(\omega^{(1)})
	\phi^{\mathrm{ GRC},\pmb{\hat{h}}}_{V,\mathrm{ \#}}
	(\omega^{(2)})
	\end{eqnarray}
	for all $\omega^{(1)},\omega^{(2)}\in\{0,1\}^E$, where
	$\mathrm{\#}$ denotes the ``free'' and ``max'' wired 
	boundary conditions.
	For details, see Theorem 2.3, item (c), p. 20 in 
	\cite{Grimmett2}.
	It is also well known that \eqref{cond. latice} is a 
	consequence of 
	\begin{eqnarray}\label{divisão de medidas}
	{
		\phi^{\mathrm{ GRC},\pmb{\hat{h}'}}_{V,\mathrm{ \#}}
		(\zeta^e)
		\over 
		\phi^{\mathrm{ GRC},\pmb{\hat{h}'}}_{V,\mathrm{ \#}}
		(\zeta_{(e)})
	}
	\geqslant 
	{
		\phi^{\mathrm{ GRC},\pmb{\hat{h}}}_{V,\mathrm{ \#}}
		(\xi^e)
		\over 
		\phi^{\mathrm{ GRC},\pmb{\hat{h}}}_{V,\mathrm{ \#}}
		(\xi_{(e)})
	},
	\end{eqnarray}
	for any $\xi\preceq \zeta$ and $e\in E$, where 
	$\xi_{(e)}$ ($\xi^e$) is the configuration that agrees with 
	$\xi$ 
	in all edges, except in $e$ where its value is zero (one).
	We shall remark that the notations $\xi_e$ and $\xi_{(e)}$
	have different meaning.
	
	Without loss of generality, we can assume that
	$\xi$ and $\zeta$ are of the form 
	\[
	\xi\equiv(*,\ldots,*,\aunderbrace[D]{0}_{\mathrm{ e-th}},*,
	\ldots,*)
	\quad \text{and} \quad  
	\zeta\equiv(*',\ldots,*',\aunderbrace[D]{0}_{\mathrm{ e-th}},
	*',\ldots,*'),
	\]
	with $\xi\preceq \zeta.$
	Let 
	$k'\equiv \underset{\{i,j\}:\zeta_{ij}=1}{\prod}r_{ij}$ and
	$k\equiv \underset{\{i,j\}:\xi_{ij}=1}{\prod}r_{ij}$. 
	From the definitions we get that 
	\[
	\zeta_{(e)}=\zeta, \quad \xi_{(e)}=\xi, \quad
	\prod_{\{i,j\}:\zeta^e_{ij}=1}r_{ij}=r_e k'
	\quad \text{and} \quad
	\prod_{\{i,j\}:\xi^e_{ij}=1}r_{ij}=r_e k.
	\]
	Therefore
	\begin{eqnarray*}
	{
		B_{\pmb{j}}(\zeta^e)
		\over 
		B_{\pmb{j}}(\zeta_{(e)})
	}
	=
	{
		\prod_{\{i,j\}:\zeta^e_{ij}=1}r_{ij}
		\over 
		\prod_{\{i,j\}:\zeta_{ij}=1}r_{ij}
	}
	=
	r_e
	=
	{
		\prod_{\{i,j\}:\xi^e_{ij}=1}r_{ij}
		\over 
		\prod_{\{i,j\}:\xi_{ij}=1}r_{ij}
	}
	=
	{
		B_{\pmb{j}}(\xi^e)
		\over 
		B_{\pmb{j}}(\xi_{(e)})
	}.
	\end{eqnarray*}
	So it follows from the equation above that 
	\eqref{divisão de medidas} is a consequence of
	\begin{eqnarray}\label{novotheta}
	{
		\Theta^{\pmb{\hat{h}'}}_{V,\mathrm{ \#}}(\zeta^e)
		\over 
		\Theta^{\pmb{\hat{h}'}}_{V,\mathrm{ \#}}(\zeta_{(e)})
	}
	\geqslant 
	{
		\Theta^{\pmb{\hat{h}}}_{V,\mathrm{ \#}}(\xi^e)
		\over 
		\Theta^{\pmb{\hat{h}}}_{V,\mathrm{ \#}}(\xi_{(e)})
	},
	\end{eqnarray}
	for both ``free'' and ``max'' wired boundary conditions.
	
	If  $e=\{x,y\}$ and  $x\leftrightarrow y$ in $\xi$,
	then \eqref{novotheta} is an equality.
	On the other hand, if 
	$x\not\leftrightarrow y$ in $\xi$, then there are two 
	connected components
	$A\equiv C(x,\xi)$ and $B\equiv C(y,\xi)$ containing the 
	vertices $x$ and $y$, respectively.
	If $e$ is an open edge in $\xi$, then the components $A$ and 
	$B$ are connected and will be denoted by $C\equiv A\cup B$. 
	So $\vert C\vert=\vert A\vert +\vert B\vert$, from where we 
	deduce that 
	\begin{eqnarray*}
		{
			\Theta^{\pmb{\hat{h}'}}_{V,\mathrm{ \#}}(\zeta^e)
			\over 
			\Theta^{\pmb{\hat{h}'}}_{V,\mathrm{ \#}}(\zeta_{(e)})
		}
		\cdot
		{
			\Theta^{\pmb{\hat{h}}}_{V,\mathrm{ \#}}(\xi_{(e)})
			\over
			\Theta^{\pmb{\hat{h}}}_{V,\mathrm{ \#}}(\xi^e)
		}
		=
		{
			\Theta^{\pmb{\hat{h}'}}_{V,\mathrm{ \#}}(C)
			\over 
			\Theta^{\pmb{\hat{h}}}_{V,\mathrm{ \#}}(C)
		}
		\cdot
		{
			\Theta^{\pmb{\hat{h}}}_{V,\mathrm{ \#}}(A)
			\Theta^{\pmb{\hat{h}}}_{V,\mathrm{ \#}}(B)
			\over 
			\Theta^{\pmb{\hat{h}'}}_{V,\mathrm{ \#}}(A)
			\Theta^{\pmb{\hat{h}'}}_{V,\mathrm{ \#}}(B)
		},
	\end{eqnarray*}
	for either free or max wired boundary conditions.
	In order to prove \eqref{novotheta}, it is enough to prove 
	that  
	\begin{eqnarray}\label{maior=1}
	{
		\Theta^{\pmb{\hat{h}'}}_{V,\mathrm{ \#}}(C)
		\over 
		\Theta^{\pmb{\hat{h}}}_{V,\mathrm{ \#}}(C)
	}
	.
	{
		\Theta^{\pmb{\hat{h}}}_{V,\mathrm{ \#}}(A)
		\Theta^{\pmb{\hat{h}}}_{V,\mathrm{ \#}}(B)
		\over 
		\Theta^{\pmb{\hat{h}'}}_{V,\mathrm{ \#}}(A)
		\Theta^{\pmb{\hat{h}'}}_{V,\mathrm{ \#}}(B)
	}
	\geqslant 1.
	\end{eqnarray}
	To establish the above inequality, we analyze separately 
	the ``free'' and ``max'' wired boundary condition cases. 
\noindent
\\\\
{\sc Free boundary condition case.}
\\
Keeping the notation used in the proof of the FKG inequality, 
for each $m\in\{1,\ldots,q\}$
we define
	\[
	{
	a_m
	\equiv 
	\exp 
	\big( \beta\underset{i\in A}{\sum}h_{i,m}
	\big), 
	\ 
	b_m
	\equiv 
	\exp 
	\big(\beta\underset{i\in B}{\sum}h_{i,m}
	\big)
	\ \ \text{and} \ \
	c_m
	\equiv 
	\exp
	\big(
	\beta\underset{i\in C}{\sum}h_{i,m}
	\big).
	}
	\]
	Similarly we define $a'_{m}, b'_{m}$ and $c'_{m}$ by 
	replacing $(h_{i,m})$ for $(h'_{i,m})$. 
	With this notation, \eqref{maior=1} reads
	\begin{multline}\label{des1}
	{ \textstyle
	(\sum_{j=1}^q q_j a'_j)
	(\sum_{k=1}^q q_k b'_k)
	(\sum_{l=1}^q q_l a_l b_l)
	}
	\\[0,2cm]
	\leqslant
	{\textstyle 
	(\sum_{j=1}^q q_j a_j)
	(\sum_{k=1}^q q_k b_k)
	(\sum_{l=1}^q q_l a'_l b'_l).
	}
	\end{multline}
	The proof of \eqref{des1} is divided in two steps. 
\noindent
\\\\
\textbf{Step 1}:(move the primes from $a_j$'s) we claim that 
	\begin{multline}\label{des2}
	{\textstyle
	(\sum_{j=1}^q q_j a'_j)
	(\sum_{k=1}^q q_k b'_k)
	(\sum_{l=1}^q q_l a_l b_l)
	}
	\\[0,2cm]
	\leqslant 
	{\textstyle
	(\sum_{j=1}^q q_j a_j)
	(\sum_{k=1}^q q_k b'_k)
	(\sum_{l=1}^q q_l a'_l b_l)
	}.
	\end{multline}
	In fact, we first remark that without loss of generality
	we can assume that $h_{i,l}-h_{i,j}>0$, 
	$\forall i\in \mathbb{V}$. 
	From the hypothesis we have  
	$\pmb{\hat{h}}\prec\pmb{\hat{h}'}$, so
	we get 
	$\forall \ l,j=1,\ldots,q$ and $\forall \ i\in \mathbb{V}$
	that  $h_{i,l}-h_{i,j}\leqslant h'_{i,l}-h'_{i,j}$.
	From the last inequality, it follows that
	\begin{eqnarray}\label{implicacion}
	{a_l\over a_j}\leqslant {a'_l\over a'_j},
	\quad \text{which implies} \quad 
	a'_ja_l-a_ja'_l \leqslant 0.
	\end{eqnarray}
	On the other hand, since $h_{i,l}-h_{i,j}>0$,   
	we have $b_l-b_j>0$. Putting together the last two 
	inequalities yields
	\begin{eqnarray*}
	(a'_ja_l-a_ja'_l)(b_l-b_j)\leqslant 0,
	\end{eqnarray*}
	and we conclude that 
	$
	a'_j a_l b_l \leqslant \left[a'_ja_l-a_ja'_l\right]b_j 
	+ 
	a'_l a_j b_l
	\leqslant a_j a'_l  b_l,
	$
	where in the last inequality we have used \eqref{implicacion}.
	By multiplying the above inequality for $q_j q_k q_l b'_k$ 
	and then summing over $j,k,l=1,\ldots, q$, we prove the claim.
\noindent
\\\\
\textbf{Step 2}:(move the primes from $b_k$'s) we claim that
	\begin{multline}\label{des3}
	{\textstyle
	(\sum_{j=1}^q q_j a_j)
	\left(\sum_{k=1}^q q_k b'_k\right)
	\left(\sum_{l=1}^q q_l a'_l b_l\right)
	}
	\\[0,2cm]
	\leqslant 
	{\textstyle
	(\sum_{j=1}^q q_j a_j)
	(\sum_{k=1}^q q_k b_k)
	(\sum_{l=1}^q q_l a'_l b'_l).
	}
	\end{multline}
	The proof is similar to the one given for the Step 1.
	We assume that 
	$h_{i,l}-h_{i,k}>0$, $\forall \ i\in \mathbb{V}$
	and prove in place of \eqref{implicacion} that 
	$b'_k b_l - b_k b'_l \leqslant 0$, proceeding similarly to 
	reach the conclusion.	
	
	Finally, by piecing together the inequalities \eqref{des2} 
	and \eqref{des3}, 
we obtain \eqref{des1}.
\noindent
\\\\
{\sc Max wired boundary condition case.}
\\
We first observe that if 
$
\mathrm{m}
\in
{\cap}_{i\in\mathbb{V}}\mathcal{Q}_{i,m}(\pmb{\hat{h}})
$
and 
$
\mathrm{\widetilde{m}}
\in 
{\cap}_{i\in\mathbb{V}}
\mathcal{Q}_{i,\widetilde{m}}(\pmb{\hat{h}'})
$,
then $\mathrm{m}=\mathrm{\widetilde{m}}$.
Given two connected components $A$ and $B$, if 
$A\cap V^c=\emptyset$
and $B\cap V^c=\emptyset$, then the inequality follows from 
the free boundary condition case. The remaining cases 
will be analyzed by considering the following cases:
\begin{center} 
\begin{minipage}{\linewidth}
\makebox[\linewidth]{%
\includegraphics[scale=0.79,keepaspectratio=true]
{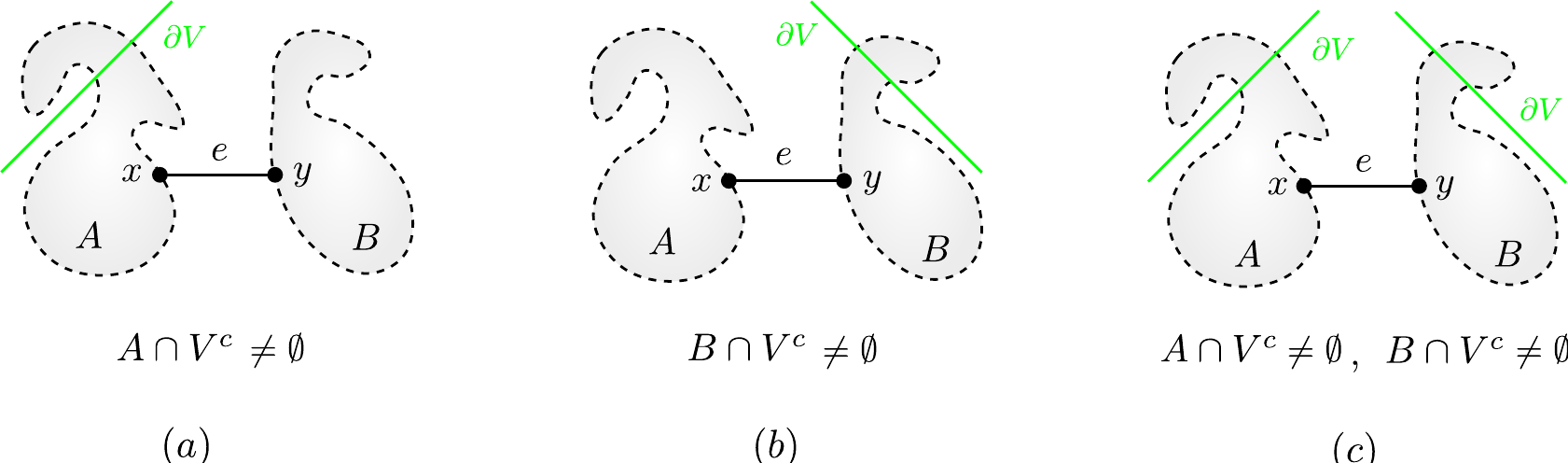}}
%\captionof{figure}{}
\end{minipage}
\end{center}
In the case $(a)$, \eqref{maior=1} is equivalent to 
the inequality
\begin{equation}\label{fin1}
\begin{array}{c}
a'_\mathrm{m}
\left(\sum_{k=1}^q q_k b'_k\right)a_\mathrm{m} b_\mathrm{m}
\leqslant 
a_\mathrm{m} \left(\sum_{k=1}^q q_k b_k\right)a'_\mathrm{m} 
b'_\mathrm{m}.
\end{array}
\end{equation}
To see that this inequality holds, it is sufficient 
to observe that the ordering 
between the magnetic fields implies 
$
b_\mathrm{m} b'_k\leqslant b'_\mathrm{m} b_k.
$
Multiplying this last inequality by $a'_m a_m$ 
and summing over $k=1,\ldots,q$, we obtain \eqref{fin1}.

In the case $(b)$ the inequality \eqref{maior=1} reduces to 
\begin{equation*}
\label{fin2}
\begin{array}{c}
(\sum_{j=1}^q q_j a'_j)b'_\mathrm{m} a_\mathrm{m} b_\mathrm{m} 
\leqslant 
(\sum_{j=1}^q q_j a_j)b_\mathrm{m} a'_\mathrm{m} b'_\mathrm{m}.
\end{array}
\end{equation*}
Now we use that the magnetic field ordering
implies that
$
a_\mathrm{m} a'_j\leqslant a'_\mathrm{m} a_j
$
and then proceed similarly to the previous case.
Finally, in the case $(c)$ the inequality \eqref{maior=1}
is equivalent to 
$
a'_\mathrm{m} b'_\mathrm{m} a_\mathrm{m} b_\mathrm{m}
=  
a'_\mathrm{m} b'_\mathrm{m} a'_\mathrm{m} b'_\mathrm{m},
$
which is trivial.			
\end{proof}
Our next result is the monotonicity, in the FKG 
sense, with respect to the  coupling constants $\pmb{J}$
in the special case where $J_{ij}\equiv J,$
$\forall$ $i,j\in\mathbb{V}$.
\begin{theorem}\label{teorema ge}
	Suppose that $0\leqslant J_1<J_2$ are two coupling constants. 
	For each finite 
	$V\subset\mathbb{V}$ 
	denote by  $\phi^{\mathrm{ GRC},J_k}_{V,\mathrm{ \#}}$ , k=1,2;
	the measure defined by the  weights \eqref{free peso} or by 
	the weights \eqref{def.rcm-wired condition} with 
	$
	\mathrm{m}
	\in
	{\cap}_{i\in\mathbb{V}}
	\in
	\mathcal{Q}_{i,\mathrm{ max}}$$(\pmb{\hat{h}}).
	$
	Then
\[
	\phi^{\mathrm{ GRC},J_1}_{V,\mathrm{ max}}(f)
	\leqslant 
	\phi^{\mathrm{ GRC},J_2}_{V, \mathrm{ free}}(f),
\]
where $f$ is a cylindrical increasing function, and
$\mathrm{\#}$ stands for ``$\mathrm{free}$'' or ``$\mathrm{max}$''.
\end{theorem}%%%%%%%%%%%%%%%%%%%%%%%%%%%%%%%%%%%%
\begin{proof}
By assuming $J_1<J_2$, we get that
$e^{q\beta J_1}-1<e^{q\beta J_2}-1$.
For any configuration $\omega\in\{0,1\}^{\mathbb{B}(V)}$,
we define the function $g:\{0,1\}^{\mathbb{B}(V)}\to\mathbb{R}$ 
by
\[
g(\omega)
\equiv
\left[
	\dfrac{e^{q\beta J_1}-1}{e^{q\beta J_2}-1}
\right]^{o(\omega)}
\times 
\displaystyle
{
\prod_{\substack{C(\omega): \\ \mathbb{V}(C(\omega))\cap 
\partial V\not= \emptyset}}
}
\exp
\big(
\beta\sum_{i\in C(\omega)}h_{i,\mathrm{ max}} 
\big),
\]
where $o(\omega)$ denotes the numbers of open edges in $\omega$. 
One can easily see that the function $g$ is decreasing
since 
$g$ is composed by the product of  non-negative decreasing 
functions.

Let $f:\{0,1\}^{\mathbb{B}(V)}\to\mathbb{R}$ be a cylindrical 
increasing arbitrary function.
Since  
$J_{ij}\equiv J,$ we have the following expression for
Bernoulli factor:
$B_{J_k}(\omega)=\left(e^{q\beta J_k}-1\right)^{o(\omega)}$, 
$k=1,2$.
From the definition of the expected values we obtain
\begin{multline}\label{r_1}
\phi^{\mathrm{ GRC},J_1}_{V,\mathrm{ max}}(f)
=
{
1
\over 
Z^{\mathrm{ GRC},J_1}_{V,\mathrm{ max}}	
}
\sum_{\omega\in\{0,1\}^{\mathbb{B}(V)}}
f(\omega)
{
\left(e^{q\beta J_1}-1\right)^{o(\omega)}
\displaystyle
\prod_{C(\omega)} \Theta_{V,\mathrm{ max}}(C(\omega))
}
\\[0,3cm]
=
{
1
\over 
Z^{\mathrm{ GRC},J_2}_{V,\mathrm{ max}}	
}
\sum_{\omega\in\{0,1\}^{\mathbb{B}(V)}}
f(\omega) g(\omega)
%{
\left(e^{q\beta J_2}-1\right)^{o(\omega)}
\\[0,2cm]
\hspace*{4,2cm}
\times
\prod_{{\substack{C(\omega): \\ \mathbb{V}(C(\omega))\cap 
\partial V= \emptyset}}}
\sum_{p=1}^q
e^{\beta\sum_{i\in C(\omega)}h_{i,p}}
\times
{ 
Z^{\mathrm{ GRC},J_2}_{V,\mathrm{ free}}
\over 
Z^{\mathrm{ GRC},J_1}_{V,\mathrm{ max}}
}
\\[0,3cm]
=
\phi^{\mathrm{ GRC},J_2}_{V,\mathrm{ free}}(f\cdot g)
\times
{
Z^{\mathrm{ GRC},J_2}_{V,\mathrm{ free}}
\over 
Z^{\mathrm{ GRC},J_1}_{V,\mathrm{ max}}
},
\end{multline}
\noindent
\\[0,2cm]
where
$Z^{\mathrm{ GRC},J}_{V,\mathrm{ \#}}$ 
denotes the normalization constant of the measure
$\phi^{\mathrm{ GRC},J}_{V,\mathrm{ \#}}$
and
$\mathrm{\#}$ stands for \textrm{``free''} or \textrm{``max''}.
By taking $f\equiv1$ in \eqref{r_1} 
we get the following equality
\[
\phi^{\mathrm{ GRC},J_2}_{V,\mathrm{ free}}(g)
=
{
Z^{\mathrm{ GRC},J_1}_{V,\mathrm{ max}}
\over 
Z^{\mathrm{ GRC},J_2}_{V,\mathrm{ free}}
}.
\]
Using the last equation, \eqref{r_1} 
and the strong FKG property 
(Theorem \ref{Strong FKG Property}) we  finally conclude that
\[
\phi^{\mathrm{ GRC},J_1}_{V,\mathrm{ max}}(f)
=
{
\phi^{\mathrm{ GRC},J_2}_{V,\mathrm{ free}}(f\cdot g)
\over 
\phi^{\mathrm{ GRC},J_2}_{V,\mathrm{ free}}(g)
}
\stackrel{\mathrm{ FKG}}{\leqslant} 
\phi^{\mathrm{ GRC},J_2}_{V,\mathrm{ free}}(f).
\]
\end{proof}
\begin{remark}
Note that Theorem \ref{teorema ge} can be extended 
using 
Item $(iii)$ of Theorem \ref{limite termodinamico} for 
any pair of $\mathrm{GRC}$ Gibbs measures at $J=J_1,$ resp. 
$J=J_2$. As a particular case, we obtain the following corollary.
\end{remark}
\begin{corollary}[Monotonicity in coupling constant] 
	\label{coro ge0}
	Suppose that $0\leqslant J_1<J_2$ are two coupling constants. 
	For each finite $V\subset\mathbb{V}$ 
	denote by  $\phi^{\mathrm{ GRC},J_k}_{V,\mathrm{ \#}}$ , 
	k=1,2;
	the measure defined by the  weights \eqref{free peso} or by 
	the weights \eqref{def.rcm-wired condition} with 
	$
	\mathrm{m}\in {\cap}_{i\in\mathbb{V}}
	\in\mathcal{Q}_{i,\mathrm{ max}}
	$
	$(\pmb{\hat{h}}).$
	Then
\[
	\phi^{\mathrm{ GRC},J_1}_{V,\mathrm{ \#}}(f)
	\leqslant 
	\phi^{\mathrm{ GRC},J_2}_{V, \mathrm{ \#}}(f),
\]
where $f$ is a cylindrical increasing function and
$\mathrm{\#}$ stands for ``$\mathrm{free}$'' or ``$\mathrm{max}$''.
\end{corollary}
\section{GRC model and quasilocality}
In what follows we study the quasilocality of the 
random-cluster model in non-homogeneous magnetic field. 
The next lemma tells us that the specifications
$\{\phi_{\mathbb{B}}^{\mathrm{{GRC}}}\}$ are 
almost surely quasilocal
(see \cite{Georgii88,Pfister}). 
To give a precise statement of this lemma,
we need to introduce some notation:
\[
\mathscr{M}(\Delta,\Lambda)
\equiv
\Big\{ 
\!
\omega\in \{0,1\}^{\mathbb{E}} 
: 
\forall x,y\in\Lambda, \ x\leftrightarrow\Delta^c 
\ \text{and} \
y\leftrightarrow\Delta^c 
\Rightarrow 
x \xleftrightarrow[\vspace*{-1cm}\mathbb{B}_0(\Delta)]{}y
\!
\Big\}
\]
where $\Lambda\subset\Delta$ are finite subsets in $\mathbb{V}$.
The following lemma is an adaptation of Lemma VI.2 in 
\cite{BBCK00} for our model.
\begin{lemma}[Quasilocality]\label{quasilocality}
Let $\mathbb{B}\subset\mathbb{B}_0(\mathbb{E})$ be a finite 
set and $f$ a cylindrical function depending only 
on the edges in $\mathbb{B}$. 
Then, for each pair of  finite subsets 
$(\Delta,\Lambda)$ with 
$\mathbb{V}(\mathbb{B})\subset \Lambda\subset \Delta$, 
the function
\[
\omega\mapsto \mathds{1}_{\mathscr{M}(\Delta,\Lambda)}(\omega) 
\phi_{\mathbb{B}}^{\mathrm{ GRC}}(f\vert \omega_{\mathbb{B}^c})
\]
is quasilocal. If in addition 
$\phi\in \mathscr{G}^{\mathrm{{GRC}}}_{\mathrm{{lim}}}$ 
or $\phi\in \mathscr{G}^{\mathrm{{GRC}}}$ 
have at most one infinite connected component and 
$\Lambda\subset \mathbb{V}$, then
\[
\phi(\mathscr{M}(\Delta,\Lambda))\uparrow 1, 
\quad \text{whenever} \ \ \Delta\uparrow\mathbb{V}.
\]
\end{lemma}
\begin{proof}
Recalling the definition of
$
\phi_{\mathbb{B}}^{\mathrm{ GRC}}
(\cdot |\omega_{\mathbb{B}^c})
$,
we note that it is enough to prove that the function  
\begin{eqnarray}\label{quasilocalidade}
\omega\mapsto \mathds{1}_{\mathscr{M}(\Delta,\Lambda)}(\omega)
\phi_{\mathbb{B}}^{\mathrm{ GRC}}(\overline{\omega}_\mathbb{B} |
\omega_{\mathbb{B}^c}),
\quad 
\forall \ \overline{\omega}_\mathbb{B}\in \{0,1\}^{\mathbb{B}}
\end{eqnarray}
is quasilocal.
In the sequel we shall prove the quasilocality of
the mapping defined in \eqref{quasilocalidade}.
Let $\tilde{\Delta}$ be a finite subset of $\mathbb{V}$ 
such that $ \Delta\subset\tilde{\Delta}$. Consider the following 
configurations:
\[
\begin{array}{c}
\omega
\equiv
(*,\ldots,*,\aunderbrace[D]{0}_{b\mathrm{-th}},*,\ldots),
\qquad
\omega^b
\equiv
(*,\ldots,*,\aunderbrace[D]{1}_{b\mathrm{-th}},*,\ldots)
\end{array}
\]
where $*$ is an arbitrary element in $\{0,1\}$ and 
$b\in\mathbb{B}(\tilde{\Delta})^c$.
Suppose that $\omega\in\mathscr{M}(\Delta,\Lambda)$ and that 
there exists a connected component $C^{*}$
connecting $\Lambda$ to $\mathbb{B}(\tilde{\Lambda})^c$ in
$\omega$. By definition, we have that 
$\omega^b\in \mathscr{M}(\Delta,\Lambda)$ and  
the connected component $C^{*}$ is unique.
Let us consider two cases:
\[
	\mathrm{1)} \ 
	\mathbb{V}(C^{*})\cap\mathbb{V}(\{b\})=\emptyset
	\quad 
	\mathrm{and}
	\quad
	\mathrm{2)} \ 
	\mathbb{V}(C^{*})\cap\mathbb{V}(\{b\})\neq \emptyset.
\]
In the first case we trivially have
$
|W_\mathbb{B}^{\mathrm{ GRC}}(\overline{\omega}_\mathbb{B} | 
\omega^b_{\mathbb{B}^c})
-
W_\mathbb{B}^{\mathrm{ GRC}}(\overline{\omega}_\mathbb{B} | 
\omega_{\mathbb{B}^c})|=0.
$
The second case is more elaborate. 
We consider separately two cases.
We first assume that there is some 
$\varepsilon>0$ such that 
$\varepsilon<|h_{i,\mathrm{ max}}-h_{i,m}|$
for all $i\in \mathbb{V}$ and $m\in \{1,\ldots,q\}$.
For this case 
let us denote by $C^*_b$ the connected component 
$\mathbb{V}(C^{*})\cap\mathbb{V}(\{b\})$.
Then
\begin{multline*}
|W_\mathbb{B}^{\mathrm{ GRC}}(\overline{\omega}_\mathbb{B} |
\omega^b_{\mathbb{B}^c})
-
W_\mathbb{B}^{\mathrm{ GRC}}(\overline{\omega}_\mathbb{B} | 
\omega_{\mathbb{B}^c})|
\leqslant 
\overline{B}_J(\omega) k(\omega)\times
\\[0,3cm]
\times
\sum_{m=1}^qq_m
\Big|
\exp
\big(
-\beta\sum_{i\in C^*_b}(h_{i,\mathrm{ max}}-h_{i,m})
\big)
-
\exp
\big(
-\beta\sum_{i\in C^*}(h_{i,\mathrm{ max}}-h_{i,m})
\big)
\Big|,
\end{multline*}
where
\[
k(\omega)\equiv \prod_{\substack{C(\omega):\mathbb{V}(C)
\cap\mathbb{V}
(\mathbb{B})\not=\emptyset \\|C|<\infty}}
\,
\sum_{m=1}^q q_m
\exp 
\big(
-\beta\sum_{i\in C}(h_{i,\mathrm{ max}}-h_{i,m})
\big)
<
\infty
\]
and
\[
\overline{B}_J(\omega)
=
\prod_{\{i,j\}\in\mathbb{B}:\omega_{ij}=1}
r_{ij}.
\]
Suppose that 
$
m
\not\in 
\underset{i\in\mathbb{V}}{\cap}\mathcal{Q}_{i,\mathrm{ max}}
(\pmb{\hat{h}})
$,
then we have 
\[
d(b,\Lambda)\leqslant |\mathbb{V}(C^*_b)|
\leqslant 
\sum_{i\in C^{*}_b}
\frac{1}{\varepsilon}
(h_{i,\mathrm{ max}}-h_{i,m}).
\]
The last inequalities imply that if $d(b,\Lambda)\to \infty$, 
then
\[
\sum_{i\in C^{*}_b}(h_{i,\max}-h_{i,m})\to \infty
\qquad \text{and} \qquad
\sum_{i\in C^{*}}(h_{i,\mathrm{ max}}-h_{i,m})\to \infty,
\]
whenever 
$
m
\not\in 
\underset{i\in\mathbb{V}}{\cap}
\mathcal{Q}_{i,\mathrm{ max}}(\pmb{\hat{h}})
$.
Therefore, whenever $d(b,\Lambda)\to\infty$, we have that 
\begin{eqnarray}\label{condicao fraca de quasilocalidade}
|W_\mathbb{B}^{\mathrm{ GRC}}(\overline{\omega}_\mathbb{B} |
	\omega^b_{\mathbb{B}^c})
-
W_\mathbb{B}^{\mathrm{ GRC}}(\overline{\omega}_\mathbb{B} |
	\omega_{\mathbb{B}^c})|
\to 0.
\end{eqnarray}
In the case 
$\liminf_{i\in \mathbb{V}} |h_{i,\mathrm{ max}}-h_{i,m}|=0$
it is enough to analyze whether 
\[
\sum_{i\in C^*}(h_{i,\mathrm{ max}}-h_{i,m})
\]
is finite or not. If it is infinite, then the result
is trivial. Otherwise we use the continuity of the 
exponential function and a suitable choice of $b$
so that $d(b,\Lambda)\to\infty$.
Now we consider two different configurations:
\[
\hat{\omega}
\equiv 
(*,\ldots,*,\aunderbrace[D]{0}_{b\mathrm{-th}},*,\ldots,*,
\aunderbrace[D]{1}_{ b'\mathrm{-th}} ,*,\ldots)
\]
and
\[
\tilde{\omega}
\equiv 
(*,\ldots,*,\aunderbrace[D]{1}_{b\mathrm{-th}},*,\ldots,*,
\aunderbrace[D]{0}_{ b'\mathrm{-th}} ,*,\ldots),
\]
where $*$ is arbitrary in $\{0,1\}$ and 
$b, b'\in\mathbb{B}(\tilde{\Delta})^c$ .
We also denote
\[
\omega
\equiv 
(*,\ldots,*,\aunderbrace[D]{0}_{b\mathrm{-th}},*,\ldots,*,
\aunderbrace[D]{0}_{ b'\mathrm{-th}} ,*,\ldots,*),
\]
with $*$ arbitrary in $\{0,1\}$.
See that 
$\hat{\omega}=\omega^{b'}$ and $\tilde{\omega}=\omega^b$. 
Then by 
\eqref{condicao fraca de quasilocalidade} and the triangle 
inequality we have
\begin{multline*}
\hspace*{3,7cm}
|W_\mathbb{B}^{\mathrm{ GRC}}(\overline{\omega}_\mathbb{B} |
\hat{\omega}^b_{\mathbb{B}^c})
-
W_\mathbb{B}^{\mathrm{ GRC}}(\overline{\omega}_\mathbb{B} |
\tilde{\omega}_{\mathbb{B}^c})|
\\
\leqslant
\\
|W_\mathbb{B}^{\mathrm{ GRC}}(\overline{\omega}_\mathbb{B} |
\omega^{b'}_{\mathbb{B}^c})
-
W_\mathbb{B}^{\mathrm{ GRC}}(\overline{\omega}_\mathbb{B} |
\omega_{\mathbb{B}^c})|
+
|W_\mathbb{B}^{\mathrm{ GRC}}(\overline{\omega}_\mathbb{B} |
\omega^b_{\mathbb{B}^c})
-
W_\mathbb{B}^{\mathrm{ GRC}}(\overline{\omega}_\mathbb{B} |
\omega_{\mathbb{B}^c})|
\to 0,
\end{multline*}
when $d(b',\Lambda), d(b,\Lambda)\to \infty.$
Following this reasoning,  we prove that for any two distinct 
configurations
$\hat{\omega}, \tilde{\omega}$ em $\mathbb{B}(\tilde{\Delta})^c$ 
we have
\[
\left|W_\mathbb{B}^{\mathrm{ GRC}}(\overline{\omega}_\mathbb{B}|
\hat{\omega}^b_{\mathbb{B}^c})
-
W_\mathbb{B}^{\mathrm{ GRC}}(\overline{\omega}_\mathbb{B} |
\tilde{\omega}_{\mathbb{B}^c})\right|
\to 0,
\quad \text{whenever} \quad \min_{b\in A}d(b,\Lambda)\to\infty
\]
and $A\equiv \{e\in\mathbb{B}(\tilde{\Delta})^c: 
\hat{\omega}_e\not=\tilde{\omega}_e\}$, thus proving the
quasilocality of the application \eqref{quasilocalidade}. 

Finally, in order to prove the second statement, it is enough to 
notice that 
$
\{
\mathscr{M}(\Delta,\Lambda):\Delta\subset \mathbb{V} \
\text{finite}
\}
$ 
is an increasing sequence of events.
\end{proof}
\begin{lemma}[Subsets of Gibbs measures]
\label{continencia-Gibbs}
Let $q\in \mathbb{Z}^+$, $\beta\geqslant 0$, 
$\pmb{J}=(J_{ij}:\{i,j\}\in \mathbb{E})$ 
$\in [0,\infty)^\mathbb{E}$, 
$
\pmb{\hat{h}}
=
(h_{i,p}\in\mathbb{R}: i\in\mathbb{V},\ p=1,\ldots, q)
$ 
and $\{q_p: p=1,\ldots,q\}$ satisfying \eqref{qmmaiorque1}.
If
$\phi\in \mathscr{G}^{\mathrm{ GRC}}_{\mathrm{ lim}}$ 
and has at most one infinite connected component then  
$
\mathscr{G}^{\mathrm{ GRC}}_{\mathrm{ lim}}
\subset 
\mathscr{G}^{\mathrm{ GRC}}.
$
\end{lemma}
\begin{proof}
The proof of this lemma follows from Lemma \ref{quasilocality} 
and the almost sure quasilocality. See \cite{BBCK00}. 
\end{proof}
The next theorem gives sufficient conditions for quasilocality 
of the specifications under a geometric assumption
of almost sure existence of an infinite connected component 
in the graph, thus facilitating many technical calculations.
\begin{theorem}[Conditional expectations for GRC]\label{vec}
Let
$\beta \geqslant 0, J_{ij}\geqslant 0,$ 
$h_{i,m}\in\mathbb{R}$, $\forall i,j\in \mathbb{V}$ and 
$q_m>0, m=1,\ldots,q$ 
satisfying \eqref{qmmaiorque1}. 
If $\phi\in\mathscr{G}^{\mathrm{{GRC}}}$ 
and has at most one infinite connected component almost surely,  
$\mathbb{B}\subset\mathbb{B}_0(\mathbb{V})$, and $f$ is a 
cylindrical function depending 
on the configuration $\omega_{\mathbb{B}}$, then
\[
\phi(f\vert \mathscr{F}_{\mathbb{B}^c})
=
\phi^{\mathrm{ GRC}}_{\mathbb{B}}(f\vert \omega_{\mathbb{B}^c}),
\quad \phi \mathrm{-a.s.}
\]
\end{theorem}
\noindent
\textit{Sketch of the Proof.}
The idea is the same as the one employed in the proof of the 
Theorem III.4 in \cite{BBCK00}.
For the sake of completeness, we sketch a proof.		
Let $\mathbb{B}_1,$ $\mathbb{B}_2$ be finite sets of bonds with 
$\mathbb{B}_1\cap \mathbb{B}_2=\emptyset$ and $f$ and $g$ be 
bounded cylinder functions depending only on the bonds in 
$\mathbb{B}_1$ and $\mathbb{B}_2$, 
respectively. Using the DLR equation \eqref{Gibbs-rcm} and the 
consistence of the specifications 
$\{\phi_\mathbb{B}^{\mathrm{ GRC}}\}$, for 
$\mathbb{B}\supset \mathbb{B}_1\cap\mathbb{B}_2$
we can easily obtain that
\begin{eqnarray}\label{c31}
\phi(gf)=\lim_{\mathbb{B}\uparrow\mathbb{E}}
\int \phi_\mathbb{B}^{\mathrm{ GRC}}
\big(g\phi_{\mathbb{B}_1}^{\mathrm{ GRC}}(f|\cdot) 
\big|\omega_{\mathbb{B}^c} \big)
\phi(\mathrm{d}\omega).
\end{eqnarray}
Let $\Delta\supset \mathbb{V}(\mathbb{B}_1)$, since both $g$ and 
$
\mathds{1}_{\mathscr{M}(\Delta, \mathbb{V}(\mathbb{B}_1))}(\cdot)
\phi_{\mathbb{B}}^{\mathrm{ GRC}}(f\vert \cdot)
$ 
are quasilocal, the function
$
g\cdot \mathds{1}_{\mathscr{M}(\Delta, \mathbb{V}(\mathbb{B}_1))}
\phi_{\mathbb{B}}^{\mathrm{ GRC}}(f\vert \cdot)
$ can be 
approximated  by local functions. Then
by DLR equation \eqref{Gibbs-rcm}, we have 
\begin{multline}\label{c33}
\phi
\big(
g\cdot \mathds{1}_{\mathscr{M}(\Delta, \mathbb{V}(\mathbb{B}_1))}
\phi_{\mathbb{B}_1}^{\mathrm{ GRC}}(f\vert \cdot)
\big)
\\[0,2cm]
=
\lim_{\mathbb{B}\uparrow\mathbb{E}}
\int \phi_\mathbb{B}^{\mathrm{ GRC}}
\big(
g\cdot \mathds{1}_{\mathscr{M}(\Delta, \mathbb{V}(\mathbb{B}_1))}
\phi_{\mathbb{B}_1}^{\mathrm{ GRC}}(f|\cdot) 
\big|\omega_{\mathbb{B}^c} 
\big)
\phi(\mathrm{d}\omega).
\end{multline}
From Lemma \ref{quasilocality} we get  
$
\phi\big(\mathscr{M}(\Delta, \mathbb{V}(\mathbb{B}_1))\big)
\uparrow 
1
$ 
whenever 
$\Delta \uparrow \mathbb{V}$. Since $f$ and $g$ are bounded, 
using the Dominated Convergence Theorem we have
\[
\lim_{\Delta\uparrow\mathbb{V}}
\int \!\!\phi_\mathbb{B}^{\mathrm{ GRC}}
\big(
g\cdot \mathds{1}_{\mathscr{M}(\Delta, \mathbb{V}(\mathbb{B}_1))}
\phi_{\mathbb{B}_1}^{\mathrm{ GRC}}(f|\cdot) 
\big|\omega_{\mathbb{B}^c} 
\big)
\phi(\mathrm{d}\omega)
=\!\!
\int \!\!\phi_\mathbb{B}^{\mathrm{ GRC}}
\big(g
\phi_{\mathbb{B}_1}^{\mathrm{ GRC}}(f|\cdot) 
\big|\omega_{\mathbb{B}^c} \big)
\phi(\mathrm{d}\omega)
\]
and 
\[
\lim_{\Delta\uparrow\mathbb{V}}
\phi
\big(
g\cdot \mathds{1}_{\mathscr{M}(\Delta, \mathbb{V}(\mathbb{B}_1))}
\phi_{\mathbb{B}_1}^{\mathrm{ GRC}}(f\vert \cdot)
\big)
=
\phi
\big(
g
\phi_{\mathbb{B}_1}^{\mathrm{ GRC}}(f\vert \cdot)
\big).
\]
Combining the above limits, together with the items \eqref{c31} 
and \eqref{c33}, we have 
\[
\phi(gf)=\phi
\big(
g\phi_{\mathbb{B}_1}^{\mathrm{ GRC}}
(f\vert \cdot)
\big)
\] 
for all bounded $g$ depending only 
on the configurations $\omega_{\mathbb{B}_1^c}$. 
From the almost sure uniqueness of conditional expectation 
with respect to $\phi$, the proof follows.
\qed

Using the general theory of thermodynamic formalism, one can prove
the following lemma.
\begin{lemma}[Monotonicity in the volume, 
\cite{BBCK00}]\label{monotonicidade no volume}
	
Let $q\in \mathbb{Z}^+$, $\beta\geqslant 0$, 
$\pmb{J}=(J_{ij}:\{i,j\}\in \mathbb{E})$ 
$\in [0,\infty)^\mathbb{E}$, 
	the magnetic field be $\pmb{\hat{h}}$ and the sequence
	$\{q_p: p=1,\ldots,q\}$ 
	satisfy \eqref{qmmaiorque1}.	
	If $\Lambda\subset V$ are finite subsets
	of $\mathbb{V}$, then for any cylindrical increasing function 
	$f$ we have
\[
	\phi_{\Lambda,\mathrm{ free}}^{\mathrm{ GRC}}(f)
	\leqslant 
	\phi_{V,\mathrm{ free}}^{\mathrm{ GRC}}(f)
	\quad 
	\text{and} 
	\quad
	\phi_{\Lambda,\mathrm{ max}}^{\mathrm{ GRC}}(f)
	\geqslant 
	\phi_{V,\mathrm{ max}}^{\mathrm{ GRC}}(f).
\]
\end{lemma}
\begin{remark} When $q_p=1$, for all $p=1,\ldots,q$ in 
\eqref{definição de general boundary condition} ,
then we call the model simply the $\mathrm{RC}$ model.
In this case, 
we define the set of Gibbs measures
$\mathscr{G}^{\mathrm{ RC}}$ and 
$\mathscr{G}^{\mathrm{ RC}}_{\mathrm{ lim}}$
similarly to \eqref{Gibbs-rcm} and \eqref{Gibbs-limit}.
\end{remark}
From now on, 
the study turns to some fundamental properties 
of the $\mathrm{RC}$ model.
The following theorem is valid only for the random-cluster model.
\begin{theorem}\label{principal-marginal}
	Let $q\in \mathbb{Z}^+$, $\beta\geqslant 0$, 
	$\pmb{J}=(J_{ij}:\{i,j\}\in \mathbb{E})$ 
	$\in [0,\infty)^\mathbb{E}$ and $\pmb{\hat{h}}$ be a
	magnetic field as previously defined.
	Given 
	$
	\nu\in\mathscr{G}^{\mathrm{{ES}}}$, let 
	$\phi_\nu$  
	denote its edge-marginal.
	Then for any cylindrical increasing function 
	$f$ we have
$
	\phi_{\nu}(f)
	\leqslant 
	\phi_{\mathrm{ max}}^{\mathrm{ RC}}(f).
$
\end{theorem}
\begin{proof}
For more details see the proof of Theorem III.2 reference 
\cite{BBCK00}.
\end{proof}
\section{Uniqueness of the infinite connected component}
We have so far developed the theory of the random-cluster model 
with non-uniform magnetic field for countably infinite graphs.
We are interested in the situation in which the infinite 
connected component is (almost surely) unique, as is commonly 
the case for an ``amenable graph''. The amenability hypothesis
is important for the uniqueness of the infinite connected 
component in several models, \cite{Burton,G06,Hag06}. When the 
graph is non-amenable, the non-uniqueness of the infinite 
connected component is known for several models including the 
Bernoulli percolation and null magnetic field random-cluster 
model, see \cite{BS11,G06,Hag06} and references therein. 
Therefore, from now on we assume tacitly that the lattice 
$\mathbb{L}$ is amenable, that is, 
$
\inf 
\{|\partial_{\mathbb{E}} V|/|V|\}
=
0,
$
where the infimum ranges over all finite connected subsets $V$ 
of $\mathbb{V}$, and $\partial_{\mathbb{E}}V$ is the set of edges 
with one end-vertex in V and one in $\mathbb{V}\setminus V$.

In what follows we denote by
$N_\infty$ 
the random variable that counts the number of 
infinite connected components in both sample spaces 
$\Omega\equiv \{0,1\}^{\mathbb{E}}$ and $\Sigma_{q}\times\Omega$. 
\begin{theorem}[Uniqueness of the infinite connected component]
\label{teo-unicidade-aglo-infinito}
Let $\beta>0$ be the inverse temperature 
and  $\pmb{\hat{h}}$ a magnetic field.
Then
\[
\phi_{\mathrm{ max}}^{\mathrm{ GRC}, \pmb{\hat{h}} }
(N_{\infty}\leqslant 1)
=
\phi_{\mathrm{ free}}^{\mathrm{ GRC}, \pmb{\hat{h}} }
(N_{\infty}\leqslant 1)
=1.
\]
\end{theorem}
\begin{proof}
We only present the argument for 
$
\phi_{\mathrm{ max}}^{\mathrm{GRC}, \pmb{\hat{h}}}
\in
\mathscr{G}^{\mathrm{ GRC}}
$,
since for the free boundary condition case the proof works 
similarly. 
Let $\Lambda\subset \Delta$ be finite subsets of $\mathbb{V}$ 
and
$D_{\Lambda,\Delta}$ the set of all $\omega\in \Omega$
with the property: there exist two points 
$u,v\in \partial \Lambda$
such that both $u$ and $v$ are joined to $\partial\Delta$ by 
paths using 
$\omega$-open edges of 
$\mathbb{E}_{\Delta}\setminus \mathbb{E}_{\Lambda}$,
but $u$ is not joined to $v$ by a path using $\omega$-open edges 
of $\mathbb{E}_{\Delta}$. For any fixed configuration 
$\eta\in\Omega$ the mapping  
$
\omega
\mapsto
\mathds{1}_{D_{\Lambda,\Delta}}
(\omega_{\mathbb{E}_{\Lambda}} 
\eta_{\mathbb{E}_{ \mathbb{V}\setminus\Lambda}} )
$
is decreasing. Because of the definition of $D_{\Lambda,\Delta}$
we can abuse notation and simply write 
$
\mathds{1}_{D_{\Lambda,\Delta}}
(\omega_{\mathbb{E}_{\Lambda}} 
\eta_{\mathbb{E}_{ \Delta\setminus\Lambda}} ).
$
Given $\epsilon>0$ small enough,
we consider the external magnetic field 
$
\epsilon\pmb{\hat{h}}\equiv(
\epsilon h_{i,p}, \forall \ i\in\mathbb{V}, p=1,\ldots,q)
$.
A straightforward computation shows that 
$\epsilon\pmb{\hat{h}}\prec \pmb{\hat{h}}$, where the partial 
order is given by \eqref{relação de ordem nos campos}. 
If $V$ contains $\Delta$ then it follows from the 
Theorem \ref{teo-monotonicidade-campo-externo-grc}
that 
\[
\phi_{V, \mathrm{ max}}^{\mathrm{ GRC}, \pmb{\hat{h}} }
(
	\mathds{1}_{D_{\Lambda,\Delta}}(
		\cdot\,
		\eta_{\mathbb{E}_{ \Delta\setminus\Lambda}} 
	) 
)
\leqslant
\phi_{V, \mathrm{ max}}^{\mathrm{ GRC}, \epsilon\pmb{\hat{h}} }
(
	\mathds{1}_{D_{\Lambda,\Delta}}(
	\cdot\,
	\eta_{\mathbb{E}_{ \Delta\setminus\Lambda}} 
	)
).
\]
By summing over all 
$\eta_{\mathbb{E}_{ \Delta\setminus\Lambda}}$
the above inequality we get that
\[
\phi_{V, \mathrm{ max}}^{\mathrm{ GRC}, \pmb{\hat{h}} }
(
D_{\Lambda,\Delta}
)
\leqslant
\phi_{V, \mathrm{ max}}^{\mathrm{ GRC}, \epsilon\pmb{\hat{h}} }
(
D_{\Lambda,\Delta}
).
\]
Taking $\epsilon\to 0$, $V\uparrow\mathbb{V}$ and using 
the continuity of $\mathds{1}_{D_{\Lambda,\Delta}}$ and the 
Theorem \ref{limite termodinamico} we get from the last 
inequality that 
\begin{align*}
\phi_{\mathrm{ max}}^{\mathrm{ GRC}, \pmb{\hat{h}} }
(D_{\Lambda,\Delta})
\leqslant
\phi_{\mathrm{ max}}^{\mathrm{ GRC}, \pmb{\hat{0}} }
(D_{\Lambda,\Delta}).
\end{align*}
Since 
$
\cap_{\Delta \supset \Lambda} D_{\Lambda,\Delta}
\uparrow
\{N_{\infty}>1\},
$
when $\Lambda\uparrow \mathbb{V}$
it follows from the continuity of the measure
and Theorem III.3 of \cite{BBCK00}
that 
$\phi_{\mathrm{ max}}^{\mathrm{ GRC}, \pmb{\hat{h}} }
(N_{\infty} >1)=0$.
\end{proof}

To state our next theorem, which is the main theorem 
of the next section, we need to introduce the following 
parameters:
	\begin{align*}
		P_\infty(\beta,\pmb{J},\pmb{\hat{h}})
		\equiv
		\sup_{x\in\mathbb{V} }\,
		\sup_{\phi\in\mathscr{G}^{\mathrm{ GRC}}}
		\phi(|C_x|=\infty)
	\end{align*}
	and
	\begin{align*}
		\widetilde{P}_\infty(\beta,\pmb{J},\pmb{\hat{h}})
		\equiv
		\sup_{x\in\mathbb{V} }\,
		\inf_{\phi\in\mathscr{G}^{\mathrm{ GRC}}}
		\phi(|C_x|=\infty),
	\end{align*}
	where $C_{x}$ is the infinite connected component containing 
	the 
	vertex $x$. For the RC model, the parameters 
	$P_\infty$ and $\widetilde{P}_\infty$
	are defined similarly.
	We also define the critical parameter 
	\begin{eqnarray*}
		\beta_c(\pmb{J},\pmb{\hat{h}})
		\equiv
		\inf\{\beta > 0: P_\infty(\beta,\pmb{J},\pmb{\hat{h}})>0\}.
	\end{eqnarray*}
	To lighten the notation we introduce 
	for each $\mathrm{m}\in \{1,\ldots, q\}$ fixed, the event
	\[
	\mathscr{A}_{\geqslant 1,\mathrm{ m}}^\infty
	\!
	\equiv
	\!
	\left\{
	\!
	(\sigma,\omega)\in\Sigma_{q}\times\Omega 
	:
	\!\!
	\begin{array}{c}
	N_\infty(\sigma,\omega)\geqslant 1 \ 
	\text{and all vertices in any infinite} 
	\\
	\text{connected} \  \text{ component satisfies} \ 
	\sigma_x=\text{m}
	\!	
	\end{array}
	\right\}.
	\]	
\section{Uniqueness and phase transition}
Now we are ready to state and prove one of the main theorems
of this paper.
We emphasize that this theorem was inspired by Theorem II.5 in 
\cite{BBCK00}.
\begin{theorem}[Uniqueness and phase transition]
\label{teo-unicidade-gibbs-percolacao}
	Fix $q\in\mathbb{Z}^+$, $\beta\geqslant 0$, 
	a magnetic field
	$
	\pmb{\hat{h}}
	=
	(h_{i,p}\in\mathbb{R}: i\in\mathbb{V}, p=1,\ldots, q)
	$ 
	and $\{q_p: p=1,\ldots,q\}$ 
	satisfying \eqref{qmmaiorque1}.
\begin{itemize}
	\item[(i)] 
	For all $\pmb{J}\geqslant 0$  
	$
	\mathrm{(}J_{ij}
	\geqslant 
	0, \ \forall\,\{i,j\}\in\mathbb{E}\mathrm{)}
	$, 
	there is at most one probability measure $\mu_0$ in
	$\mathscr{G}^{\mathrm{ ES}}_0 \equiv
	\{\nu\in\mathscr{G}^{\mathrm{ ES}}:\nu(N_\infty=0)=1\}$.
	
	\item[(ii)] 
	If $P_\infty(\beta,\pmb{J},\pmb{\hat{h}})=0$, 
	then
	$\vert \mathscr{G}^{\mathrm{ ES}} \vert = 
	\vert \mathscr{G}^{\mathrm{ RC}}\vert =1.$ 
	In particular, if $\beta<\beta_c(\pmb{J},\pmb{\hat{h}}),$
	then
	$\vert \mathscr{G}^{\mathrm{ ES}} \vert = 
	\vert \mathscr{G}^{\mathrm{ RC}}\vert =1.$ 
\end{itemize}
\begin{itemize} 
	\item[(iii)]
	If $\pmb{J}\geqslant 0$ is an uniform coupling constant  
	$
	\mathrm{(}J_{ij}\equiv J\geqslant 0, \ \forall$ 
	$\{i,j\}\in\mathbb{E}\mathrm{)}$,
	then
	$
	P_\infty(\beta,J,\pmb{\hat{h}})
	=
	\sup_{x\in\mathbb{V}}\
	\phi_{\mathrm{ max}}^{\mathrm{ GRC}}
	(|C_x|=\infty)
	$
	\ and 
	$
	\widetilde{P}_\infty(\beta,J,\pmb{\hat{h}})
	=
	\sup_{x\in\mathbb{V}}\
	\phi_{\mathrm{ free}}^{\mathrm{ GRC}}
	(|C_x|=\infty).
	$
	
	\item[(iv)] 
	Let $\pmb{J}\geqslant 0$ 
	$\mathrm{(}J_{ij}\equiv J\geqslant 0, \ \forall$ 
	$\{i,j\}\in\mathbb{E}\mathrm{)}$.
	If $P_\infty(\beta,\pmb{J},\pmb{\hat{h}})>0$ then
	the states $\phi_{\mathrm{ m}}^{\mathrm{ ES}},$ 
	$
	\mathrm{m}
	\in
	\cap_{i\in\mathbb{V}}\mathcal{Q}_{i,\mathrm{ max}}
	(\pmb{\hat{h}})
	$ 
	are
	extremal $\mathrm{ES}$ Gibbs states with 
	$
	\phi_{\mathrm{ m}}^{\mathrm{ ES}}
	(\mathscr{A}_{\geqslant 1,\mathrm{ m}}^\infty)
	=
	1.
	$
	Moreover under the strong assumption 
	$\widetilde{P}_\infty(\beta,\pmb{J},\pmb{\hat{h}})>0$
	we have that  $\vert \mathscr{G}^{\mathrm{ ES}} \vert > 1$.
		
	\item[(v)]
	Let $\pmb{J}\geqslant 0$ 
	$\mathrm{(}J_{ij}\equiv J\geqslant 0, \ \forall$ 
	$\{i,j\}\in\mathbb{E}\mathrm{)}$.
	If $\beta<\beta_c$, then 
	$
	P_\infty(\beta,J,\pmb{\hat{h}})
	=
	\widetilde{P}_\infty(\beta,J,\pmb{\hat{h}})
	=
	0
	$,
	while both $P_\infty(\beta,J,\pmb{\hat{h}}) > 0$
	and $\widetilde{P}_\infty(\beta,J,\pmb{\hat{h}})>0$ 
	whenever $\beta>\beta_c$.
\end{itemize}
\end{theorem}
\begin{proof} The whole proof follows closely reference 
\cite{BBCK00}.
$(i)$ We prove that
	$\mathscr{G}^{\mathrm{ ES}}_{0}=
	\{\phi_{\mathrm{ free}}^{\mathrm{ ES}}\}$,
	In fact, let $\nu\in \mathscr{G}^{\mathrm{ ES}}_{0}$ 
	and
	$\{\Delta_n:n\in\mathbb{N}\}$ be a cofinal sequence of 
	subsets of $\mathbb{V}$. Then the sequence of random sets  
	$\{\Lambda_n: n\in\mathbb{N}\}$ 
	defined by 
	$\Lambda_n(\omega)\equiv\{x\in\Delta_n:
	x\not\leftrightarrow\Delta_n^c \ \text{in} \ \omega\}$ is 
	also
	increasing. Note that the set $\Lambda_n$ is well 
	defined due to the absence
	of infinite connected components.
	By Theorem \ref{limite termodinamico}, given $\epsilon>0$, 
	we can take $\Delta$ 
	big enough so that
	for each function $f$ with support in 
	$(\Delta,\mathbb{B}_0(\Delta))$ we have
\begin{eqnarray}\label{2m1}
	|\phi_{V,\mathrm{ free}}^{\mathrm{ ES}}(f)
	-
	\phi_{\mathrm{ free}}^{\mathrm{ ES}}(f)
	|\leqslant \epsilon, \quad \forall \ V\supset\Delta.
\end{eqnarray}
	On the other hand, we have that
\begin{multline}
	\nu(f)
	=
	\nu(f \mathds{1}_{\{\Lambda_n(\cdot)\not\supset \Delta\}}) + 
	\nu(f \mathds{1}_{\{\Lambda_n(\cdot)\supset \Delta\}})
	\nonumber
	\\
	=
	\nu(f \mathds{1}_{\{\Lambda_n(\cdot)\not\supset \Delta\}}) + 
	\sum_{\overline{\Lambda}_n(\cdot)\supset \Delta}
	\nu(f \mathds{1}_{\{\Lambda_n(\cdot)
	=
	 \overline{\Lambda}_n\}}).
	\nonumber
\end{multline}
	By using the DLR equations and their equivalent version 
	of conditional expectations for the specification 
	$\{\phi^{\mathrm{ ES}}_{\Lambda,\mathbb{B}_0(\Lambda)}\}$,
	we can rewrite the above expression as
\begin{multline*}
	=
	\nu(f \mathds{1}_{\{\Lambda_n(\cdot)\not\supset \Delta\}}) + 
	\sum_{\overline{\Lambda}_n(\cdot)\supset \Delta}
	\nu
	\big(\phi_{\overline{\Lambda}_n,\mathbb{B}_0
	(\overline{\Lambda}_n)}^{\mathrm{ ES}}
	(
	f \mathds{1}_{\{\Lambda_n(\cdot)= \overline{\Lambda}_n\}}| 
	\cdot
	)
	\big)
	\\
	=
	\nu(f \mathds{1}_{\{\Lambda_n(\cdot)\not\supset \Delta\}}) + 
	\sum_{\overline{\Lambda}_n(\cdot)\supset \Delta}
	\nu
	\big(
	\nu
	(
	f \mathds{1}_{\{\Lambda_n(\cdot)= \overline{\Lambda}_n\}}| 
	\mathscr{F}_{\overline{\Lambda}_n^c,\mathbb{B}_0
	(\overline{\Lambda}_n)^c}
	)
	\big).
\end{multline*}
	For each fixed $n$, the random variable 
	$\mathds{1}_{\{\Lambda_n(\cdot)= \overline{\Lambda}_n\}}$ 
	depends 
	only on the states of the sites and edges
	in $(\overline{\Lambda}_n,\mathbb{B}_0(\overline{\Lambda}_n))$, 
	so we have that this random variable is 
	independent of the
	$\sigma$-algebra 
	$\mathscr{F}_{\overline{\Lambda}_n^c,
	\mathbb{B}_0(\overline{\Lambda}_n)^c}$. Hence the latter
	expression can be rewritten as
\begin{multline*}
	=
	\nu(f \mathds{1}_{\{\Lambda_n(\cdot)\not\supset \Delta\}}) + 
	\sum_{\overline{\Lambda}_n(\cdot)\supset \Delta}
	\nu
	\big( 
	\mathds{1}_{\{\Lambda_n(\cdot)= \overline{\Lambda}_n\}}
	\nu
	(f| 
	\mathscr{F}_{\overline{\Lambda}_n^c,\mathbb{B}_0
	(\overline{\Lambda}_n)^c}
	)
	\big)
	\\
	=
	\nu(f \mathds{1}_{\{\Lambda_n(\cdot)\not\supset \Delta\}}) + 
	\sum_{\overline{\Lambda}_n(\cdot)\supset \Delta}
	\nu
	\big(
	\mathds{1}_{\{\Lambda_n(\cdot)= \overline{\Lambda}_n\}}
	\phi_{\overline{\Lambda}_n,
	\mathbb{B}_0(\overline{\Lambda}_n)}^{\mathrm{ ES}}
	(
	f| 
	\cdot
	)
	\big)
	\\
	=
	\nu(f \mathds{1}_{\{\Lambda_n(\cdot)\not\supset \Delta\}}) + 
	\sum_{\overline{\Lambda}_n(\cdot)\supset \Delta}
	\nu
	\big(
	\mathds{1}_{\{\Lambda_n(\cdot)= \overline{\Lambda}_n\}}
	\phi_{\overline{\Lambda}_n,\mathrm{ free}}^{\mathrm{ ES}}
	(f)
	\big),
\end{multline*}
	where in the second equality we have used again the
	equivalent version of conditional expectation
	for specification 
	$\{\phi^{\mathrm{ ES}}_{\Lambda,\mathbb{B}_0(\Lambda)}\}$, 
	and in the last one we use the
	definition of the measure
	$\phi_{\overline{\Lambda}_n,\mathrm{ free}}^{\mathrm{ ES}}$. 
	So we have the identity
\begin{eqnarray}\label{2m2}
	\nu(f)
	=
	\nu(f \mathds{1}_{\{\Lambda_n(\cdot)\not\supset \Delta\}}) + 
	\sum_{\overline{\Lambda}_n(\cdot)\supset \Delta}
	\nu
	\big(
	\mathds{1}_{\{\Lambda_n(\cdot)= \overline{\Lambda}_n\}}
	\phi_{\overline{\Lambda}_n,\mathrm{ free}}^{\mathrm{ ES}}
	(f)
	\big).
\end{eqnarray}
	Combining the identities \eqref{2m1} and \eqref{2m2} we have
\begin{multline}\label{2m3}
	\nu(f \mathds{1}_{\{\Lambda_n(\cdot)\not\supset \Delta\}}) 
	+
	[\phi_{\mathrm{ free}}^{\mathrm{ ES}}(f)-\epsilon]
	\nu(f \mathds{1}_{\{\Lambda_n(\cdot)\supset \Delta\}}) 
	\leqslant 
	\nu(f)
	\\[0,2cm]
	\leqslant 
	\nu(f \mathds{1}_{\{\Lambda_n(\cdot)\not\supset \Delta\}}) 
	+
	[\phi_{\mathrm{ free}}^{\mathrm{ ES}}(f)+\epsilon]
	\nu(f \mathds{1}_{\{\Lambda_n(\cdot)\supset \Delta\}}). 
\end{multline}
	Since the sequence $\{\Lambda_n(\omega):n\in\mathbb{N}\}$ 
	is increasing, the sequence $\{A_n:n\in\mathbb{N}\}$, with 
	$A_n\equiv\{\Lambda_n(\cdot)\supset \Delta\}$,
	is also increasing. Therefore $\mathds{1}_{A_n}\uparrow 1$. 
	Since $f$ is bounded, taking
	$n\uparrow\infty$ in \eqref{2m3} and using the 
	Dominated Convergence Theorem  yields
\[
	|\nu(f)-\phi_{\mathrm{ free}}^{\mathrm{ ES}}(f)|
	\leqslant \epsilon.
\]
	Since $\epsilon>0$ was arbitrary, we conclude the proof 
	of this item.
	\\\\
$(ii)$
	If $P_\infty(\beta,\pmb{J},\pmb{\hat{h}})=0$, 
	then
\begin{eqnarray}\label{sm4}
	\phi(|C_x|=\infty)=0,
	\quad
	\forall \ \phi\in \mathscr{G}^{\mathrm{ RC}}
	\ \text{and} \ x\in\mathbb{V}.
\end{eqnarray}
	By the uniqueness of the infinite connected component
	(Theorem \ref{teo-unicidade-aglo-infinito})  
	and Lemma \ref{continencia-Gibbs}, it 
	follows that property \eqref{sm4} 
	holds for $\phi_{\mathrm{ max}}^{\mathrm{ RC}}$.
	If $\phi_\nu$ denotes the edge-marginal of 
	$\nu\in\mathscr{G}^{\mathrm{ ES}}$,
	by Theorem \ref{principal-marginal} we have
\[
	0
	=
	\phi_{\mathrm{ max}}^{\mathrm{ RC}}(N_\infty>0)
	\geqslant 
	\phi_\nu(N_\infty>0)= \nu(N_\infty>0),
\] 
	which implies that 
	$\nu\in 
	\mathscr{G}^{\mathrm{ ES}}_0
	$.
	Therefore it follows from the item
	$(i)$ that $\nu=\phi^{\mathrm{ ES}}_{\mathrm{ free}}$. 
	That is,
	$\mathscr{G}^{\mathrm{ ES}}=\mathscr{G}^{\mathrm{ ES}}_0=
	\{\phi^{\mathrm{ ES}}_{\mathrm{ free}}\}.$

	On the other hand, 
	if we denote by 
	$
	\mathscr{G}^{\mathrm{ RC}}_{0}
	\equiv
	\{\phi \in\mathscr{G}^{\mathrm{ RC}}:\phi(N_\infty=0)=1\},
	$ 
	we have from \eqref{sm4} and
	Theorem \ref{limite termodinamico}, 
\[
0
=
\phi_{\mathrm{ max}}^{\mathrm{ RC}}(N_\infty>0)\geqslant 
\phi(N_\infty>0), \quad \forall \ \phi
\in 
\mathscr{G}^{\mathrm{ RC}}.
\]
	Thus $\phi\in \mathscr{G}^{\mathrm{ RC}}_0$. 
	By repeating the proof of
	item $(i)$, using Theorem \ref{vec} and the 
	DLR equations \eqref{Gibbs-rcm}, we have that
	$\mathscr{G}^{\mathrm{ RC}}=
	\{\phi^{\mathrm{ RC}}_{\mathrm{ free}}\}.$
	\\\\
	$(iii)$ 	
	Using Item $(iii)$ of Theorem \ref{limite termodinamico} 
	gives
	\[
	P_\infty(\beta,J,\pmb{\hat{h}})
	\leqslant
	\sup_{x\in\mathbb{V}}
	\phi_{\mathrm{ max}}^{\mathrm{ GRC}}
	(|C_x|=\infty)
	\quad \!\!
	\text{and}
	\quad \!\!
	\widetilde{P}_\infty(\beta,J,\pmb{\hat{h}})
	\geqslant
	\sup_{x\in\mathbb{V}}
	\phi_{\mathrm{ free}}^{\mathrm{ GRC}}
	(|C_x|=\infty).
	\]
	To prove that the equality is attained,
	it is enough to show that 
	$
	\phi_{\mathrm{ max}}^{\mathrm{ GRC}}
	\in
	\mathscr{G}^{\mathrm{ GRC}}
	$
	and
	$
	\phi_{\mathrm{ free}}^{\mathrm{ GRC}}
	\in
	\mathscr{G}^{\mathrm{ GRC}}
	$,
	respectively.	
	By using Theorem \ref{teo-unicidade-aglo-infinito}, we have 
	that
	$
	\phi_{\mathrm{ max}}^{\mathrm{ GRC}, \pmb{\hat{h}} }
	(N_{\infty}\leqslant 1)
	=
	\phi_{\mathrm{ free}}^{\mathrm{ GRC}, \pmb{\hat{h}} }
	(N_{\infty}\leqslant 1)
	=
	1
	$, thus we conclude from Lemma \ref{continencia-Gibbs} that
	$
	\phi_{\mathrm{ max}}^{\mathrm{ GRC}}
	\in
	\mathscr{G}^{\mathrm{ GRC}}
	$
	and
	$
	\phi_{\mathrm{ free}}^{\mathrm{ GRC}}
	\in
	\mathscr{G}^{\mathrm{ GRC}}
	$.
\\\\
$(iv)$
Using the same technique employed by \cite{BBCK00}, one can prove 
that 
\begin{eqnarray*}
\phi_{\mathrm{ max}}^{\mathrm{ RC}}
(x \leftrightarrow \infty)
=
\lim_{V\uparrow\mathbb{V}}
\phi_{V,\mathrm{ max}}^{\mathrm{ RC}}
(x \leftrightarrow V^c).
\end{eqnarray*}
As a consequence, we have that, 
for all 
$
\mathrm{m}
\in 
\cap_{i\in\mathbb{V}}\mathcal{Q}_{i,\mathrm{ max}}(\pmb{\hat{h}})
$,
\begin{eqnarray*}
\phi_{\mathrm{ m}}^{\mathrm{ ES}}
(x \leftrightarrow \infty)
=
\lim_{V\uparrow\mathbb{V}}
\phi_{V,\mathrm{ m}}^{\mathrm{ ES}}
(x \leftrightarrow V^c).
\end{eqnarray*}
Combining the two last identities with the trivial fact
\begin{eqnarray*}
\phi_{V,\mathrm{ m}}^{\mathrm{ ES}}
(x \leftrightarrow V^c,\sigma_x=\mathrm{\tilde{m}})
=
\phi_{V,\mathrm{ m}}^{\mathrm{ ES}}
(x \leftrightarrow V^c)
\delta_{\mathrm{m},\mathrm{\tilde{m}}},
\end{eqnarray*}
and taking the thermodynamic limit, 
we have
\begin{eqnarray}\label{sm7}
\phi_{\mathrm{ m}}^{\mathrm{ ES}}
(\sigma_x=\mathrm{\tilde{m}}|x \leftrightarrow \infty)
=
\delta_{\mathrm{m},\mathrm{\tilde{m}}},
\quad
\forall
\mathrm{m}
\in 
\cap_{i\in\mathbb{V}}\mathcal{Q}_{i,\mathrm{ max}}(\pmb{\hat{h}}).
\end{eqnarray}

Now we prove that the state $\phi_{\mathrm{ m}}^{\mathrm{ ES}}$
is extremal whenever 
$
\mathrm{m}
\in 
\cap_{i\in\mathbb{V}}\mathcal{Q}_{i,\mathrm{ max}}(\pmb{\hat{h}})
$.
To this end, let us assume that 
$
\phi_{\mathrm{ m}}^{\mathrm{ ES}}
(
   \mathscr{A}_{\geqslant 1,\mathrm{ m}}^\infty
)
=
1
$,
this will be proved below.
Suppose that $\phi_{\mathrm{ m}}^{\mathrm{ ES}}$ is not extremal -
then there are two Gibbs measures in 
$\mathscr{G}^{\mathrm{ ES}}$ so that 
\begin{eqnarray}\label{sm8}
\phi_{\mathrm{ m}}^{\mathrm{ ES}}
=
t
\phi_1^{\mathrm{ ES}}
+
(1-t)
\phi_2^{\mathrm{ ES}},
\quad \text{with} \ 
\phi_i^{\mathrm{ ES}}
(\mathscr{A}_{\geqslant 1,\mathrm{ m}}^\infty)=1
\
\text{and}
\ t\in(0,1).
\end{eqnarray}
If $\phi_i^{\mathrm{ RC}}$ denotes the RC marginal of
$\phi_i^{\mathrm{ ES}}$, it follows from  
Lemma VIII.1 in \cite{BBCK00} that
$\phi_i^{\mathrm{ RC}}\in\mathscr{G}^{\mathrm{ RC}}$,
$i=1,2$. This implies that
\begin{eqnarray*}
\phi_{\mathrm{ max}}^{\mathrm{ RC}}
=
t
\phi_1^{\mathrm{ RC}}
+
(1-t)
\phi_2^{\mathrm{ RC}},
\quad
t\in(0,1).
\end{eqnarray*}
By stochastic domination one can prove that 
$\phi_{\mathrm{ max}}^{\mathrm{ RC}}$
is an extremal probability measure, so 
$
\phi_1^{\mathrm{ RC}}
=
\phi_2^{\mathrm{ RC}}
=
\phi_{\mathrm{ max}}^{\mathrm{ RC}}
$.
Using Lemma VIII.3 in \cite{BBCK00}, this fact implies 
$\phi_1^{\mathrm{ ES}}=\phi_2^{\mathrm{ ES}}$, hence
the extremality of $\phi_{\mathrm{ m}}^{\mathrm{ ES}}$ is proved.

Finally, we prove that 
$
\phi_{\mathrm{ m}}^{\mathrm{ ES}}(
\mathscr{A}_{\geqslant 1,\mathrm{ m}}^\infty
)
=
1
$.
Since $P_\infty(\beta,\pmb{J},\pmb{\hat{h}})>0$, 
we get from Item $(iii)$  of Theorem \ref{limite termodinamico} 
that
\begin{align}\label{sm10}
0
<
P_\infty(\beta,\pmb{J},\pmb{\hat{h}})
\leqslant 
\sup_{x\in\mathbb{V}}
\phi_{\mathrm{ max}}^{\mathrm{ RC}}(x\leftrightarrow\infty)
\leqslant
\phi_{\mathrm{ max}}^{\mathrm{ RC}}
(
N_{\infty}\geqslant 1
).
\end{align}
Since $\phi_{\mathrm{ max}}^{\mathrm{ RC}}$ is an extremal 
Gibbs state and $\{N_{\infty}\geqslant 1\}$ is a tail event it 
follows
from the inequality \eqref{sm10} and the uniqueness of the 
infinite connected component 
(Theorem \ref{teo-unicidade-aglo-infinito}) 
that
\[
1=
\phi_{\mathrm{ max}}^{\mathrm{ RC}}
(
N_{\infty}= 1
)
=
\phi_{\mathrm{ m}}^{\mathrm{ ES}}
(N_{\infty}= 1),
\quad \forall  
\mathrm{m}
\in 
\cap_{i\in\mathbb{V}}\mathcal{Q}_{i,\mathrm{ max}}(\pmb{\hat{h}}).
\]
The previous equation
together with the identity \eqref{sm7}
implies, for each
$
\mathrm{m}
\in 
\cap_{i\in\mathbb{V}}\mathcal{Q}_{i,\mathrm{ max}}(\pmb{\hat{h}})
$, 
that 
$
\phi_{\mathrm{ m}}^{\mathrm{ ES}}
(\mathscr{A}_{\geqslant 1,\mathrm{ m}}^\infty)
=1
$.

We now prove the second statement of Item $(iv)$. 
As long as the set 
$
\cap_{i\in\mathbb{V}}\mathcal{Q}_{i,\mathrm{ max}}
(\pmb{\hat{h}})
$
has more than one element, the result follows from the first 
statement 
of the Item $(iv)$. Otherwise, without loss of generality, 
we can assume that  
$
\cap_{i\in\mathbb{V}}\mathcal{Q}_{i,\mathrm{ max}}(\pmb{\hat{h}})
=
\{1\}
$.
Let $\phi^{\tiny\mathrm{RC}} \in \mathscr{G}^{\mathrm{RC}}$ 
be a spin-marginal of $\phi_{\mathrm{ 2}}^{\mathrm{ES}}$, then 
\begin{multline*}
0
<
\widetilde{P}_\infty(\beta,\pmb{J},\pmb{\hat{h}})
\le 
\sup_{x\in\mathbb{V}}
\phi^{\mathrm{ RC}}(x\leftrightarrow\infty)
=
\sup_{x\in\mathbb{V}}
\phi_{\mathrm{ 2}}^{\mathrm{ ES}}(x\leftrightarrow\infty)
\\[0,3cm]
=
\sup_{x\in\mathbb{V}}
\phi_{\mathrm{ 2}}^{\mathrm{ ES}}
(
x \leftrightarrow \infty,
\sigma_x=2
)
\le 
\phi_{\mathrm{ 2}}^{\mathrm{ ES}}
(\mathscr{A}_{\geqslant 1,\mathrm{ 2}}^\infty).
\end{multline*}
Since $\phi_{\mathrm{ 1}}^{\mathrm{ ES}}
(\mathscr{A}_{\geqslant 1,\mathrm{ 1}}^\infty)=1$
and 
$
\mathscr{A}_{\geqslant 1,\mathrm{ 1}}^\infty
\cap 
\mathscr{A}_{\geqslant 2,\mathrm{ 1}}^\infty
=
\emptyset
$, 
it follows from the above inequality that 
$
\phi_{\mathrm{ 1}}^{\mathrm{ ES}}
\neq 
\phi_{\mathrm{ 2}}^{\mathrm{ ES}}.
$
\\\\
$(v)$  By Item $(iii)$ and 
Corollary \ref{coro ge0}, we get that the maps 
$J\mapsto P_\infty(\beta,J,\pmb{\hat{h}})$ 
and $J\mapsto\widetilde{P}_\infty(\beta,J,\pmb{\hat{h}})$ are 
increasing, and so are the maps
$\pmb{\hat{h}}\mapsto P_\infty(\beta,J,\pmb{\hat{h}})$ 
and 
$\pmb{\hat{h}}\mapsto\widetilde{P}_\infty(\beta,J,\pmb{\hat{h}})$, 
with respect to the partial 
order \eqref{relação de ordem nos campos}.  
From the definition, one has 
\begin{eqnarray}\label{ar1}
\widetilde{P}_\infty(\beta,J,\pmb{\hat{h}})
\leqslant
P_\infty(\beta,J,\pmb{\hat{h}}),
\quad
\forall \beta,J\ \text{and}\ \pmb{\hat{h}}.
\end{eqnarray}
From Item $(iii)$ and Theorem
\ref{teorema ge}, we get 
$
P_\infty(\beta,J_1,\pmb{\hat{h}})
\leqslant
\widetilde{P}_\infty(\beta,J_2,\pmb{\hat{h}}) 
$
for all $J_1<J_2$.
By Item $(iii)$, we have that $P_\infty$
and $\widetilde{P}_\infty$ are thermodynamical limits.
Using the form of the Hamiltonian of this model and
the monotonicity properties proved above, we 
get, for all $\beta_1<\beta_2$, that 
\begin{align}\label{ar2}
P_\infty(\beta_1,J,\pmb{\hat{h}})
=
P_\infty(1,\beta_1 J, \beta_1 \pmb{\hat{h}})
&\leqslant
\widetilde{P}_\infty(1,\beta_2 J,\beta_1\pmb{\hat{h}}) 
\nonumber
\\
&\leqslant
\widetilde{P}_\infty(1,\beta_2 J,\beta_2\pmb{\hat{h}})
\nonumber
\\
&=
\widetilde{P}_\infty(\beta_2,J,\pmb{\hat{h}}).
\end{align}
Combining inequalities \eqref{ar1} and \eqref{ar2} yields $(v)$.
\end{proof}

We now consider the $q$-state Potts model where each value 
of the spin is coupled to a distinct and site dependent
external field. The formal Hamiltonian of the model is
\begin{equation}\label{hamiltoniano-potts-BBCK-generalizado}
H(\hat{\sigma})
=
-J\sum_{ \{i,j\} } 
\delta_{\hat{\sigma}_i}\delta_{\hat{\sigma}_j}
-\sum_{p=1}^q \sum_{i} {h_{i,p}\over q}
\delta_{\hat{\sigma}_{i,p}}.
\end{equation}
Let $\mathscr{G}^{\mathrm{{Spin}}}$ denote the set 
of all spin Gibbs states, defined by means of the DLR condition
and the above Hamiltonian (appropriately modified to incorporate
boundary conditions).

\begin{theorem}\label{teo-isomorfismo-GES-GISING}
	Let $\Pi_S:\mathscr{G}^{\mathrm{ ES}}\rightarrow 
	\mathscr{G}^{\mathrm{{Spin}}}$
	denote the mapping that assigns the spin-marginal to infinite
	volume $\mathrm{ES}$ measure. Then $\Pi_S$ is a linear 
	isomorphism.
\end{theorem}
\begin{proof}
	A direct proof of this theorem can be found in \cite{BBCK00}.
\end{proof}
\section{Application - Ising model with power law decay external 
field}  
In this section we apply the results above obtained to prove
the uniqueness of the Gibbs measures, at any positive 
temperature,
for the Ising model in 
$\mathbb{L}\equiv (\mathbb{Z}^d,\mathbb{E}^d)$, 
where $\mathbb{E}^d$ is the set of the nearest neighbors in the 
$d$-dimensional hypercubic lattice, with the Halmiltonian given 
by 
\begin{equation}
\label{hamiltoniano-BCCP}
\mathscr{H}^{\mu,\mathrm{{Ising}}}_{\pmb{h},V}(\sigma)
\equiv 
- \sum_{
	\substack{ i,j\in V  \\ \{i,j\}\in \mathbb{E} }
}
J\, \sigma_i \sigma_j
- 
\sum_{i\in V} h_i\, \sigma_i
-
\!\!\!\! 
\sum_{
	\substack{ i\in V, \ j\in \partial V \\ \{i,j\}\in \mathbb{E} }
}
\!\!\!\!\!\!
J\, \sigma_i \mu_j,
\end{equation}
where $\alpha\geqslant 0$ and $h^*>0$ and 
\[
h_i
=
\begin{cases}
\frac{h^*}{\|i\|^{\alpha}},& \text{if}\ i\neq 0;
\\[0.3cm]
h^*,& \text{if}\ i= 0.
\end{cases}
\]

From now on, we write $\mathscr{G}^{\mathrm{{Spin}}}_{\beta}$
instead of $\mathscr{G}^{\mathrm{{Spin}}}$ to make clear 
its dependence on the inverse temperature. 

By Proposition \ref{Medeq}, 
it follows that the set 
$\mathscr{G}^{\mathrm{{Spin}}}_{\beta}$  
(defined in the last section) 
is precisely 
$
\mathscr{G}^{\mathrm{{Potts}}}_{2\beta}(\pmb{J},\pmb{h})
=
\mathscr{G}^{\mathrm{{Ising}}}_{\beta}(\pmb{J},\pmb{h})
$,
the set of the Gibbs measures of the above Ising model,
if we take in the Hamiltonian 
\eqref{hamiltoniano-potts-BBCK-generalizado}
$q=2$ and the magnetic field given by
\[
\pmb{\hat{h}} 
=
\big( (h_{i,1},h_{i,2}) \in \mathbb{R}^2 : 
\  h^{*}/\|i\|^{\alpha} = h_{i,1}=-h_{i,2},\  
\forall i\in \mathbb{V}\ 
\big).
\]

In order to apply the previous results to study the uniqueness 
of this Ising model with magnetic field decaying to zero with 
polynomial rate $0\leqslant \alpha <1$, we will consider in this 
section the GRC model defined in 
\eqref{definição de general boundary condition}
with $q=2$, the constants $q_p\equiv 1$ and the magnetic field 
$\pmb{\hat{h}}$ as above.

In \cite{BCCP15} the authors proved that for any 
$\alpha \in [0,1)$ there is a positive inverse 
temperature $\beta_{\alpha}<+\infty$ so that, for any $\beta>0$ 
such that $\beta_{\alpha}<\beta$, the set of the Gibbs measures
for the Ising model defined by \eqref{hamiltoniano-BCCP}
is a singleton. By the Dobrushin Uniqueness Theorem, 
we know that for any $\beta< 1/(2d J)$ the set of 
Gibbs measures for this Ising model at these inverse temperatures
is also singleton. In the reference \cite{BCCP15} it was 
conjectured that the set of the Gibbs measures for this model 
with $\alpha \in [0,1)$ is a singleton for any $\beta>0$.
In this work we settle this conjecture.
\begin{center} 
	\begin{minipage}{\linewidth}
	\makebox[\linewidth]{%
	\includegraphics[scale=0.8,keepaspectratio=true]{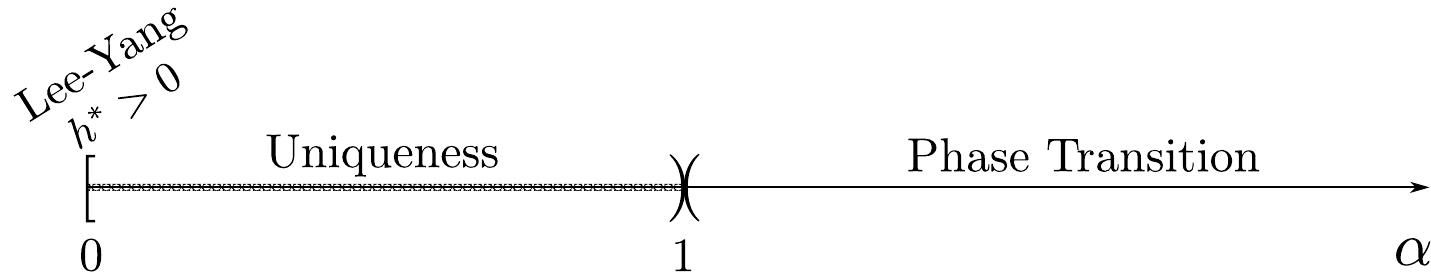}
	}
	\captionof{figure}
	{Uniqueness and non-uniqueness interval for the ferromagnetic 
	Ising model with magnetic field $h_i = h^*/\|i\|^{\alpha}$.}
	\label{cfl}
	\end{minipage}
\end{center}

Suppose that 
\begin{eqnarray*}
\beta_c(\pmb{J},\pmb{\hat{h}})
\equiv
\inf\{\beta > 0: P_\infty(\beta,\pmb{J},\pmb{\hat{h}})>0\} 
=
+\infty.
\end{eqnarray*}
In this case, it follows from Item $(ii)$ 
of Theorem \ref{teo-unicidade-gibbs-percolacao} 
that for any $\beta>0$ we have  
$\vert \mathscr{G}^{\mathrm{ ES}}_{\beta} \vert =1$.  
By Theorem \ref{teo-isomorfismo-GES-GISING} 
we get that 
$\vert \mathscr{G}^{\mathrm{{Spin}}}_{\beta}\vert =1$.

Suppose that 
$\beta_c(\pmb{J},\pmb{\hat{h}})<+\infty$. 
By using once more Item $(ii)$ 
of Theorem \ref{teo-unicidade-gibbs-percolacao},
we obtain the uniqueness for 
$\beta < \beta_c(\pmb{J},\pmb{\hat{h}})$, 
that is, $\vert \mathscr{G}^{\mathrm{{Spin}}}_{\beta}\vert =1$
for such values of $\beta$. If 
$
\beta 
> 
\max\{
\beta_\alpha, 
\beta_c(\pmb{J},\pmb{\hat{h}}) 
\}
$
it was proved in \cite{BCCP15} that 
$\vert \mathscr{G}^{\mathrm{{Spin}}}_{\beta}\vert =1$.
We claim that 
$
\widetilde{P}_\infty(\beta,J,\pmb{\hat{h}})
\equiv  0
$
for any $\beta >0$.
Indeed, take $\beta>\beta_{\alpha}$ 
if 
$
\widetilde{P}_\infty(\beta,J,\pmb{\hat{h}})>0
$, so by Item $(iv)$ of Theorem 
\ref{teo-unicidade-gibbs-percolacao}
we have at least two ES Gibbs measures and 
by Theorem \ref{teo-isomorfismo-GES-GISING} two  
Gibbs measures for the Ising model \eqref{hamiltoniano-BCCP}
which contradicts  \cite{BCCP15}.
Therefore 
$\widetilde{P}_\infty(\beta,J,\pmb{\hat{h}})=0$
whenever $\beta>\beta_{\alpha}$. 
Since the mapping 
$\beta \mapsto \widetilde{P}_\infty(\beta,J,\pmb{\hat{h}}) $
is increasing, the claim follows.

From Item $(v)$ we have for any 
$\beta>\beta_c(\pmb{J},\pmb{\hat{h}})$ that 
$
P_\infty(\beta,J,\pmb{\hat{h}})
>0
$
and 
$
\widetilde{P}_\infty(\beta,J,\pmb{\hat{h}})
> 0,
$
but this contradicts the above claim. Therefore 
we have proved that for any $\alpha \in [0,1)$
that
$
\beta_c(\pmb{J},\pmb{\hat{h}})=+\infty,
$
which implies by Theorem \ref{teo-unicidade-gibbs-percolacao}
that 
$\vert \mathscr{G}^{\mathrm{{Spin}}}_{\beta}\vert =1$
for any $\beta>0$.

\paragraph{Acknowledgements.}
The authors thank Aernout van Enter for fruitful discussions.
We also acknowledge the anonymous referees for their helpful
comments, suggestions and references provided in their reports.
Leandro Cioletti is supported by FEMAT and Roberto Vila is 
supported by CNPq.

\bibliographystyle{alpha}

\vspace{2cm}
\raggedleft
{ \sc 
  Departamento de Matemática
  \smallskip
  
  Universidade de Brasília
  \smallskip

  Brasília, Brazil
  \smallskip
}
  
\texttt{cioletti@mat.unb.br} 
\\
\texttt{rovig161@mat.unb.br} 
 
\end{document}